\documentclass[12pt]{report}

\usepackage{amsthm}
\usepackage{amsmath}
\usepackage[frenchb]{babel}
\usepackage[autolanguage]{numprint}
\usepackage[T1]{fontenc}

\def\N{\mbox{I\hspace{-.15em}N}}

\def\R{\mbox{l\hspace{-.47em}R}}

\newcommand{\om}{(\Omega,{\cal F},P)}
\newcommand{\va}{variable al\'{e}atoire }
\newcommand{\vas}{variables al\'{e}atoires }

\newcommand{\smm}{(M,\Sigma,\mu)}
\newcommand{\E} { {\cal{E}} \!\mbox{xp}}
\newcommand{\Een}{\Bigl(E,{\cal E}_n,\nu \Bigl) }
\newcommand{\Ee}{\Bigl(E,{\cal E},\nu \Bigl) }

\newcommand{\Ex}{{\cal E}xp\Bigl( X\Bigl) }
\newcommand{\Exk}{{\cal E}xp\Bigl( (X_ k)_{k\in \N^*}\Bigl) }

\newcommand{\omr}{(\R,{\cal B}(\R),P)}
\newcommand{\fd}{F(X_o)+\Delta }
\newcommand{\ev}{espace vectoriel }
\newcommand{\evs}{espaces vectoriels }
\newcommand{\sev}{sous-espace vectoriel }
\newcommand{\sevs}{sous-espaces vectoriels }
\newcommand{\F} { {\cal{F}}}
\newtheorem{de}{\Large{Définition}}[section]
\newtheorem{pp}{\Large{Proposition}}[section]

\newtheorem{lm}{\Large{Lemme}}[section]
\newtheorem{co}{\Large{Corollaire}}[section]
\newtheorem{thm}{\Large{Théorème}}[section]
\newcommand{\dps}{\displaystyle }
\newcommand{\n}{\parallel }
\title{Construction et \'etude d'une int\'egrale stochastique}
\author{Ludovic VALET}
\begin{document}
\maketitle
\tableofcontents
\newpage
\part{INDÉPENDANCE DE VARIABLES ALÉATOIRES}
\chapter{Système de \vas indépendantes}
\section{Introduction}
On va construire, pour différentes familles de mesures,
des suites de \vas indépendantes.\\
Pour cela on va transporter des suites de Rademacher généralisée sur
notre espace mesuré via les fonctions de répartitions des mesures
considérées. Dans un premier temps on regardera le cas de mesures
discrètes et dans un deuxième temps on s'intéressera à  des mesures
diffuses pour ensuite généraliser au cas de mesures avec un seul saut.
\vspace{1,54cm}
\begin{de}
On dit qu'une suite, $(b_k)_k$, de \vas est indépendante si, pour tous
entiers naturels
$k_1<\cdots<k_N$, les \vas $b_{k_1},\ldots, b_{k_N}$ sont
mutuellement, ou globalement, indépendantes.
\end{de}
\section{Mesures discrètes}
Soient :
\begin{itemize}
\item
$a_1,\ldots,a_n$ des éléments distincts d'un certain ensemble E,
\item
$\Bigl(b_k\Bigl)_{k=1,\ldots,N}$ : $\Bigl( \Omega=\{ a_1,\ldots,a_n \} ,
{\cal{F}} = {\cal{P}}(\Omega)\Bigl)\to \Bigl(\R ,{\cal{B}}(\R) \Bigl)$ des
applications mesurables.
\end{itemize}
Il existe des partitions $\dps\bigsqcup_{j=1}^{n_k}A_j^k = \Omega$ telle que :
$$b_k = \sum_{j=1}^{n_k} x_j^k 1_{A_j^k}.$$
Posons  $\dps I_j^k = \{l\  | a_l\in A_j^k\}$.
\begin{pp}
Supposons :
\begin{enumerate}
\item
$\forall (j_1,\ldots,j_N)\in\{1,\ldots,n_1\}\times\cdots\times
\{1,\ldots,n_N\}$ :
$$\dps\bigcap_{k=1}^N I_{j_k}^k \not= \emptyset$$
\hspace{-0,7cm} et soit $(p_1,\ldots,p_n)\in]0,1[^n$ vérifiant :
\item
$\forall (j_1,\ldots,j_N)\in\{1,\ldots,n_1\}\times\cdots\times
\{1,\ldots,n_N\}$ :
$$\dps\prod_{k=1}^N \sum_{l\in I_{j_k}^k} p_l =   \sum_{\tiny{  l\in
\bigcap_{k=1}^N I_{j_k}^k} }p_l.$$
\end{enumerate}
Alors  les \vas $X_1,\ldots,X_N$, sont mutuellement indépendantes sur
$L^2\Bigl(\Omega,{\cal{F}},P\Bigl)$, avec $P=\dps\sum_{k=1}^np_k
\delta_{a_k}$.
\end{pp}
\subsection*{Remarque}
Cela oblige $2^{N-1}\leq n $ et fournit une caractérisation de l'indépendance
globale.
\subsection*{\underline{\textbf{Démonstration}}}
\begin{enumerate}
\item[]
\small
Soit
$\Gamma_1\times\cdots\times\Gamma_N\subset
\{x_1^1,\ldots,x_{n_1}^1\}\times\cdots\times\{x_1^N,\ldots,x_{n_N}^N\}.$
$$P
\Bigl(
 \bigl[
 ( X_1,\ldots,X_N)\in\Gamma_1\times\cdots\times\Gamma_N
 \bigl]
\Bigl ) =
P\Bigl(X_1^{-1}(\Gamma_1)\bigcap\ldots\bigcap X_N^{-1}(\Gamma_N)\Bigl)$$
$$= P\Bigl(B_1\bigcap\cdots\bigcap B_N\Bigl)$$
où
$ \dps B_{\nu} =
\bigcup_{k_{\nu}\in J_{\nu}} A_{k_{\nu}}^{\nu}$ avec
$J_{\nu}\subset\{1,\ldots,n_{\nu} \}$,
$ J_{\nu}  = \{l|x_l^{\nu}\in\Gamma_{\nu}  \}$.
$$\begin{array}{lcl}
P\Bigl(B_1\bigcap\cdots\bigcap B_N\Bigl)
 & = &
\dps\sum_{
{(j_1,\ldots,j_N)}\atop{\in J_1\times\cdots\times J_N}  } P\Bigl(
A_{j_1}^1\bigcap\cdots\bigcap A_{j_n}^N\Bigl)\\
 &=&
\dps \sum_{  {(j_1,\ldots,j_N)}\atop{\in J_1\times\cdots\times J_N}  }
\dps\sum_{l\in\bigcap_{k=1}^N I_{j_k}^k  } p_l\\
& =&
\dps \sum_{  {(j_1,\ldots,j_N)}\atop{\in J_1\times\cdots\times J_N}  }
\prod_{k=1}^N\sum_{l\in I_{j_k}^k } p_l\\
&=&
\dps \sum_{  {(j_1,\ldots,j_N)}\atop{\in J_1\times\cdots\times J_N}  }
\prod_{k=1}^NP\Bigl( A_{j_k}^k \Bigl)\\
&=&
 \dps\prod_{k=1}^NP\Biggl(\bigcup_{j_k\in J_k} A_{j_k}^k \Biggl)\\
& = &
\dps\prod_{k=1}^N  P\Bigl(B_k\Bigl).
\end{array}$$
\end{enumerate} \normalsize$\triangle$
\subsection*{Remarque}
En fait il suffit de regarder ce qu'il se passe sur les atomes de la tribu
engendrée
par
$b_1,\ldots,b_N$, puisque :
$$P\Bigl( A_{k_1}^1\bigcap\cdots\bigcap A_{k_N}^N \Bigl) =
\dps\sum_{ l\in\bigcap_{j=1}^NI_{k_j}^j }p_l
=\prod_{i=1}^N\sum_{ l\in I_{k_j}^j } p_l
= \prod_{j=1}^NP\Bigl(A_{k_j}^j\Bigl).$$
\subsection*{Exemple}
\begin{itemize}
\item
$\Omega=\{0,1\}^{2^N}$
\item
$b_k =\Bigl( b_k(i)\Bigl)_{i=1\ldots,2^N}$ $ k=1\ldots,N$ avec : \\
$ b_k(i)=\left\{\begin{array}{ccc}
   1 & \mbox{si} & 2j2^{N-k}\leq i\leq(2j+1)2^{N-k}\\
     &           &   j=0,\ldots,2^k-1\\
    0 & \mbox{si} & (2j+1)2^{N-k}\leq i\leq(2j+2)2^{N-k}\\
      &           &   j=0,\ldots,2^k-1
\end{array}\right.$
\item
$b_k:\Omega    \to  \{-1,1\}$ $\left\{\begin{array}{lcl}
 1     &  \to & 1\\
 0     &   \to &  -1
  \end{array}\right.$
\item
$\dps p_k=\frac{ 1 }{2^N}$
\end{itemize}
$b_0 = 1$\\
$b_1=\Bigl(
\underbrace{1,\ldots,1}_{2^{N-1}},\underbrace{0,\ldots,0}_{2^{N-1}}
\Bigl)$\\
$b_2 =
\Bigl(\underbrace{1,\ldots,1}_{2^{N-2}}\underbrace{0,\ldots,0}_{2^{N-2}}
\underbrace{1,\ldots,1}_{2^{N-2}}\underbrace{0,\ldots,0}_{2^{N-2}}\Bigl)$\\
$A_{k_0}^0\bigcap A_{k_1}^1\bigcap\cdots\bigcap A_{k_N}^N =
A_{k_1}^1\bigcap\cdots\bigcap A_{k_N}^N =a_k$ $\Rightarrow
 P\Bigl( A_{k_1}^1\bigcap\cdots\bigcap A_{k_N}^N\Bigl) = \dps\frac{ 1
}{2^N}$\\
et $P\Bigl( A_{k_1}^1\Bigl)\cdots P\Bigl( A_{k_N}^N\Bigl) =
\dps\frac{ 2^{N-1} }{2^N}\times\cdots\times\frac{ 2^{N-1} }{2^N}=\frac{ 1
}{ 2^N}$
\section{Mesures diffuses}
Dans toute la suite on va considérer un espace mesuré $\om$ où $\Omega =
[0,1]$, ${\cal{F}=B}([0,1])$ et où la mesure P sera précisée dans chaque cas. On
appellera $F$ la fonction de répartition de la mesure P.
\subsection{ P=$\lambda$, mesure de Lebesgue}
$F(x) = x$
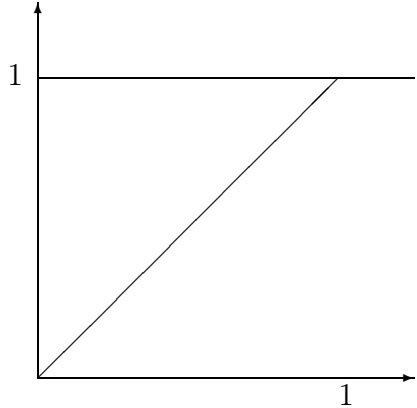
\begin{center}
\unitlength=1cm
\begin{picture}(5,5)
\put(0,0){  \vector(1,0){5}  }
\put(0,0){  \vector(0,1){5}  }
\put(0,0){  \line(1,1){4,0}  }

\put(0,4){  \line(1,0){ 5 }  }
\put(-.4,3.9){  1   }

\put(4,-.35){ 1 }
\end{picture}
\end{center}
Posons
$$\left\{\begin{array}{lcl}
                    \widetilde{r}_0 &  =  &1\\
                     \widetilde{r}_k &  =  & \dps\sum_{j=0}^{2^k-1}(-1)^j
                                 1_{]\frac{ j }{2^k} ,\frac{j+1}{2^k}]}\ \
                                 k\geq1
                \end{array}\right.$$
C'est le système de Rademacher, ce sont des \vas globalement indépendantes
sur
$\om$. On peut le généraliser de la façon suivante :\\
On pose $r_0=1$, après avoir choisi $\alpha_1\in]0,1[$, on prend
$$r_1=1_{]0,\alpha_1]}-1_{]\alpha_1,1]}.$$
Ensuite on découpe chacun des intervalles, après avoir préalablement fixé
 un
$\alpha_2\in]0,1[$, proportionnellement à $\alpha_1$ et $1-\alpha_1$ :
$$r_2 =
1_{]0,\alpha_1\alpha_2]}-1_{]\alpha_1\alpha_2,\alpha_1]}+
1_{]\alpha_1,\alpha_1+(1-\alpha_1)\alpha_2]}-1_{]\alpha_1+(1-\alpha_1)
\alpha_2,1]}.$$
On va construire de façon récurrente notre suite :\\
étant donné $ r_k = \dps\sum_{j=0}^{2^k-1}(-1)^j1_{]a^k_j,a^k_{j+1}]}$, on
choisit
$\alpha_{k+1}\in ]0,1[$ et on découpe chacun des intervalles
proportionnellement à
$\alpha_{k+1}$ et $1-\alpha_{k+1}$ :
$$r_{k+1} = \dps\sum_{j=0}^{2^k-1}
\Biggl(
 1_{\bigl]
     a_j^k , a_j^k +\alpha_{k+1}\bigl(a_{j+1}^k-a_j^k\bigl)
    \bigl]}-
 1_{\bigl]
     a_j^k+\alpha_{k+1}\bigl(a_{j+1}^k-a_j^k\bigl) , a_{j+1}^k
    \bigl] }
\Biggl).$$
En posant :
$$\left\{\begin{array}{lcl}
                   a_{2j}^{k+1}  &  =  &a_j^k\\
                   a_{2j}^{k+1}  &  =  & \alpha_{k+1}(a_{j+1}^k-a_j^k)
                   \end{array}\right.$$
on a :
$$ r_{k+1} = \dps\sum_{j=0}^{2^{k+1}-1}(-1)^j1_{]a^{k+1}_j,a^{k+1}_{j+1}]}.$$
Introduisons maintenant, pour N fixé dans $\N$, quelques notations :
\begin{itemize}
\item
$\beta_0^i(N) = \alpha_i$
\item
$\beta_1^i(N) = 1-\alpha_i\ $ pour $i=1,\ldots,N$
\item
$\beta_{\nu}(N) = \beta_{j_1}^1(N)\cdots\beta_{j_N}^N(N)$\\
où $(j_1,\ldots,j_N)$
est le $\nu^{\mbox{ième}}$ terme pour l'ordre lexicographique dans
$\{0,1\}^N$.\\
\hspace*{-1,1cm}  Dans ces conditions :
\item
$a_{\nu}^N-a_{\nu-1}^N =\beta_{\nu}(N)$ (par récurrence).
\end{itemize}
\begin{lm}
Etant donnés des entiers naturels $1\leq k_1<\cdots<k_N$ :\\
$$r_{k_N}\cdots r_{k_1} =\dps\sum_{i_{N-1}=0 }^{ 2^{k_1}-1  } (-1)^{ i_{N-1} }
\sum_{i_{N-2}=i_{N-1}d_2 }^{ (i_{N-1}+1)d_2-1 } (-1)^{ i_{N-2} }
\cdots
\sum_{i_{0}=i_{1}d_N }^{
(i_{1}+1)d_N-1}(-1)^{i_{0}}
1_{\bigl]a_{i_0}^{k_N} , a_{i_0+1}^{k_N}\bigl]}$$
où $d_j = 2^{k_j-k_{j-1}}$.
\end{lm}
\subsection*{\underline{\textbf{Démonstration}}}
\begin{enumerate}
\item[]
\small
Par récurrence \\
$k_1<k_2$\\
$$ r_{k_1}r_{k_2} =\sum_{j_1=0}^{2^{k_1}-1}\sum_{j_2=0}^{2^{k_2}-1}
(-1)^{j_1+j_2} 1_{  \bigl] a_{j_1}^{k_1},a_{j_1+1}^{k_1}\bigl]\bigcap\bigl]
a_{j_2}^{k_2},a_{j_2+1}^{k_2}\bigl] }$$
$$1_{  \bigl] a_{j_1}^{k_1},a_{j_1+1}^{k_1}\bigl]\bigcap\bigl]
a_{j_2}^{k_2},a_{j_2+1}^{k_2}\bigl] } =
\left\{\begin{array}{ll}
1_{  \bigl] a_{j_2}^{k_2},a_{j_2+1}^{k_2}\bigl] }  \mbox{ pour } &
     j_2\in\{j_1d_2,j_1d_2+1,\ldots, (j_1+1)d_2-1\}\\
0 &\mbox{ ailleurs}
\end{array}\right.$$
On suppose que c'est vrai jusqu'au rang N :\\
$$r_{k_{N+1}}r_{k_N}\cdots r_{k_1} =
\dps\sum_{j_{N+1}=0}^{2^{k_{N+1}-1} }
(-1)^{j_N}1_{ \bigl] a_{j_{N+1}}^{k_{N+1}},a_{j_{N+1}+1}^{ k_{N+1} }\bigl] }
\sum_{i_{N-1}=0 }^{ 2^{k_1}-1  } (-1)^{ i_{N-1} }
\cdots
          \sum_{i_{0}=i_{1}d_N }^{ (i_{1}+1)d_N-1}(-1)^{i_{0}}
1_{\bigl]a_{i_0}^{k_N} , a_{i_0+1}^{k_N}\bigl]}$$
\vspace{0,5cm}
$$\dps1_{ \bigl] a_{j_{N+1}}^{k_{N+1}},a_{j_{N+1}+1}^{ k_{N+1} }\bigl] }
 r_{k_N}\cdots r_{k_1}
=
\dps\sum_{i_{N-1}=0 }^{ 2^{k_1}-1  } (-1)^{ i_{N-1} }
            \cdots
          \sum_{i_{0}=i_{1}d_N }^{ (i_{1}+1)d_N-1}(-1)^{i_{0}}
1_{\bigl]a_{i_0}^{k_N} , a_{i_0+1}^{k_N}\bigl]
\bigcap\bigl] a_{j_{N+1}}^{k_{N+1}},a_{j_{N+1}+1}^{
k_{N+1} }\bigl]  }$$
\vspace{0,5cm}
$$1_{\bigl]a_{i_0}^{k_N} , a_{i_0+1}^{k_N}\bigl]
\bigcap\bigl] a_{j_{N+1}}^{k_{N+1}},a_{j_{N+1}+1}^{
k_{N+1} }\bigl] } = \left\{\begin{array}{lcl}
1_{ \bigl] a_{j_{N+1}}^{k_{N+1}},a_{j_{N+1}+1}^{ k_{N+1} }\bigl] } & si &
j_{N+1}\in\{ i_0d_{k_{N+1} },\ldots\\
&&\hspace*{1cm}\ldots,(i_0+1)d_{k_{N+1} }-1\}\\
0 & \mbox{ sinon}
\end{array}\right.\vspace{0,5cm}$$
On réinjecte dans la somme précédente et, à une réindexation près, on a le
résultat.\\
\end{enumerate} \normalsize$\triangle$
\begin{pp}
$(r_k)_{k\in\N}$ est une suite indépendante sur $\om$.
\end{pp}
\subsection*{\underline{\textbf{Démonstration}}}
\begin{enumerate}
\item[]
\small
Commençons par fixer quelques notations.
N étant un entier strictement positif, prenons
$(k_1,\ldots,k_N)\in\N^N$
tel que :
$0<k_1<\cdots < k_N$
et
$\bigl(\varepsilon_1,\ldots,\varepsilon_N\bigl)\in\{-1,1\}^N.$\\
On aura également besoin de :
$$\varphi_j(x)=\frac{1}{2}\Bigl[(1-\alpha_j)(1-x)+(1+x)\alpha_j\Bigl]
\;(j=1,\ldots,N).$$
On voit facilement que
$$P\Bigl(\bigl[b_k = \varepsilon_k \bigl]\Bigl) =\varphi_k(\varepsilon_k),$$
puisque
$$P\Bigl(\bigl[b_k = \varepsilon_k \bigl] \Bigl)
=\varphi_k(\varepsilon_k)\dps\sum_{\nu=0}^{2^{k-1}-1}\beta_{\nu}(k-1).$$
D'autre part on peut montrer que :
$$\begin{array}{lcl}
P\Bigl(\bigl[b_{k_1} = \varepsilon_{k_1},b_{k_2} = \varepsilon_{k_2}
\bigl]\Bigl)
& = &\varphi_{k_1}(\varepsilon_{k_1})\dps
\sum_{\nu_1=0}^{2^{k_1-1}-1}\Biggl[\beta_{\nu_1}(k_1-1)\varphi_{k_2}
(\varepsilon_{k_2})
\sum_{\nu_2=0}^{2^{k_2-1}-1} \beta_{\nu_2}(k_2-1) \Biggl]\\
& = &\varphi_{k_1}(\varepsilon_{k_1})\varphi_{k_2}(\varepsilon_{k_2})\\
\end{array}$$
De proche en proche on a :\\
$$P\Bigl(\bigl[ b_{k_1} = \varepsilon_{k_1},\ldots, b_{k_N} = \varepsilon_{k_N}
\bigl]\Bigl)
=$$
$$\varphi_{k_1}(\varepsilon_{k_1})\dps
\sum_{\nu_1=0}^{2^{k_1-1}-1}\Biggl[\beta_{\nu_1}(k_1-1)
\Biggl[\cdots \Biggl[
\varphi_{k_N}(\varepsilon_{k_N})
\sum_{\nu_N=0}^{2^{k_N-1}-1} \beta_{\nu_N}(k_N-1) \Biggl]\cdots  \Biggl]$$
$$=\varphi_{k_1}(\varepsilon_{k_1}) \cdots \varphi_{k_N}(\varepsilon_{k_N}).$$
\end{enumerate} \normalsize$\triangle$
\subsection{ P=$\mu$ est une mesure diffuse sur $\R$}
Notons F la fonction de répartition de $\mu$.
\begin{pp}
$(r_k\circ F)_k$ et $(\widetilde{r_k}\circ F)_k$ sont des suites
indépendantes sur
$\Bigl(\R,{\cal{B}}(\R),\mu\Bigl)$.
\end{pp}
\subsection*{\underline{\textbf{Démonstration}}}
\begin{enumerate}
\item[]
\small
Prenons des entiers naturels $0<k_1<\cdots<k_N$, notons $b_k$ pour $r_k$ ou
$\widetilde{r_k}$, alors :
$$b_k :\Bigl([0,1],{\cal{B}}([0,1]),\lambda\Bigl)\to\{-1,1\}$$
et
$$F : \Bigl(\R, { \cal{B} }(\R) , P \Bigl) \to [0,1].$$
D'une part pour $\varepsilon\in\{-1,1\}$ :
$$ P\Bigl( \bigl[b_k\circ F=\varepsilon\bigl]\Bigl)
= P \Bigl( \{x\in\R|F(x) \in\bigl(b_k\bigl)^{-1}(\varepsilon)\} \Bigl)
= P \Bigl( \bigcup_j F^{-1}\bigl(\widetilde{B}_j\bigl)\Bigl).$$
Où
$$\widetilde{B}_j = \left\{\begin{array}{lll}
                               ]a^k_{2j},a^k_{2j+1}] & \mbox{si} &\varepsilon=1\\
                               ]a^k_{2j+1},a^k_{2j+2}] & \mbox{si} &\varepsilon=-1.
                          \end{array}\right.$$
Les $\widetilde{B}_j$ sont disjoints donc :\\
$$P\Bigl( \bigl[b_k\circ F=\varepsilon\bigl]\Bigl)
=\sum_j P \Bigl (F^{-1}\bigl(\widetilde{B}_j\bigl)\Bigl)$$
$$=\sum_j  P \Bigl( \bigl]F^{-1}(a_j),F^{-1}(a_{j+1})\bigl]\Bigl)$$
où
$$F^{-1}(x) = \inf\{y|F(y)=x\}$$
et donc\\
$$P\Bigl( \bigl[b_k\circ F=\varepsilon\bigl]\Bigl)
=\sum_j \lambda\Bigl(]a_j,a_{j+1}]\Bigl)
=\lambda\Bigl(\bigcup_j]a_j,a_{j+1}]\Bigl)
=\cdots=
\lambda\Bigl( \bigl[b_k\circ F=\varepsilon\bigl]\Bigl).$$
D'autre part pour $(\varepsilon_1,\ldots,\varepsilon_N)  \in\{-1,1\}^N$ :
$$P\Bigl(\bigl[b_{k_1}\circ F=\varepsilon_1,\ldots,b_{k_N}\circ F=\varepsilon_N
\bigl]\Bigl)
=P\Bigl(\{x\in\R| F(x)\in\dps\bigcup_{\nu=1}^pB_{j_\nu}\}\Bigl),$$
où
$B_{j_\nu} = ]a_{j_\nu}^{k_N},a_{{j_\nu}+1}^{k_N}]$, pour certains
 ${j_\nu}$, sont
des intervalles disjoints, donc :
$$P\Bigl(\bigl[b_{k_1}\circ F=\varepsilon_1,\ldots,b_{k_N}\circ
F=\varepsilon_N\bigl]\Bigl)
=\dps\sum_{\nu=1}^p\lambda\Bigl({B}_{ j_{\nu} }\Bigl)
=\lambda\Bigl(\dps\bigcup_{\nu=1}^p{B}_{ j_{\nu} }\Bigl)$$
$$= \lambda\Bigl(
\bigl[ b_{k_1}=\varepsilon_1,\ldots,b_{k_N}=\varepsilon_N\bigl]\Bigl).
$$
\end{enumerate} \normalsize$\triangle$
\section{Mesures avec sauts}
\subsection{ Cas d'un seul saut}
On suppose que F est continue en dehors de $X_0$ (plus précisément sur
$[0,X_0]\bigcup]X_0,1]$)
 où elle fait un saut de hauteur
$\Delta$ avec $F(X_0)+\Delta<1$.
\begin{pp}
Soit $b_k = \dps\sum_{j=0}^{r_k}(-1)^j1_{A_j^k}$, où
$A_j^k=]a_j^k,a_{j+1}^k]$,
l'une des suites  $(r_k)_k$ ou  $(\widetilde{r}_k)_k$.\\
Si pour tous entiers naturels $k_1<\cdots<k_N$ on a :
$$\Bigl(\exists j\in\N\Bigl)\Bigl( ]F(X_0),F(X_0)+\Delta]\equiv I\subset
A_{j}^{k_N}\Bigl).$$
Alors la suite $\bigl(b_k\circ F\bigl)_k$ est indépendante sur $\omr$.
\end{pp}
\subsection*{Remarque}
Les conditions de la proposition excluent implicitement le cas
$]F(X_0),F(X_0)+\Delta]    =     A_{j}^{k_N}$
(car si tel était le cas, au rang $k_N+1$, il existerait un indice
$j_{\nu}$ tel que
$I =A_{j_{\nu}}^{k_N+1}\bigcup A_{j_{\nu}}^{k_N+1}$).
\subsection*{\underline{\textbf{Démonstration}}}
\begin{enumerate}
\item[]
\small
Les calculs sont sensiblement les mêmes que dans la démonstration de la
proposition
1.3.2. On doit simplement distinguer le cas où $I\subset
{B}_{j_0}$ pour un certain
$j_0\in\{0,\ldots,2^{k_N}-1\}$.\vspace{0,5cm}\\
$$F^{-1}\bigl({B}_{j_0}\bigl)
= F^{-1}\bigl(
]a_{j_0},F(X_0)]\bigcup]F(X_0),F(X_0)+\Delta]\bigcup]F(X_0)+\Delta,a_{{j_0}+1}
]\bigl)$$
$$=]F^{-1}(a_{j_0}),X_0]\bigcup]X_0,X_0]\bigcup]X_0,F^{-1}(a_{{j_0}+1})]$$
$$=]F^{-1}(a_{j_0}),F^{-1}(a_{{j_0}+1})].$$
Finalement :
$$\mu\Bigl(F^{-1}\bigl({B}_{j_0}\bigl) \Bigl)
 = \lambda\Bigl({B}_{j_0} \Bigl).$$
\end{enumerate} \normalsize$\triangle$
\subsection*{Exemples}
1)
Considérons $F(x)=\left\{\begin{array}{lll}
              x              &  \mbox{si}  &  0\leq x\ \leq X_0=\frac{1}{2}\\
              x+\frac{1}{4}  &  \mbox{si}  &  \frac{3}{4}\leq x\leq 1
                \end{array}\right.$
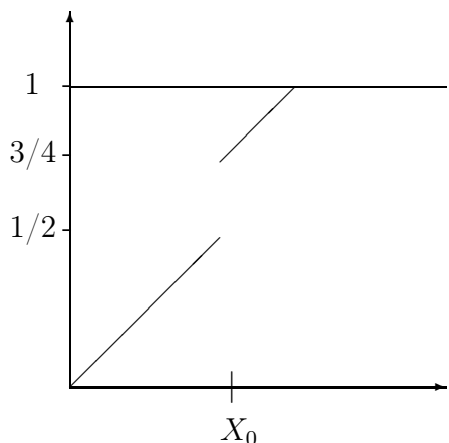
\begin{center}
\unitlength=1cm
\begin{picture}(5,5)
\put(0,0){  \vector(1,0){5}  }
\put(0,0){  \vector(0,1){5}  }
\put(0,0){  \line(1,1){  2 }  }
\put(2,3){  \line(1,1){ 1 }  }
\put(0,4){  \line(1,0){ 5 }  }
\put(-0.1,2){  -  }
\put(-.8,2){  1/2   }
\put(-0.1,3){  -  }
\put(-.8,3){  3/4 }
\put(-0.1,3.915){  -  }
\put(-.6,3.9){  1  }
\put(2.1,-.1){ $|$ }
\put(2,-.7){ $X_0$ }
\end{picture}
\end{center}
$$r_1=1_{]0,\frac{1}{2}]}-1_{]\frac{1}{2},1]}$$
$$r_2=1_{]0,\frac{1}{4}]}-1_{]\frac{1}{4},\frac{1}{2}]}
+1_{]\frac{1}{2},\frac{3}{4}]}-1_{]\frac{3}{4},1]}$$
$$E\bigl(r_1\circ F \cdot r_2\circ F\bigl) = \frac{1}{4}
\not=E\bigl(r_1\circ F\bigl) E\bigl(r_2\circ F\bigl)=-\frac{1}{16}$$
$r_1\circ F$
et $r_2\circ F$ ne sont pas indépendantes. \vspace{0,9cm}\\
\\
2)
Considérons $F(x)=\left\{\begin{array}{lll}
              x              &  \mbox{si}  &  0\leq x\ \leq X_0=\frac{3}{8}\\
              x+\frac{1}{4}  &  \mbox{si}  &  \frac{5}{8}\leq x\leq 1
                \end{array}\right.$
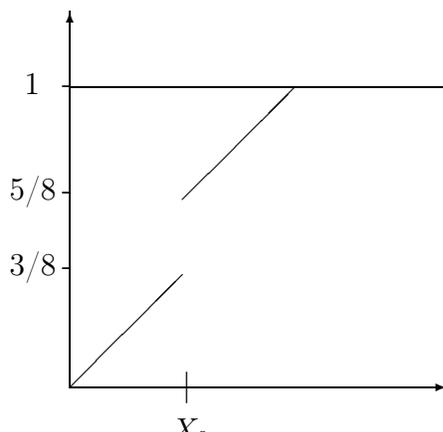
\begin{center}
\unitlength=1cm
\begin{picture}(5,5)
\put(0,0){  \vector(1,0){5}  }
\put(0,0){  \vector(0,1){5}  }
\put(0,0){  \line(1,1){  1.5 }  }
\put(1.5,2.5){  \line(1,1){ 1.5 }  }
\put(0,4){  \line(1,0){ 5 }  }
\put(-.8,1.5){  3/8  }
\put(-.1,1.5){  -  }
\put(-.8,2.5){  5/8 }
\put(-.1,2.5){  -  }
\put(-0.1,3.915){  -  }
\put(-.6,3.9){  1  }
\put(1.5,-.1){ $|$ }
\put(1.35,-.7){ $X_0$ }
\end{picture}
\end{center}
$E\bigl(r_1\circ F \cdot r_2\circ F\bigl) = \frac{1}{4}
\not=E\bigl(r_1\circ F\bigl) E\bigl(r_2\circ F\bigl)=0$, $r_1\circ F$ et
$r_2\circ F$ ne sont pas indépendantes. \\
\subsection*{i) Une première construction}
On va construire une suite de Rademacher généralisée en choisissant les
$\alpha_k$ de façon à découper l'intervalle $I\equiv]F(X_0),F(X_0)+\Delta]$
 après
$F(X_0)+\Delta$ (pour appliquer la proposition 1.4.1.).\vspace{0,5cm}\\
$$b_0=1$$
on choisit $\alpha_1>F(X_0)+\Delta$ dans ]0,1[,\\
$$b_1= 1_{]0,\alpha_1]}-1_{]\alpha_1,1]} =
1_{]0,F(X_0)+\Delta]\bigcup]F(X_0)+\Delta,\alpha_1]}-1_{]\alpha_1,1]}$$
on choisit $\alpha_2\in]0,1[$ tel que $\alpha_1\alpha_2>F(X_0)+\Delta$,\\
$$b_2=1_{]0,F(X_0)+\Delta]\bigcup]F(X_0)+\Delta,\alpha_1\alpha_2]}-
1_{]\alpha_1\alpha_2,\alpha_1]}+1_{]\alpha_1,\alpha_1+(1-\alpha_1)\alpha_2]}
-1_{]\alpha_1+(1-\alpha_1)\alpha_2,1]}.$$
\subsubsection*{CONDITION TECHNIQUE}
\vspace{0,9cm}
\begin{tabular}{||p{15cm}  }
On prend une suite $(\alpha_k)_{k\geq1}$ vérifiant :\\
$\forall N\in\N-\{0,1\} :$
$$(C1)\hspace{01,5cm}\prod_{k=1}^N\alpha_k>F(X_0)+\Delta.$$
\end{tabular}
\vspace{0,9cm}
\\
Supposons, qu'au rang N, l'intervalle contenant I soit de la forme :
$$]0,\alpha_1\cdots \alpha_N] = ]a,b].$$
On le découpe en deux, comme dans le paragraphe 1.3.1, proportionnellement
à $\alpha_{N+1}$ et $1-\alpha_{N+1}$ :
$$]a,a+(a-b)\alpha_{N+1}]\bigcup]a+(a-b)\alpha_{N+1},b].$$
Pour respecter le modèle  de construction imposé il faut :\\
$$a+(a-b)\alpha_{N+1}>F(X_0)+\Delta.$$
En remplaçant a et b par leur valeur, (C1) donne le  résultat, à savoir :
$I\subset]0,\alpha_1\cdots\alpha_{N+1} [.$
\begin{lm}On utilise les même notations\\
1) Si $\forall k\geq1$
$\alpha_k=\alpha\in]0,1[$ alors on ne peut pas construire la suite
$(b_k)_{k\geq0}$.\\
2) La suite dont le terme générique est :
$$ \alpha_k =\left\{\begin{array}{ll}
    F(X_0)+\Delta+a & k=1\\
    1-\frac{ a^{k-1}(1-a) }{F(X_0)+\Delta+a^{k-1} }& k\geq2
\end{array}\right.$$
avec $a\in]0,1-F(X_0)-\Delta[$
convient à la construction des $b_k$ suivant le modèle imposé.
\end{lm}
\subsection*{\underline{\textbf{Démonstration}}}
\begin{enumerate}
\item[]
\small
On va poser $p(N) = \dps\prod_{i=1}^N\alpha_i$.\\
1)\\
Si $\alpha_k=\alpha$ $\forall k\geq1$ alors $p(N) = \alpha^N$ et
$p=\dps\lim_{N\to\infty} p(N) = 0$. Comme $\fd>0$,
la condition (C1) n'est pas
vérifiée.
\vspace{0,7cm}\\
2)
$$\alpha_k=1-\dps\frac{ a^{k-1}(1-a) }{F(X_0)+\Delta+a^{k-1} }
 = \dps\frac{ \fd +a^k }{ \fd+a^{k-1}  }$$
$$p(N) = (\fd +a )\dps\frac{ \fd+a^2}{\fd+a  }\cdots\dps\frac{
\fd+a^N}{\fd+a^{N-1}  }  =\fd+a^N.$$
p(N) est une suite décroissante vers $\fd$, (C1) est donc vérifiée.
\end{enumerate} \normalsize$\triangle$
\vspace{0,5cm}
\begin{pp}
La suite $\bigl(b_k\circ F\bigl)_k$ est indépendante sur $\om$.
\end{pp}
\subsection*{\underline{\textbf{Démonstration}}}
\begin{enumerate}
\item[]
\small
Elle découle immédiatement de la nature du découpage en utilisant la
proposition 1.4.1.\\
\end{enumerate} \normalsize$\triangle$
\subsection*{Remarque}
L'inconvénient de cette méthode réside dans le fait que l'intervalle
$[0,\fd]$ n'est pas découpé. Même si $\fd$ est "très petit", il existe
un rang pour
lequel $F^{-1}\bigl(\bigl[0,\fd]\bigl)$ est le plus grand intervalle.
\subsection*{ii) Une deuxième construction}
On va construire une suite de Rademacher généralisée en choisissant les
$\alpha_k$
de façon à découper l'intervalle $I\equiv]F(X_0),F(X_0)+\Delta]$ avant
$F(X_0)$ (pour appliquer la proposition 1.3.2.).\\
$$b_0=1$$
on choisit $\alpha_1<F(X_0)$ dans ]0,1[,
$$b_1= 1_{]0,\alpha_1]}-1_{]\alpha_1,1]} =
1_{]0,\alpha_1]}-1_{]\alpha_1,F(X_0)]\bigcup]F(X_0),1]}$$
pour $\alpha_2\in]0,1[$ tel que :
$$\alpha_1+(1-\alpha_1)\alpha_2<F(X_0)$$
$$b_2=1_{]0,\alpha_1\alpha_2]}-
1_{]\alpha_1\alpha_2,\alpha_1]}+1_{]\alpha_1,\alpha_1+(1-\alpha_1)\alpha_2]}-
1_{]\alpha_1+(1-\alpha_1)\alpha_2,F(X_0)]\bigcup]F(X_0),1]}.$$
\subsubsection*{CONDITION TECHNIQUE}
\vspace{0,9cm}
\begin{tabular}{||p{15cm}  }
On prend une suite $(\alpha_k)_{k\geq1}$ vérifiant :\\
$\forall N\in\N-\{0,1\} :$
$$(C2)\hspace{01,5cm}\alpha_1+\dps\sum_{k=2}^N \alpha_k\prod_{ 1\leq j\leq k-1
}(1-\alpha_j)<F(X_0).$$
\end{tabular}
\vspace{0,9cm}\\
Supposons, qu'au rang N, l'intervalle contenant I soit de la forme :
$$\Bigl]
\alpha_1+\dps\sum_{k=2}^N \alpha_k\prod_{ 1\leq j\leq k-1 }(1-\alpha_j),1\Bigl]
=]a,b].$$
On le découpe en deux, avec la méthode du 1.3, proportionnellement à
$\alpha_{N+1}$
et
$1-\alpha_{N+1}$ :
$$]a,a+(a-b)\alpha_{N+1}]\bigcup]a+(a-b)\alpha_{N+1},b].$$
Pour respecter le modèle  de construction imposé il faut :\\
$$a+(a-b)\alpha_{N+1}<F(X_0).$$
En remplaçant a et b par leur valeur, (C2) donne le  résultat, à savoir :
$$I\subset\left]\alpha_1+\dps\sum_{k=2}^N \alpha_k\prod_{ 1\leq j\leq k-1
}(1-\alpha_j)\;,\;1\right].$$
\begin{lm}On conserve le mêmes hypothèses\\
1) Si $\forall k\geq1$ $\alpha_k=\alpha\in]0,1[$ alors on ne peut pas
construire la suite $(b_k)_{k\geq0}$.\\
2) Etant donné $a\in]0,1[$, la suite $\bigl(\alpha_k\bigl)_{k\geq1}$ définie
par
$\alpha_k =a^k$ vérifie la condition (C2) si :
$$0\leq\frac{  a  }{  1-a(a+1)  }<F(X_0).$$
\end{lm}
\subsection*{\underline{\textbf{Démonstration}}}
\begin{enumerate}
\item[]
\small
On va étudier
$$ S=\alpha_1+\dps\sum_{k=2}^{\infty} \alpha_k\prod_{ 1\leq j\leq k-1}(1-\alpha_j).
$$
1) Si $\forall k$ $\alpha_k=\alpha$ alors
$$ S = \alpha\Bigl(\dps\sum_{j=2}^{\infty}(1-\alpha)^{k-1}+1\Bigl) = 1.$$
La condition (C2) s'écrit
$$1<F(X_0)$$
ce qui est faux.\\[1cm]
2) On remplace $\alpha_k$ par $a^k$ :
$$ S = a +  \dps\sum_{j=1}^{\infty}a^{k+1}\prod_{ 1\leq j\leq k-1}(1-a^j).$$
On va utiliser l'identité :
$$ (1-x_1)\cdots(1-x_k) =1+ \sum_{j=1}^k(-1)^j\sum_{1\leq i_1<\cdots<i_j\leq k}
x_{i_1}\cdots x_{i_j}$$
pour majorer le produit.
$$ S=a +  \dps\sum_{j=1}^{\infty}a^{k+1}\Biggl(
\sum_{j=1}^k(-1)^j\sum_{1\leq i_1<\cdots<i_j\leq k}a^{{i_1}+\cdots+ i_j}
+1\Biggl)$$
$$\leq \dps\frac{ a }{ 1-a} + \dps\sum_{j=1}^{\infty}a^{k+1}
\sum_{j=1}^k {k\choose j }a^{\frac{ j(j+1) }{ 2 } },$$
car
$$1+\cdots+j\leq{i_1}+\cdots+ i_j.$$
Comme
$$j\leq\frac{ j(j+1) }{ 2 }$$
et
$$a^{\frac{ j(j+1) }{ 2 }}\leq a^{j}$$
on a :
$$ S\leq \dps\frac{ a }{ 1-a} + \dps\sum_{j=1}^{\infty}a^{k+1}
\sum_{j=1}^k {k\choose j }a^j
=\dps\frac{ a }{ 1-a} + \dps\sum_{j=1}^{\infty}a^{k+1}\Bigl[ (1+a)^k-1\Bigl]
= \frac{ a }{ 1-a(1+a)  }.$$
Si on note S(N) les sommes partielles alors $\Bigl(S(N)\Bigl)_N$ est une
suite
croissante. Par hypothèse $S\leq F(X_0)$ et la condition (C2) est donc
vérifiée.\\
\end{enumerate} \normalsize$\triangle$
\subsection*{Remarque}
Posons $f(a) =\frac{ a }{ 1-a(1+a)  }$ sur $]0,\frac{\sqrt5-1}{2}[$.
C'est une fonction continue croissante de 0 à $+\infty$;
donc
$\{x|f(x)<F(X_0)\}\not=\emptyset$, on choisit $a$ dans cet ensemble.
\begin{pp}
La suite $\bigl(b_k\circ F\bigl)_k$ est indépendante sur $\om$.
\end{pp}
\subsection*{\underline{\textbf{Démonstration}}}
\begin{enumerate}
\item[]
\small
Elle découle immédiatement de la nature du découpage en utilisant la
proposition 1.4.1.\\
\end{enumerate} \normalsize$\triangle$
\subsection*{Remarque}
L'inconvénient de cette méthode réside dans le fait que l'intervalle
$][F(X_0),1]$ n'est pas découpé. Même si $F(X_0)$ est proche de 1, il existe un rang
pour lequel $F^{-1}\Bigl(]F(X_0),1]\Bigl)$ est le plus grand intervalle.
\subsection*{iii) Une troisième construction}
On va construire une suite de Rademacher généralisée en choisissant les
$\alpha_{2k}$  (resp. les $\alpha_{2k+1}$)  de façon à découper l'intervalle
contenant I après $\fd$ (resp. avant $F(X_0)$) :
$$b_0=1$$
on choisit $\alpha_1<F(X_0)$ dans ]0,1[,
$$b_1= 1_{]0,\alpha_1]}-1_{]\alpha_1,1]} =
1_{]0,\alpha_1]}-1_{]\alpha_1,F(X_0)]\bigcup]F(X_0),1]}$$
on choisit $\alpha_2\in]0,1[$ tel que :
$$\alpha_1+(1-\alpha_1)\alpha_2>\fd $$
$$b_2=1_{]0,\alpha_1\alpha_2]}-
1_{]\alpha_1\alpha_2,\alpha_1]}
+1_{]\alpha_1,F(X_0)]\bigcup]F(X_0),\alpha_1+(1-\alpha_1)\alpha_2]}
-1_{]\alpha_1+(1-\alpha_1)\alpha_2,1]}.$$
\subsubsection*{CONDITION TECHNIQUE}
\vspace{0,9cm}
\begin{tabular}{||p{15cm}  }
On prend une suite $(\alpha_k)_{k\geq1}$ vérifiant :\\
$$(C3)\hspace{1,5cm}\left\{\begin{array}{llll}
g(k)\equiv &\dps\sum_{j=0}^{ \frac{k-1}{2} }D(2j+1)  &  <F(X_0)  & \mbox{k
impair}\\
d(k)\equiv &\dps\sum_{j=0}^{\frac{k}{2} -1 }D(2j+1) +U(k)  & >\fd &\mbox{k pair }
\end{array}\right.$$
Où
$$             D(2j+1)
=\dps\prod_{i=0}^{j-1}\bigl(1-\alpha_{2i+1}\bigl)
\Bigl(\prod_{i=1}^j\alpha_{2i}\Bigl)\alpha_{2j+1}$$
$$U(2j) = \dps\prod_{i=0}^{j-1}\bigl(1-\alpha_{2i+1}\bigl)
\Bigl(\prod_{i=1}^j\alpha_{2i}\Bigl).$$
\end{tabular}
\vspace{0,9cm}\\
Supposons qu'au rang 2N (resp. 2N+1), l'intervalle contenant I soit de la
forme :
$$]g\bigl(2(N-1)+1\bigl),d(2N)]$$
(resp. $]g\bigl(2N+1\bigl),d(2N)] $).\vspace{0,5cm}\\
On fait le découpage au rang 2N+1 :\\
$$]g(2N-1),g(2N-1)+\Bigl(d(2N)-g(2N-1)\Bigl)\alpha_{2N+1}]\bigcup$$
$$]g(2N-1)+\Bigl(d(2N)-g(2N-1)\Bigl)\alpha_{2N+1}  ,  d(2N)]$$
$$\Bigl(d(2N)-g(2N-1)\Bigl)\alpha_{2N+1} = U(2N)\alpha_{2N+1} =D(2N+1).$$
Donc
$$g(2N-1)+\Bigl(d(2N)-g(2N-1)\Bigl)\alpha_{2N+1} = g(2N+1),$$
comme
$$g(2N+1)<F(X_0),$$
l'intervalle contenant I est :
$$]g(2N+1),d(2N)]$$
On fait le découpage au rang 2N+2 :
$$]g(2N+1),g(2N+1)+\Bigl(d(2N)-g(2N+1)\Bigl)\alpha_{2N+2}]\bigcup$$
$$]g(2N+1)+(d(2N)-g(2N+1))\alpha_{2N+2}  ,  d(2N)].$$
Avec
$$\Bigl(d(2N)-g(2N+1)\Bigl)\alpha_{2N+2} =
U(2N)\Bigl(1-\alpha_{2N+1}\Bigl)\alpha_{2N} = U\Bigl(2(N+1)\Bigl)$$
et
$$g(2N+1)+U\bigl(2(N+1)\bigl) = d\bigl( 2(N+1) \bigl).$$
Comme
$$d\bigl(2(N+1) \bigl)\geq\fd$$
l'intervalle contenant I est
$$](g(2N+1),d\bigl(2(N+1) \bigl)].$$
Ce qui termine la récurrence pour le cas pair.
On vérifie aisément le résultat pour le cas impair.\\
\subsection*{Remarque}
Pour $k\geq1$ :\\
$$g(2k+1) =
\alpha_1+(1-\alpha_1)\alpha_2\alpha_3+(1-\alpha_1)
\dps\sum_{j=2}^{k}\alpha_{2j}\alpha_{2j+1}
    \prod_{i=1}^{j-1}(1-\alpha_{2i+1})\alpha_{2i}$$
Pour $k\geq2$ :
$$d(2k) = \alpha_1+
\dps\sum_{j=2}^{k-1}\Bigl(\alpha_{2j+1}
    \prod_{i=0}^{j-1}(1-\alpha_{2i+1})
    \prod_{i=1}^{j}\alpha_{2i}\Bigl)
+ \prod_{i=0}^{k-1}(1-\alpha_{2i+1})
    \prod_{i=1}^{k}\alpha_{2i}$$
$$ = g(2k+1)+(1-\alpha_1)\prod_{i=1}^{k}(1-\alpha_{2i+1})\alpha_{2i}.$$
Récapitulons:\\
$$\left\{\begin{array}{lcl}
g(1)  & = & \alpha_1\\
  & \\
k\geq1 &  &\\
g(2k+1) & = &\alpha_1+(1-\alpha_1)\alpha_2\alpha_3+(1-\alpha_1)
\dps\sum_{j=2}^{k}\alpha_{2j}\alpha_{2j+1}
    \prod_{i=1}^{j-1}(1-\alpha_{2i+1})\alpha_{2i}\\
d(2) & = &\alpha_1 +(1-\alpha_1)\alpha_2\\
  & &\\
k\geq2   & &\\
  & \\
d(2k) &  = &
g(2k+1)+(1-\alpha_1)\dps\prod_{i=1}^{k}(1-\alpha_{2i+1})\alpha_{2i}
\end{array}\right.$$
\subsection*{Monotonie}
$\Bigl( g(2k+1) \Bigl)_{k}$ est une suite croissante, majorée donc
convergente.\\
Pour $k\geq2$ :
$$d(2k+2)-d(2k) =
(1-\alpha_1)\Bigl(\prod_{i=1}^{k}(1-\alpha_{2i+1})\alpha_{2i}\Bigl)
(\alpha_{2k+2}-1)\leq0. $$
$\Bigl( d(2k) \Bigl)_{k}$ est une suite décroissante, minorée donc
convergente.\\
Posons
$$\left\{\begin{array}{l}
g =\dps\lim_{k\to\infty} g(2k+1)\\
d =\dps\lim_{k\to\infty} d(2k).
\end{array}\right.$$
\vspace{0,5cm}
\begin{lm}On conserve les mêmes notations\\
1)  Si pour tout entier naturel k, $\alpha_k =\alpha$ alors on ne peut pas
construire la suite
$\Bigl( b_k\Bigl) _{k\geq0}$.\\
2) La suite $\Bigl( \alpha_k\Bigl)_{k\geq1}$ définie par :
$$\left\{\begin{array}{lcl}
\alpha_1 & = &1-A(\widetilde{\Delta}-a)\\
k\geq1 & & \\
\alpha_{2k}  & =&
\dps\frac{  \widetilde{\Delta} +\frac{ a}{2k}      }
         {  \widetilde{\Delta}+\frac{ a}{2k-1}    }
\\
  &    &    \\
\alpha_{2k+1}  & =& 1-
\dps\frac{  \widetilde{\Delta} +\frac{ a}{2k+1}   }
        {   \widetilde{\Delta}+\frac{a}{2k}         }
\end{array}\right..$$
Avec :                       $\bullet$  $A = \dps\frac{1}{\ln2}$\\
\hphantom{Avec :} $\bullet$
$\widetilde{\Delta} = \biggl( \Delta +\dps\frac{1}{10^p} \biggl)\ln2$\\
\hphantom{Avec :} $\bullet$
$a =1-\biggl( \fd+\dps\frac{1}{10^{p+1}}\biggl)$.\\
p étant choisit pour que $\fd+\dps\frac{1}{10^{p+1}}<1$,
vérifie la condition (C3)
\end{lm}
\subsection*{\underline{\textbf{Démonstration}}}
\begin{enumerate}
\item[]
\small
1) Supposons $\forall k\geq1$ $ \alpha_k=\alpha$\\
$$g(2k+1) = \alpha +\dps\sum_{j=1}^k(1-\alpha)^j\alpha^{j+1}$$
$$=\alpha+\alpha^2(1-\alpha)
\dps\frac{ 1-\Bigl( (1-\alpha)\alpha\Bigl)^k }
                 { 1-  (1-\alpha)\alpha                           } .$$
$$ \alpha \in]0,1[\Rightarrow
\dps\lim_{k\to\infty}  g(2k+1)
=\frac{ \alpha }{1- (1-\alpha)\alpha  }.$$
Par des calculs similaires on obtient :
$$\lim_{k\to\infty}  d(2k)
=\frac{ \alpha }{Ê1- (1-\alpha)\alpha  }.$$
A la limite on aurait :
$$\fd\geq \frac{ \alpha }{Ê1- (1-\alpha)\alpha  }\geq F(X_0)$$
ce qui est impossible car $\Delta >0$.\\[1cm]
2)
On va vérifier que $d\geq\fd$ et $g\leq F(X_0)$ de cette façon (C3)  sera
vérifiée.\\
On remplace $\alpha_k$ par sa valeur dans les expressions donnant $g$ et $d$
pour obtenir :
$$g = 1-A\Bigl(  \widetilde{\Delta} +a\ln2\Bigl)
= F(X_0) - \dps\frac{  9  }{  10^{p+1}   }$$
$$d = g+A \widetilde{\Delta}
= \fd + \dps\frac{  1  }{  10^{p+1}   }.$$
\end{enumerate} \normalsize$\triangle$
\subsection*{Remarque}
On a :
$$\Bigl( F(X_0) -g\Bigl)\bigvee\Bigl(  d-(\fd)\Bigl) \leq\frac{  1  }{ 10^{p+1}
}.$$
\vspace{0,5cm}
\begin{pp}
La suite $\bigl(b_k\circ F\bigl)_k$ est indépendante sur $\om$.
\end{pp}
\subsection*{\underline{\textbf{Démonstration}}}
\begin{enumerate}
\item[]
\small
Elle découle immédiatement de la nature du découpage en utilisant
la proposition 1.4.1.\\
\end{enumerate} \normalsize$\triangle$
\subsection*{Remarques}
1) En prenant $\widetilde{\Delta} =\Delta\ln2$\\
\hphantom{1) En prenant }$a =1-\Bigl(  \fd\Bigl).$\\
On obtient                         $ g = 1-\Delta-a = F(X_0)$\\
\hphantom{On obtient }$ d = 1-a = \fd$.\\
Comme la suite$\Bigl(   g(2k+1)\Bigl)_k$ (resp. $\Bigl(  d(2k)\Bigl)_k$)
est décroissante (resp. croissante), la condition (C3) est
toujours vérifiée.\\[1,5cm]
2) Ce découpage est plus intéressant que les deux précédents car seul I
n'est pas découpé et comme
$$F^{-1}(I) =\emptyset$$
quand on fait opérer $F^{-1}$, cela n'apparaît pas.
\subsection{Cas d'un nombre fini de sauts}
%
%
%
%
%
\newpage
\chapter{Système maximum de \vas indépendantes}
\section{Position du problème}
On fixe :  $\bullet$ des éléments $a_1,\ldots, a_n$, distincts, dans un
ensemble $\Omega$\\
\hphantom{On fixe :  }$\bullet$ $\dps\mu= \sum_{k=1}^np_k\delta_{a_k}$ avec
$ \dps\sum_{k=1}^np_k = 1$, $p_k\in$ ]0,1[, $k=1,\ldots,n$\\
\hphantom{On fixe :  }$\bullet$ l'espace probabilisé $\om$.\\
On se pose la question suivante : étant donnée une
\va c=$\dps\sum_{k=1}^nc_k1_{  \{a_k\}  }$ (où $c_k\in\R$), quelle est la taille d'un système maximum de \vas
globalement indépendantes?
\subsection*{Modélisation}
On construit l'isomorphisme naturel :\\
$$L_{\R}^2(\mu)\to\R^n\hspace{1cm}
\dps\sum_{k=1}^nx_k1_{  \{a_k\}  }\mapsto\dps\sum_{k=1}^nc_ke_k$$
où $\Bigl(e_k\Bigl)_{k=1,\ldots,n}$ est la base canonique de $\R^n$.\\
La \va c s'écrit alors c =$\left(\begin{array}{c}
                                 c_1\\
                                 \vdots\\
                                 c_n
                                 \end{array}\right)$.
\section{Caractérisation}
\begin{thm}
Une \va  b =$\left(\begin{array}{c}
                                 b_1\\
                                 \vdots\\
                                 b_n
                                 \end{array}\right)$
 est indépendante de c ssi pour toutes fonctions mesurables f et g  sur $\Bigl(\R,B(\R)\Bigl)$  :
 $$\Bigl( f(b) , Ag(c) \Bigl)_{\R^n} =0$$
 où $A = \Bigl( a_{ij} \Bigl)_{1\leq i,j\leq n}$ avec
 $a_{ij} =\left\{\begin{array}{lcl}
 p_i^2-p_i  & si  &  i=j\\
 p_ip_j      & si  & i\not= j
 \end{array}\right.$
 \end{thm}
\subsection*{\underline{\textbf{Démonstration}}}
\begin{enumerate}
\item[]
\small
b et c sont indépendantes ssi $E\Bigl(f(b)g(c)\Bigl) = E\Bigl(f(b)\Bigl)E\Bigl(g(c)\Bigl)$ pour toutes fonctions f, g mesurables sur $\Bigl(\R,B(\R)\Bigl)$.\\
$f(b) =\dps\sum_{k=1}^nf(b_k)1_{  \{a_k\}  }$, $g(c) =\dps\sum_{k=1}^ng(c_k)1_{  \{a_k\}  }$.\\
En explicitant la première égalité  :
$\dps\sum_{k=1}^n\Bigl(p_k^2-p_k\Bigl)f(b_k)g(c_k) =\sum_{j\not= k}p_jp_kf(b_j)g(c_k).$\\
\end{enumerate} \normalsize$\triangle$
%
%
%
%
%
%
%
%
%
%
\subsection*{Remarques}
\begin{enumerate}
\item
Ce théorème n'est que la traduction sur $\R^n$ via l'isomorphisme, de ce
qu'il se passe dans $L^2(\mu)$.
\item
Le noyau de A est le $\R$-\ev engendré par $1_{\R^n}$ (on vérifie aisément
que $Ax=0$ $\iff$ $\forall i$ $\dps\sum_{k=1}^np_kx_k =x_i$).
\end{enumerate}
\section{Cas général}
\subsection{Notations}
\begin{itemize}
\item
On notera plutôt 1 que $1_{\R^n}$
\item
NCD =nombre de coordonnées distinctes de c,
\item
G=[g(c), g] c'est le $\R$-\ev engendré par les éléments g(c), g parcourant
l'espace des fonctions mesurables sur $\Bigl(\R,B(\R)\Bigl)$,
\item
F(b) =[f(b),f],
\item
$\mbox{Im} A_{\bigl|G}$ = l'image de l'opérateur A restreint au sous \ev G,
\item
$A_k$ = la $k^{\mbox{ième} }$ colonne de A,
\item
$A_k(j)$ = la $j^{\mbox{ième} }$ coordonnée du vecteur $A_k$
(le c\oe fficient $a_{jk}$ de la matrice A),
\item
$c_{u_1},\ldots,c_{u_{ \tiny{NCD} }}$ les coordonnées distinctes de c,
\item
$I_i=\{j|c_j=c_{u_{ i }  }  \}$, $i=1,\ldots,NCD$,
\item
$N_i = \sharp I_i$ pour i=1,$\ldots,$NCD,
\item
$b_{v_1},\ldots,b_{v_{ q}}$ les coordonnées distinctes de b,
\item
$J_i=\{j|b_j=b_{v_{ i }  }  \}$, $i=1,\ldots,q$,
\item
$M_i = \sharp J_i$ pour i=1,$\ldots,$q.
\end{itemize}
\subsection{Indépendance, première approche}
\begin{lm}
On a les propriétés suivantes :
\begin{enumerate}
\item dim G = NCD
\item $1\in G$
\item $1\in(\mbox{Im} A_{\bigl|G})^{\bot}$.
\end{enumerate}
\end{lm}
%
%
%
%
%
%
%
%

%
%
%
%
\subsection*{\underline{\textbf{Démonstration}}}
\begin{enumerate}
\item[]
\small
\begin{enumerate}
\item
Si nous posons, pour tout $j=1,\ldots,NCD$ :
$\varepsilon_j(c) =\dps\sum_{k\in I_j}e_k$,
on peut écrire $c = \dps\sum_{j=1}^{NCD}c_{u_j}\varepsilon_j(c)$.\\
Pour tout g, g(c)
$\in [\varepsilon_1(c),\ldots,\varepsilon_{ \tiny{NCD} }(c)]
\equiv\widetilde{G}$ donc $G\subset \widetilde{G}$.\\

Soit z $\in\widetilde{G}$, $z = \dps\sum_{k=1}^{NCD}z_k\varepsilon_k(c)$,
il existe g mesurable sur $\Bigl(\R,B(\R)\Bigl)$
(l'interpolation de Lagrange nous fournit un
polynôme qui convient) telle que $g(c_{u_k} )=z_k$ pour $k=1,\ldots,NCD$.\\
\item
Il suffit de prendre g=1.\\
\item
$\Bigl(1,Ag(c)\Bigl)_{\R^n} =
\dps\sum_{k=1}^np_kg(c_k)-\dps\sum_{j,k}p_jp_kg(c_k)=\dps\sum_{k=1}^np_kg(c_k)-\dps\Biggl(\sum_{j=1}^np_j\Biggl)\sum_{k=1}^np_kg(c_k)=0$.
\end{enumerate}
\end{enumerate} \normalsize$\triangle$
\subsection*{Remarques}
\begin{enumerate}
\item
dim $\Bigl( \mbox{Im} A_{\bigl|G}\Bigl) = NCD -1$ puisque $1\in G$.
\item
On peut raisonner pour le point 1) du lemme 2.3.1. sur l'espace probabilisé.
G est l'image, par l'isomorphisme naturel,  de l'\ev engendré par les
indicatrices des ensembles
$c^{-1}\Bigl(\{ c_{  u_{ \nu  }  } \}\Bigl)$, $\nu=1,\ldots,NCD$.
Ces ensembles étant disjoints et   au nombre de NCD, la dimension de cet
espace est NCD.
\item
 Si c a toutes ses coordonnées distinctes alors :
    \begin{itemize}
    \item Sur $\R^n$\\
dim $\Bigl( \mbox{Im} A_{\bigl|G}\Bigl) =n-1$ donc dim
$(\mbox{Im} A_{\bigl|G})^{\bot}=1$ et par suite la \va constante
égale à 1 est la seule \va
indépendante de c (pour être précis il faut dire que toutes les
\vas de l'\ev engendré par la \va 1 sont indépendantes de c).\\
   \item Sur$L^2(\mu)$\\
  La tribu engendrée par c, est la tribu discrète, donc seule la \va 1 lui
   est indépendante.
  \end{itemize}
 Dans la suite on pourra supposer $NCD<n$.
\item
Si toutes les coordonnées de c sont égales alors :
 \begin{itemize}
   \item Sur $\R^n$\\
   $\mbox{Im} A_{\bigl|G}=\{0\}$. C'est à dire, toutes les \vas sont dans
   $\Bigl(\mbox{Im} A_{\bigl|G}\Bigl)^{\bot}$ et par
   suite, lui sont indépendantes. \\
  \item Sur $L^2(\mu)$\\
   La tribu engendrée par c, est la tribu grossière qui est indépendante de
   toute tribu.
 \end{itemize}
Dans la suite on pourra supposer $NCD>1.$
%
%
%
%
%
%
%
\item
On va utiliser la caractérisation suivante de l'indépendance de b et c :\\
$$\mbox{Im} A_{\bigl|G}\subset\Bigl(F(b)\Bigl)^{\bot}.$$
\end{enumerate}
{\bf{A partir de maintenant on prend :}} $\mbox{NCD}\in\{2,\ldots,n-1\}$
\begin{lm}
Soit une \va b telle que dim F(b) = q $\not=1$. Les éléments,
$\Bigl(x_i\Bigl)_{\tiny{i=1,\ldots,n}}$ de $\Bigl(F(b)\Bigl)^{\bot}$ vérifient un système de la forme :
$$\left\{\begin{array}{l}
                                             \dps\sum_{j\in J_i}x_j  =   0\\
                                             i=1,\ldots,q
   \end{array}\right.$$
\end{lm}
\subsection*{\underline{\textbf{Démonstration}}}
\begin{enumerate}
\item[]
\small
On a vu ( point 1 de la démonstration du lemme 2.3.1.)  que le système
$\{\varepsilon_1(b),\ldots,\varepsilon_{q}(b)  \}$ forme un base de F(b)
donc :
$$x\in F^{\bot}\iff\forall k\in \{1,\ldots,q\}  \ \ \Bigl( x ,
\varepsilon_k(b) \Bigl) = 0\iff \Bigl( x , \sum_{j\in J_k}e_j\Bigl) =0$$
\hspace*{3,85cm}$\iff \left\{\begin{array}{l}
\dps\sum_{j\in J_k}x_j =0\\
k=1,\ldots,q
\end{array}\right.$\\
\end{enumerate} \normalsize$\triangle$
\begin{lm}
S'il existe k $\in\{1,\ldots,NCD\}$ tel que $N_k=1$ alors la \va constante
égale à 1 est la seule \va
indépendante de c.
\end{lm}
\subsection*{\underline{\textbf{Démonstration}}}
\begin{enumerate}
\item[]
\small
Supposons $N_k=1$ et que b$\not=1$ est indépendante de c. Alors par le
théorème 2.2.1. pour toutes fonctions boréliennes $f,g$
$\Bigl( f(b) , Ag(c) \Bigl) =0$.\\
Prenons en particulier g(c) = $e_{u_k}$ alors : \\
\begin{itemize}
\item
Pour f=I, on a $\dps\sum_{i=1}^n b_ip_i=b_1$
\item
Pour $f(x)=x^2$ on a $\dps\sum_{i=1}^n b_i^2p_i=b_1^2$
\end{itemize}
Ce qui donne :
$\dps\sum_{i=1}^n b_i^2p_i = \Bigl( \sum_{i=1}^n b_ip_i\Bigl)^2$
et donc $b_i=b_j$ pour tout (i,j)$\in \{1,\ldots,n\}^2$,
c'est à dire $b\in[1]$.\\
\end{enumerate} \normalsize$\triangle$
%
%
%
%
%
%
%
%
%
%
%
%
%
%
%
%
\subsection{Conditions nécessaires d'indépendance}
\begin{pp}
On suppose que pour tout j dans $\{1,\ldots,NCD\}$, $2\leq N_j< n$.\\
Pour tout $k\in\{1,\ldots,NCD\}$  les vecteurs  :
$$\left\{\begin{array}{l}
 V_{\nu} = \dps\sum_{i\in I_{\nu}}A_i\\
 \nu\in\{1,\ldots,NCD\}-\{k\}
    \end{array}\right.$$
forment une base de $\mbox{Im}A_{\bigl|G}$.\\
De plus, si on note $V_{\nu}(j)$ la $j^{\mbox{ième}}$
coordonnée de $V_{\nu}$ alors :
 $$V_{\nu}(j)\not=0\mbox{ pour tout }
 (\nu,j)\in\{1,\ldots,p\}\times\{1,\ldots,n\}.$$
 \end{pp}
 \begin{co}
 On conserve les notations de la proposition 2.3.1.\\
 Prenons $b\not=1$ et posons
 $$F(b)^{\bot}=\left\{ x=\left(\begin{array}{c}
                                                      x_1\\
                                                      \vdots\\
                                                      x_n
 \end{array}\right)
  \Biggl | \tiny{\left\{\begin{array}{l}
   \dps\sum_{j\in J_i}x_j  =   0    \\
   i=1,\ldots,q
    \end{array}\right.}
                         \right\}$$
Deux conditions nécessaires pour que
$\mbox{Im} A_{\bigl|G}\subset F(b)^{\bot}$ sont :\\
\begin{itemize}
\item
$M_{\nu}\geq2$ pour $\nu=1,\ldots,q$
\item
$q\leq\dps\min_{ j\in\{1,\ldots,NCD\}-\{k\} }  \Bigl\{N_{\nu},n-N_{\nu}\Bigl\}.$
\end{itemize}
\end{co}
\subsection*{\underline{\textbf{Démonstration de la proposition 2.3.1.}}}
\begin{enumerate}
\item[]
\small
$\bullet$
$\Bigl( A\varepsilon_i(c)\Bigl)_{i\not=k}$est un système générateur de
$\mbox{Im} A_{\bigl|G}$ puisque on peut prendre
$\Bigl\{1,\varepsilon_i(c),i\in\{1,\ldots,NCD\}-\{k\}\Bigl\}$ comme base
de G . On sait d'autre part que dim
$\Bigl( \mbox{Im} A_{\bigl|G}\Bigl) = NCD -1$.  Ce système forme
donc une base de $\mbox{Im} A_{\bigl|G}$.\\
$\bullet$  $V_{\nu}(j)  = \dps\sum_{i\in I_{\nu}}A_i(j) =
                    \left\{ \begin{array}{lc} p_j\Bigl(\dps\sum_{i\in I_{\nu}}p_{ i }  \Bigl)  &>0\\
                                                          \mbox{ou} &\\
                                                          p_j\Bigl(\dps\sum_{i\in I_{\nu}}p_{ i }-1 \Bigl)     &<0
                                \end{array}\right.$\\
\end{enumerate} \normalsize$\triangle$
\subsection*{\underline{\textbf{Démonstration du corollaire 2.3.1.}}}
\begin{enumerate}
\item[]
\small
$\bullet$ S'il existe $\nu$ tel que $ N_{\nu}=1$ alors l'équation
$x_{v_{\nu}} =0$ est une des équations définissant $F(b)^{\bot}$.
D'après la proposition 2.3.1. aucune coordonnée des vecteurs de la base de
$\mbox{Im} A_{\bigl|G}$ n'est nulle donc $\mbox{Im} A_{\bigl|G}$ n'est pas
contenue
dans $F(b)^{\bot}$.\\
%
%
%
%
%
%
%
%
%
%
%
%
%
%
$\bullet$ Fixons j dans $\{1,\ldots,NCD\}-\{k\}$,
$V_j =\dps\sum_{i\in I_j}A_{i}$. On additionne  $N_j$ colonnes donc
$V_j$ a $N_j$
termes $<0$ et $n-N_j$ termes$>0$.\\
En effet : $V_j=\left(\begin{array}{c}
             p_1\Bigl(\dps\sum_{i\in I_j}p_{i} -\widetilde{\varepsilon}_1(I_j)\Bigl)\\
             \vdots\\
             p_n\Bigl(\dps\sum_{i\in I_j}p_{i} -\widetilde{\varepsilon}_n(I_j)\Bigl)
             \end{array}\right)$\\
où $\widetilde{\varepsilon}_{\nu}(I_j) =\left\{\begin{array}{lcl}
              1  & \mbox{si}  & \nu\in  I_j\\
              0  &  \mbox{sinon}  &
              \end{array}\right.$\\
et les $\widetilde{\varepsilon}_{\nu}(I_j)$ valent exactement $N_{\nu}$
 fois 1, j étant fixé.\\
$\bullet$
Si on a plus de $N_j$ (resp. $n-N_j$) équations, $q>N_j$ (resp. $q>n-N_j$),
 alors il y en a une dans laquelle il n'y a que des termes  $>0$
 (resp. $<0$),
elle ne s'annule donc pas. On en déduit qu'on doit avoir
$q\leq \min (N_j,n-N_j)$ pour que $V_j$ soit dans $F(b)^{\bot}$.
Un raisonnement similaire pour chaque vecteur de la base nous donne :
$$q\leq \dps\min_{ j\in\{1,\ldots,NCD\}-\{k\} }\Bigl(N_j,n-N_j\Bigl)$$
\end{enumerate} \normalsize$\triangle$
\begin{thm}
Soit une \va c telle que dim G = p avec 1$\leq p < n$. Une condition
 nécessaire pour que les \vas c,$b_1,\ldots,b_N$, $b_i\not=1$
 $i=1,\ldots,N$, forment un système
globalement indépendant est :
$$ p\prod_{i=1}^Nn(b_i)\leq n$$
où $n(b_i) = \mbox{dim }F(b_i)$ $ i=1\ldots,N$.
\end{thm}
\subsection*{\underline{\textbf{Démonstration}}}
\begin{enumerate}
\item[]
\small
On regarde les tribus engendrées par chacune des \vas :\\
$\F_c =\sigma \Bigl( \Gamma_1^0,\ldots, \Gamma_p^0\Bigl)$ où
 $\Gamma_j^0=\{a_i | i\in I_j\}=c^{-1}\Bigl(\{c_{u_j} \}\Bigl)$.\\
Et pour $i=1,\ldots,N$ $\F_{ b_i}=\sigma \Bigl(  \Gamma_1^i,\ldots,
 \Gamma_{ n(b_i) }^i\Bigl)$\\
où $ \Gamma_{\nu}^i = b_i^{-1}\Bigl(\{ x_{r_{\nu} }^i  \}\Bigl)$,
$x_{r_{1} }^i,\ldots,x_{r_{n(b_i)} }^i$ étant les coordonnées distinctes
de $b_i$.\\
Pour avoir l'indépendance il faut que les intersections
$ \Gamma_{k_0}^0\bigcap\cdots\bigcap  \Gamma_{k_N}^N$ soient non vides pour
tout
$(k_0,\ldots,k_N)\in\{1,\ldots,p\}\times\{1,\ldots,n(b_1)\}
\times\cdots\times\{1,\ldots,n(b_N)\}$.
Ce qui nous fait $p\dps\prod_{i=1}^Nn(b_i)$ conditions à vérifier. \\
D'autre part on a $\Bigl(A_{k_0}^0\bigcap\cdots\bigcap
A_{k_N}^N\Bigl)\bigcap\bigl(A_{l_0}^0\bigcap\cdots\bigcap
A_{l_N}^N\Bigl)=\emptyset$
pour tout :\\
$\Bigl((k_0,\ldots,k_N),(l_0,\ldots,l_N)\Bigl)\in\Bigl(\{1,\ldots,p\}
\times\{1,\ldots,n(b_1)\}\times\cdots\times\{1,\ldots,n(b_N)\}\Bigl)^2$ \\
 puisque pour tout $j\in \{ 0,\ldots,N \}$ et tous $\eta\not=\nu$,
 $ \Gamma_{\eta}^j\bigcap  \Gamma_{\nu}^j=\emptyset$.\\
Ainsi pour avoir l'indépendance on doit nécessairement avoir
$p\dps\prod_{i=1}^Nn(X_i)\leq n=\sharp \Omega$.\\
\end{enumerate} \normalsize$\triangle$
%
%
%
%
%
%
%
%
%
%
%
%
%
%
%
%
%
\begin{co}
Le nombre maximun de \vas globalement indépendantes est :
$$N_{\mbox{max}}  = \max\Bigl\{ k\in\N | 2^{k-1}\leq n\Bigl\}.$$
\end{co}
\subsection*{\underline{\textbf{Démonstration}}}
\begin{enumerate}
\item[]
\small
Pour tout  ensemble, $E(N) =\{b_1,\ldots,b_N\}$, de \vas globalement
indépendantes on a :
\begin{itemize}
\item
$1\in E(N)$
\item
$\dps\prod_{i=2}^kn(b_i)\leq n$
\item
$N_{\mbox{max} } =\max \Bigl\{ N|  \dps\prod_{i=2}^Nn(b_i)\leq n
 \mbox{ avec  }b_2,\ldots,b_N  \mbox{ ind., }\not=1 \Bigl\}$
\end{itemize}
Pour toute  \va $b_i$, $n(b_i)\geq2$, donc $\dps\prod_{i=2}^Nn(b_i)
 \geq 2^{N-1}$ et par suite :\\
$N_{\mbox{max} }\leq\max\Bigl\{N|2^{N-1}\leq n\Bigl\}$.\\
Pour finir il suffit de prendre $b_i =x_1^i1_{A^i_1}+x_2^i1_{A^i_2}$
$i=1,\ldots,N$ pour avoir $N_{\mbox{max} }\geq\max\Bigl\{N|2^{N-1}\leq
n\Bigl\}$.\\
\end{enumerate} \normalsize$\triangle$
%
%
%
%
%
%
\part{DENSITÉ DE VARIABLES ALÉATOIRES}
\chapter{Polynômes de Wick}
\section{Introduction}
On fixe :  $\bullet$ (M,$\Sigma$,$\mu$) un espace probabilisé, \\
\hphantom{on fixe : } $\bullet$ f,\;g$\in L^p(M,\Sigma,\mu)$
$\forall p\geq1$ ,\\
on note : $\bullet$ $m_p =\int_Mf^pd\mu$, $m_p^a =\int_M|f|^pd\mu$,
nombres qui sont finis, \\
\hphantom{on note : }$\bullet$ $H_n=[ f^k,k=0,\ldots,n]$.
\subsection*{Remarque}
On peut se demander quelle est la dimension de $H_n$. Ici on va prendre
comme hypothèse,
dim($H_n$)= n+1, c'est le cas quand le support de f est d'intérieur non
vide. Cette condition n'est pas automatique, comme le montre l'exemple
suivant :\\
en prenant $\mu = \frac{1}{2}(\delta_{-1}+\delta_1)$,
$f=\left \{ \begin {array}{cl}
  +1 &\\
  -1  &
  \end{array}\right.$
on obtient $H_0 =[1]$, $H_1=[1,f]$, $H_k=H_1$ pour tout $k\geq2$.
\begin{de}
On définit l'opérateur $D_f$ sur $H_n$ par :
$$ D_f(\sum_{k=0}^na_kf^k) = \sum_{k=1}^nka_kf^{k-1}.$$
\end{de}
\begin{pp}
On peut définir cet opérateur sur $\dps\bigcup_n H_n$ et le prolonger dans
$L^2\smm$ :
$$D_f(\sum_{k\geq0}a_kf^k) = \sum_{k\geq1}ka_kf^{k-1},$$
sur l'ensemble des séries convergentes, $\dps\sum_{k\geq0}a_kf^k$,  dans
$L^2\smm$.
\end{pp}
\begin{pp}On a les propriétés suivantes :\\
1) $D_f$ est linéaire,\\
2) si h = af+b alors $H_n =[h^k, k=0,\ldots,n]$ et
$$D_f[ (af + b)^n] = aD_h(h^n),$$
3) on a :
$$D_f\biggl( \  \sum_{k=0}^{N_1} a_kf^k\sum_{k=0}^{N_2}b_kf^k Ê\ \biggl ) = D_f\biggl( \ \sum_{k=0}^{N_1}a_kf  ^k \  \biggl)\sum_{k=0}^{N_2}b_kf^k
+\sum_{k=0}^{N_1}a_kf^kD_f\biggl ( \ \sum_{k=0}^{N_2}b_kf^k \  \biggl).$$
\end{pp}

\subsection*{\underline{\textbf{Démonstration}}}
\begin{enumerate}
\item[]
\small
1) Simples calculs.\\
2) On développe $(af+b)^n$ avec la formule du binôme de Newton.\\
3) On développe le produit, on lui applique $D_f$ :
$$\sum_{n=0}^{N_1+N_2}\sum_{k=0}^nna_kb_{n-k}f^{n-1}$$
on écrit $n=(n-k)+k$ et on a le résultat.
\end{enumerate} \normalsize$\triangle$
%
%
%
%
%
%
%
%
%
\subsection*{Remarque}
Quand on se place sur l'espace image, $D_f$ devient une vraie dérivation. \\
 Plus précisément, si  à $g\in\bigcup H_n$ on associe le polynôme P(x),
  dans $L^2(\mu_f)$,
$\bigl(g=\dps\sum_{finie} a_kf^k\mapsto \sum a_kx^k=P(x)\bigl)$ on a
$D_f(g) =\dps\frac{dP}{dx}of$.
\begin{de}
On définit une suite $(f^{\widetilde{n} } )_{n\geq0}$, appelée
{\bf{ $n^{\mbox{ième} }$ puissance de Wick de f}}, par :
$$\left\{
\begin{array}{lcl}
f^{\widetilde{0} } & = & 1\\
D_f(f^{\widetilde{n+1} }) & = & (n+1)f^{\widetilde{n} }\mbox{, }n\geq1
\mbox{, vérifiant }
\left\{
\begin{array}{lcl}
 \int_M f^{\widetilde{n+1} }d\mu &=&0\\
  f^{\widetilde{n+1} }&\in&H_{n+1}.
  \end{array}
  \right.
\end{array}
\right.$$
\end{de}
\begin{pp}
On a :\\
$$ \mbox{si }  h = af + b\mbox{, }b\not=0\mbox{, alors }h^{\widetilde{n} }=
 a^nf^{\widetilde{n} }.$$
\end{pp}
\subsection*{\underline{\textbf{Démonstration}}}
\begin{enumerate}
\item[]
\small
Par récurrence :\\
$f^{ \widetilde{0} }   = h^{ \widetilde{0} }   =1$,  $f^{ \widetilde{1} }
=  f - m_1$,
$h^{ \widetilde{1} }   =  h - E(h) =a(f-m_1) = a f^{ \widetilde{1} }$\\
\begin{tabbing}
Supposons que $h^{\widetilde{n} }= a^nf^{\widetilde{n} }$.\\
$D_f (h^{ \widetilde{n+1} } )$ \=  $= aD_h (h^{ \widetilde{n+1} } )$
(proposition 3.1.2.)\\
  \>       $= a(n+1)h^{ \widetilde{n }} \ \
\stackrel{\mbox{H.R.}}{=} \ a^{n+1}(n+1)f^{ \widetilde{n }}$ \\
 \> $= a^{n+1}D_f( f^{ \widetilde{n+1 }})$
\end{tabbing}
$D_f(h^{ \widetilde{n+1 }} -a^{n+1} f^{ \widetilde{n+1 }} ) = 0
\Rightarrow h^{ \widetilde{n+1 } }
-a(n+1) f^{ \widetilde{n+1 }} = $ const.
;\\ finalement comme $E(h^{ \widetilde{n+1 }}-a^{n+1} f^{ \widetilde{n+1 }} )
 = 0$, const.$=0$.
\end{enumerate} \normalsize$\triangle$
\section{Exponentielle de Wick}
- On suppose qu'on a $\exp|t|\cdot |f| \in L^2\smm$ pour tout $t\in V$
(V étant un voisinage de zéro).\\
\\
- La transformée de Laplace de f est donnée par :
$\phi_f(t) = \int_M\exp(tf)d\mu$.
 Une condition suffisante pour que $S_n(t,f) = \dps \sum_{k=0}^n
\frac{t^k}{k!}f^k$ converge dans $L^1\smm$ est que :
$$\dps \sum_{n\geq0} \frac{|t|^n }{ n! }
m_n<\infty.$$
\begin{de}On définit :
$$\E (tf)\equiv\frac{ e^{tf} }{ \phi_f(t) }$$
qu'on appelle \bf{exponentielle de Wick}.
\end{de}
\subsection*{Exemples : \bf{EXPONENTIELLES DE WICK}}
1) pour $f\sim {\cal N}(0,1)$:     \hspace{2,4cm}
$ \E(tf)   = \dps e^{tf-\frac{ t^2 }{ 2 } }$,\\
on reconnaît la fonction génératrice des polynômes de Hermite,
\begin{tabbing}
2) pour $f\sim exp(\lambda)$ ($\lambda>0$) \hspace{1cm}\=:  \hspace{1cm}                                       \=  $ \E(tf) = \dps e^{tf}\frac{\lambda -t}{\lambda}$,\\
3) pour $f\sim\Gamma(a,b)$                           \>:
  \> $\E(tf) = e^{tf} \biggl( \dps\frac{b}{b-t}\biggl)^{-a},$\\  \\
4) pour $f\sim \Gamma(1/2,1/2)$                   \>:
 \>$ \E(tf)   = e^{tf}\sqrt{1-2t}$,\\ \\
5) pour $ f\sim\Gamma\Bigl(\alpha{a_1\choose  b_1},\beta{a_2\choose b_2}\Bigl)$\>:
  \>$e^{tf}\biggl(\dps \frac{ b_1 }{ b_1-\alpha t }\biggl)^{-a_1} \biggl(\dps\frac{ b_2 }{
b_2-\beta t }\biggl)^{-a_2},$\\  \\
6) pour $f\sim {\cal P}(a)$, $a>0$                               \>:                                                              \>$\E(tf) = \dps\exp\bigl(tf-a(e^t-1)\bigl)$.
\end{tabbing}
\begin{pp}
Sous les conditions énoncées plus haut, $\E (tf)$ est, dans un voisinage de
zéro,
 développable en série entière (convergente dans $L^1\smm$) et on a :
$$E(\E ) =1.$$
\end{pp}
\subsection*{Notations}
- $\E(tf)  =\dps\sum_{n\geq0}W_n(f)\frac{ t^n }{ n! }$,\\
\\
- $\phi_f(t) =\dps\sum_{n\geq0}m_n \frac{  t^n  }{  n!  },  $\\
\\
- $(\phi_f(t))^{-1} =\dps\sum_{n\geq0}a_n \frac{  t^n  }{  n!  },  $\\
\\
- $S_n(t,f) =\dps\sum_{k=0}^n\frac{t^k}{k!}f^k$,\\
\\
- Pour $x\in\R$ et $k\in\N$ on pose : $\dps {x\choose k} =
\dps\frac{ x(x-1)\cdot\ldots\cdot(x-k+1) }{k! }$,\\
\\
- On notera $[x]$ pour la partie entière du réel $x$.
\\
- On a : $\dps \sum_{k=0}^{n} {n \choose k}m_k a_{n-k}=0.$
\subsection*{\underline{\textbf{Démonstration}}}
\begin{enumerate}
\item[]
\small
Il existe un voisinage, V', de zéro sur lequel $\phi_f(t)$ ainsi que
$(\phi_f(t))^{-1}$
sont développables en série entière ($\phi_f(t)$ ne s'annule pas).
Voyons la
convergence dans $L^2\smm$ :\\
Posons $A_n =\dps \sum_{k+l=n}\frac{t^k}{k!}\frac{t^l}{l!}a_kf^l$ on a :\\
\begin{tabbing}
$\dps\sum_{n=1}^{m}|A_n|$ \= $\leq \dps\sum_{k=1}^m
\sum_{l=1}^m\frac{|t|^k}{k!}\frac{|t|^l}{l!}|a_k|\cdot|f|^l$\\
                                      \>
$\leq  \underbrace{\dps\sum_{k=1}^m\frac{|t|^k}{k!}|a_k|}_{C(t)}
\dps\sum_{l=1}^m\frac{|t|^l}{l!}|f|^l$
\end{tabbing}
$|S_n(t,f)|\leq\dps\sum_{n=1}^{m}|A_m|\leq C(t)\exp|t|\cdot|f|\Rightarrow$
 convergence dans $L^1\smm$.\\
La convergence dans $L^2\smm$ en découle car :\\
$|S_n(t,f)| - \frac{\exp(tf)}{\phi_f(t)}|^2\leq C(t)\exp(2|t||f|)$.
\end{enumerate} \normalsize$\triangle$
\begin{pp}
 On a :

$$f^{\widetilde{n} } = W_n(f) = \dps\sum_{k=0}^{n} {n \choose k}a_kf^{n-k}.$$
\end{pp}
\subsection*{\underline{\textbf{Démonstration}}}
\begin{enumerate}
\item[]
\small
$f^{ \widetilde{0} }=1$, $W_0(f) =1$; $f^{ \widetilde{1}  }= f-m_1$, $W_1(f) =f+a_1$, mais
$\dps\sum_{k=0}^n{n\choose k} a_km_{n-k} = 0$ pour $n\geq1$, donc (n=1)
$a_1 = -m_1$.\\
Supposons que $f^{ \widetilde{n} } = W_n(f) =\dps\sum_{k=0}^n{n\choose k}
 a_kf^{n-k}$
 comme $D_ff^{ \widetilde{n+1} } = (n+1) f^{ \widetilde{n} }$
 on en déduit :\\
\begin{tabbing}
$f^{ \widetilde{n+1} }$ \= $= (n+1)\dps\sum_{k=0}^n{n\choose k} a_k
\biggl(\ \
\frac{f^{ \widetilde{n-k+1} } - m_{n-k+1}  }{  n-k+1}\ \ \biggl)$\\
  \>$ = \dps \sum_{k=0}^{n+1} { n+1\choose k} a_k
 f^{ \widetilde{n-k+1} } - a_{n+1} - \dps\sum_{k=0}^n{ n+1\choose k}
 a_km_{n-k+1}$\\
  \>$ =\dps  \sum_{k=0}^{n+1} { n+1\choose k} a_k  f^{ \widetilde{n-k+1} } $
\end{tabbing}
\end{enumerate} \normalsize$\triangle$
%
%
%
%
%
%
%
\subsection*{Exemples : {\bf{ FORMULATIONS EXPLICITES}} }
1) pour $f\sim {\cal N}(0,1)$
$$W_n(f) \equiv H_n(f) =  n!\sum_{k=0}^{[n/2]} \frac{ (-1)^kf^{n-2k} }{ 2^kk!(n-2k)! }\mbox{ \  }
(n^{\mbox{ième}}\mbox{ polynôme de Hermite)},$$
2) pour $f\sim exp(\lambda)$ ($\lambda>0$) :
$$W_n(f) \equiv E^{\lambda}_n(f)  = f^n - \frac{n}{\lambda} f^{n-1},$$
3) pour $f\sim\Gamma(a,b)$ :\\
$$ W_n(f) \equiv\Gamma^{ab}_n(f) = n!\sum_{k=0}^n{a \choose n-k}
 \frac{f^k}{k! }
\Bigl( \frac{ -1 }{b} \Bigl)^{n-k},$$
4) pour $f\sim \Gamma(1/2,1/2)$ :
$$W_n(f) \equiv\Gamma_n^{\frac{1}{2} \frac{1}{2}}(f) = f^n -
\frac{ n! }{ 2^{n-1} }
\sum_{k=0}^{n-1}\frac{ 2^k }{ n-k } { 2(n-k-1) \choose n-k-1}
\frac{f^k}{ k!  },$$
5) pour $ f\sim\Gamma\Bigl(\alpha{a_1\choose b_1},\beta{a_2\choose b_2}
\Bigl)$ :
$$W_n(f) \equiv Q_n(f) = \sum_{k=0}^n(-1)^k{n\choose k}k!
\sum_{p=0}^k{a_1\choose p}
{a_2\choose k-p} \biggl( \frac{ \alpha }{  b_1 }\biggl)^p\biggl(
\frac{ \beta }{
b_2 }\biggl)^{k-p}f^{n-k},$$
6) pour $ f\sim {\cal P}(a)$, $a>0$ ($D = \frac{d}{da}$) : \\
$$ W_n(f) \equiv P^a_n(f) = \sum_{k=0}^n{ n \choose k} f^{n-k} \Bigl[\ e^{a}
\bigl( aD \bigl) ^k\cdot \bigl( e^{-a}\bigl) \Bigl].$$
\begin{pp}
On a les propriétés suivantes :\\
1) $D_f(\E (tf)) = t\E (tf)$,\\
\\
2) $\E (tf) = \dps\frac{\exp\bigl(t(f-m_1)\bigl) }{ \phi_{f-m_1}(t) },$\\
\\
3) $\E(\lambda(sf+tg)) =\E(\lambda sf)\E(\lambda tg)\dps
\frac{ \phi_f(\lambda s)
\phi_g(\lambda t) }{ \phi_{(f,g)} (\lambda s,\lambda t) }.$
\end{pp}
\subsection*{\underline{\textbf{Démonstration}}}
\begin{enumerate}
\item[]
\small
1) On utilise le prolongement de $D_f$.\\
2) Simple vérification.\\
3) Il suffit d'écrire : $\phi_{sf+tg}(\lambda) = \dps \frac{\phi_{sf+tg}
(\lambda) }
{ \phi_{sf}(\lambda)\phi_{tg}(\lambda) }
\phi_{sf}(\lambda)\phi_{tg}(\lambda)$.
\end{enumerate} \normalsize$\triangle$
\subsection*{Remarque}
$$\E(\lambda(sf+tg)) =\frac{  \exp(\lambda(sf+tg))  }{  \phi_{(f,g)}
 (\lambda   s,\lambda t) }=
 \frac{  \exp(\lambda(sf+tg))  }{  \phi_{(sf+tg)}(\lambda) }.$$
Pour avoir des précisions sur les propriétés des puissances et de
l'exponentielle
de Wick confère {\bf [1]}.
%
%
%
%
%
%
\section{Propriétés des puissances de Wick}
\subsection{Introduction}

Il est connu (voir par exemple {\bf [2]}) que les polynômes orthogonaux
classiques ont une somme
importante
de propriétés en commun dont :\\
- ils ont $\bullet$ une formule de récurrence\\
\hphantom{- ils ont }$\bullet$ une équation différentielle\\
\hphantom{- ils ont }$\bullet$ une formule de Rodriguez...\\
(en fait, toutes ces propriétés découlent, y compris les propriétés d'orthogonalité,
de la formule de Rodriguez, confère [2]).
On a vu que les puissances de Wick de la \va X (ici on s'intéresse aux
propriétés
polynômiales on notera X pour f pour le notifier) sont des polynômes en X de
formule explicite :
    $$W_n(X) = \dps \sum_{k=0}^n{ n \choose k} a_kX^{n-k}.$$
On peut se demander si cette suite de polynômes est orthogonale, cette
situation nous permettrait de puiser dans les outils de la
théorie des polynômes orthogonaux.\\
Quand ces polynômes ne sont pas orthogonaux, existe-t-il des propriétés communes à
 toutes ces familles (formule de récurrence, équation différentielle...).
\subsection{Orthogonalité}
\begin{thm}
Etant donnée une \va X telle que $\exp|t||X| \in L^2\smm$, $t\in  V$, on a :\\
$$\forall n\not=m \int_MW_n(X)W_m(X)d\mu = 0 \Leftrightarrow X\sim{\cal N}
(m_1,m_2-m_1^2).$$
\end{thm}
\subsection*{\underline{\textbf{Démonstration}}}
\begin{enumerate}
\item[]
\small
On suppose $m_1 =0$ (quitte à prendre $X-m_1$).\\
$X^{ \widetilde{1} } = X \Rightarrow \int_MXW_n(X)d\mu =0 $ pour $n\not=1$.\\
$\phi_X(t)\E(tX) = \exp(tX)\Rightarrow \dps \sum_{k=0}^n{n\choose k}W_{n-k}m_{n-k}=X^n$
\begin{tabbing}
Pour $n\not=1$ $m_{n+1} = E(X\cdot X^n)$ \= = $\dps\sum_{k=0}^n{n\choose k}m_{n-k}
E(\ XW_k(X)\ )$\\
                                                                  \> = ${n\choose 1}m_{n-1}m_2$\\
                                                                \> = $nm_{n-1}m_2 $
\end{tabbing}
$\Rightarrow\left\{
\begin{array}{lclcl}
m_{2n+1} & = &2n m_{2n-1}m_2 & \mbox{,} & m_1=0\\
m_{2n} & =& (2n -1)m_{2n-2}m_2 & \mbox{,} & m_0=1
\end{array}
\right.$\\[.5cm]
$\Rightarrow\left\{
\begin{array}{lcl}
m_{2n+1} & = &0\\
m_{2n} & =& (2n -1)(2n-3)\cdot\ldots\cdot 1m_0m_2^n = \frac{ (2n)!  }{2^nn! } m_2^n
\end{array}
\right.$\\[.5cm]
On a donc une \va qui a les moments et la transformée de Laplace d'une \va suivant
une loi ${\cal{N}}(0,m_2)$, on peut donc conclure que $X\sim{\cal{N}}(0,m_2)$.
\end{enumerate} \normalsize$\triangle$
\subsection*{Remarque}
Ce théorème règle le problème de l'orthogonalité. Dans le cas où X suit une loi normale
on a vu que $W_n = H_n$.  C'est un objet bien connu (dans le cadre
classique de la décomposition en chaos de $L^2\bigl((\R, {\cal B}(\R), \frac{1}{\sqrt{2\pi} }
e^{-x^2/2}) \  \bigl)$. On va essayer de trouver des propriétés communes à toutes
les suites de polynômes de Wick.
%
%
%
%
%
%

\subsection{Formules de récurrence}
On sait que $\E(tX) =  \dps\sum_{n\geq0}W_n(X)\frac{  t^n  }{  n!  }$ où
 $\E(tX)$ est une
fonction connue qu'on veut développer en série entière. En
fait $\E(tX)$ n'est rien d'autre que la fonction génératrice des puissances
de Wick, que nous
appellerons également polynômes de Wick.\\
On va puiser dans la théorie des fonctions génératrices
 (confère {\bf[3]}) pour trouver une formule de récurrence.
\subsection*{i) Première formule de récurrence}
\begin{thm}
Les polynômes de Wick vérifient la formule de récurrence : \\
$$\left\{
\begin{array}{lcl}
W_0(X)       &     =    & 1\\
                  &           &\\
W_1(X)       &    =     & X-m_1\\
W_n(X)       &    =     & (X-m_1)W_{n-1}(X)-\dps\frac{1}{n}\dps \sum_{k=0}^{n-2}{n \choose k}
 b_{n-k}W_k(X).
\end{array}
\right .$$
Où $b_n = \dps\sum_{k=0}^{n}{ n \choose k} k\cdot  m_ka_{n-k}$.
\end{thm}
\subsection*{\underline{\textbf{Démonstration}}}
\begin{enumerate}
\item[]
\small
On part de l'égalité:
$$e^{tX}(\phi_X(t))^{-1} = \sum_{k=0}^nW_n(X)\frac{ t^n }{ n! }\mbox{ pour } t\in V.$$
V étant un certain voisinage de zéro inclus dans le disque de convergence de
 la série.
Comme le membre de gauche est strictement positif on peut en prendre le
logarithme,
puis on dérive par rapport à t, puis on  multiplie par t pour obtenir :
$$tX -  \dps\sum_{n\geq1}n\cdot m_n\frac{ t^n }{ n!}\sum_{n\geq0}a_n
\frac{ t^n }{ n!} =
\frac{  \dps \sum_{n\geq1}nW_n\frac{ t^n }{ n!}       }{
\dps\sum_{n\geq0}W_n\frac{ t^n }{ n!}   }$$
Ce qui nous donne :
$$ \sum_{n\geq1} \ \ \biggl[ \ \ \frac{ XW_{n-1}  }{ (n-1)! }-
\frac{ 1 }{ n! }\sum_{k=0}^{n}
{ n \choose k}b_kW_{n-k}(X)\ \ \biggl]  t^n  = \sum_{n\geq1}
nW_n(X)\frac{ t^n }{
n!}$$
Où on a posé :
$$b_n = \sum_{k=0}^n  { n \choose k} k\cdot m_ka_{n-k}.$$
On identifie les c\oe fficients de chaque côté et on a le résultat attendu.
\end{enumerate} \normalsize$\triangle$
\subsection*{Exemples : \bf{FORMULES DE RECURRENCE 1}}
1) $X\sim {\cal{N}}(0,1)$ :
$$ H_n(X) = XH_{n-1}(X)- (n-1)H_{n-2}(X),$$
2) $X\sim \exp(\lambda)$ $(\lambda>0)$ :
$$E_n^{\lambda}(X) = (X-\frac{ 1                }{ \lambda}   ) E_{n-1}^{\lambda}(X) -(n-1)!
 \sum_{k=0}^{n-2} \frac{ E_k^{\lambda}(X) }{ k!\lambda^{n-k} },$$
3) $X\sim \Gamma( a , b )$ :
$$ \Gamma_n^{ab}(X) = (X-\frac{a}{b})\Gamma_{n-1}^{ab}(X) - a(n-1)!\sum_{k=0}^{n-2}
\frac{ \Gamma_k^{ab}(X) }{ k!} (\frac{1}{b}\ )^{n-k},$$
4) $X\sim \Gamma( 1/2 , 1/2 )$ :
$$\Gamma_n^{ \frac{1}{2}\frac{1}{2} }(X) = (X-1)\Gamma_{n-1}^{\frac{1}{2}\frac{1}{2} } (X)
 - (n-1)!\sum_{k=0}^{n-2}\frac{ 2^{n-k-1} }{ k! }
\Gamma_k^{\frac{1}{2} \frac{1}{2}  }(X),$$
5) $ X\sim\Gamma\Bigl(\alpha{a_1\choose b_1},\beta{a_2\choose b_2}\Bigl)$ :\\
$$Q_n(X) = \biggl( X - \alpha\frac{  a_1 }{ b_1 }- \beta\frac{ a_2 }{ b_2 } \biggl)Q_{n-1}(X)
 -   (n-1)!   \sum_{k=0}^{n-2}\frac{Q_k(X) }{ k!Ê}\biggl[ a_1\biggl(\frac{
\alpha }{b_1}\biggl)^{n-k} +a_2\biggl(\frac{\beta }{ b_2}\biggl)^{n-k}\biggl],$$
6) $ X\sim {\cal{P} }(a)$, $a>0$ :
$$ P_n^a(X) = (\ X -a\ ) P_{n-1}^a(X) - a\sum_{k=0}^{n-2}{n-1\choose k} P_k^a(X).$$
\subsection* {Remarque}
Cette formule de récurrence n'est pas très maniable, le calcul du
c\oe fficient $b_n$ est
souvent délicat, voir même impossible, ce qui nous oblige, pour obtenir la
formule
finale,  à refaire tous les calculs.
\subsection*{ii) Seconde formule de récurrence}
\begin{thm}
Les polynômes de Wick vérifient la formule de récurrence : \\
$$\left\{
\begin{array}{lcl}
W_0(X)       &     =    & 1\\
                  &           &\\
W_1(X)       &    =     & X-m_1\\
nW_n(X)       &    =     & \dps \sum_{k=1}^{n}{n \choose k} W_{n-k}(X)(kXm_{k-1} -nm_k)
\end{array}
\right .$$
\end{thm}
\subsection*{\underline{\textbf{Démonstration}}}
\begin{enumerate}
\item[]
\small
On agit de la même façon que précédemment mais en multipliant par
$\phi_X(t)\dps\sum_{k\geq0}W_n(X)\frac{ t^n }{ n! }$ au lieu de
$\dps\sum_{k\geq0}W_n(X)\frac{ t^n }{ n! }$.
\end{enumerate} \normalsize$\triangle$
\subsection*{Exemples : \bf{FORMULES DE RECURRENCE 2}}
1) $X\sim {\cal{N}}(0,1)$ :
$$ H_n(X) = -n!\sum_{k=1}^{[n/2]}\frac{H_{n-2k}(X) }{ (n-2k)!2^kk!} +(n-1)!X
\sum_{k=0}^{[(n-1)/2]}\frac{H_{n-2k-1}(X) }{ (n-2k-1)!2^kk! },$$
2) $X\sim \exp(\lambda)$ $(\lambda>0)$ :
$$E_n^{\lambda}(X) = (n-1)! \sum_{k=1}^{n} {n\choose k}
\frac{ E_{n-k}^{\lambda}(X) }{ (n-k)! }\Biggl(\frac{\lambda X -n}{ {\lambda}^{k} }\Biggl),$$
3) $X\sim \Gamma( a , b )$ :
$$ n\Gamma_n^{ab}(X) = \sum_{k=1}^{n}{n\choose k}
\frac{ a(a+1)\cdot\ldots\cdot(a+k-2) }{ b^{k-1} }\biggl(\ kX - \frac{n}{b}(a+k-1)\ \biggl)
\Gamma_{n-k}^{ab}(X),$$
4) $X\sim \Gamma( 1/2 , 1/2 )$ :
$$n\Gamma_n^{\frac{1}{2}\frac{1}{2} }(X) = \sum_{k=1}^n{n\choose k}{2k-2\choose k-1}
\frac{(k-1)!}{2^{k-1} }\biggl(\ kX-2n(k-\frac{1}{2} )\ \biggl)
\Gamma_{n-k}^{\frac{1}{2}\frac{1}{2} }(X),$$
5) $ X\sim\Gamma\Bigl(\alpha{a_1\choose b_1},\beta{a_2\choose b_2}\Bigl)$ : \\
\hspace*{4cm}$nQ_n(X) = \dps\sum_{k=0}^n{ n\choose k}Q_{n-k}(X)k!\ \cdot\Biggl\{$\\
$\dps\sum_{p=0}^{k-1}\Biggl[{ a_2+p-1\choose p}{ a_2+k-p-2\choose k-p-1}\biggl(
\dps\frac{ \alpha }{ b_1}
\biggl)^p\biggl(\dps\frac{\beta}{b_2}\biggl)^{k-1-p}\cdot\biggl( X -n
\Bigl(\dps \frac{a_2+k-p-1}{ k-p }\Bigl)\ \biggl)\Biggl]$\\
\hspace*{12cm}$-n\dps{ a_1+k-1\choose k}\biggl(\dps\frac{\alpha}{b_1}\biggl)^k\Biggl\},$
6) $ X\sim {\cal{P}}(a)$, $a>0$ :
$$ nP_n^a(X) = e^{-a}\sum_{k=1}^n{n\choose k}
\biggl(\ kX (aD)^{k-1}\ (e^a) - n(aD)^{k}\ (e^a)\ \biggl)P_{n-k}^a(X).$$
\subsection*{Remarque}
Cette formule de récurrence est beaucoup plus maniable, il n'y a plus de
calculs
intermédiaires à effectuer. Cela vient du fait qu'on n'utilise plus le
développement
en série de l'inverse de la transformée de Laplace, objet souvent difficile
à déterminer.
Cependant, dans le cas ${\cal{N}}$(0,1), la formule obtenue est moins
intéressante. Ainsi, en fonction des situations, on pourra choisir l'une où
l'autre de ces formules de récurrence.
\subsection{Equations différentielles}
\subsection*{i) Première équation différentielle}
\begin{thm}
Les polynômes de Wick vérifient une équation différentielle donnée, pour $n\geq0$, par :
$$nW_n(X)  -(X-m_1) W'_n(X) + \dps\sum_{k=2}^{n} \frac{ b_{k} }{ k! }W_n^{(k)}(x) = 0.$$
\end{thm}
\subsection*{\underline{\textbf{Démonstration}}}
\begin{enumerate}
\item[]
\small
Comme $W_n'(X) = nW_{n-1}(X)$ on a pour $n\geq k$  $W_k(X) = \frac{k!}{n!}W_n^{(n-k)}(X)$
 on remplace cette valeur dans la formule de récurrence et  on  a le
résultat.
\end{enumerate} \normalsize$\triangle$
\subsection*{Exemples : {\bf EQUATIONS DIFFERENTIELLES 1}}
1) $X\sim {\cal{N}}(0,1)$ :
$$ nH_n(X) - XH_{n}'(X)+H_{n}''(X)=0.$$
2) $X\sim \exp(\lambda)$ $(\lambda>0)$ :
$$nE_n^{\lambda}(X) = (X-\frac{ 1               }
{ \lambda}   ) (\ E_{n}^{\lambda}\ )'(X) - \sum_{k=2}^{n}
 \frac{ (\ E_n^{\lambda}\ )^{(k)}(X) }{\lambda^{k} }$$
3) $X\sim \Gamma( a , b )$ :
$$n \Gamma_n^{ab}(X) = (X-\frac{a}{b})(\ \Gamma_{n}^{ab}\ )'(X) - a
\sum_{k=2}^{n}(\  \Gamma_n^{ab}(X) \ )^{(k)} \biggl(\frac{1}{b}\ \biggl)^{k}.$$
4) $X\sim \Gamma( 1/2 , 1/2 )$ :
$$n \Gamma_n^{ \frac{1}{2} \frac{1}{2}  }(X) = (X-1)
(\ \Gamma_{n}^{\frac{1}{2}\frac{1}{2} } \ )'(X)  - \sum_{k=2}^{n}2^{k-1}
(\Gamma_n^{\frac{1}{2}\frac{1}{2} } \ )^{ (k) }(X).$$
5) $X \sim\Gamma\Bigl(\alpha{a_1\choose b_1},\beta{a_2\choose b_2}\Bigl) $ :
$$nQ_n(X) = \biggl( X - \alpha\frac{  a_1 }{ b_1 }- \beta\frac{ a_2 }{ b_2 } \biggl)Q_{n}'(X)   -
\sum_{k=2}^{n}Q_n^{(k)}(X) \biggl[ a_1\biggl(\frac{
\alpha }{b_1}\biggl)^{k} +a_2\biggl(\frac{\beta }{ b_2}\biggl)^{k}\biggl],$$
6) $ X\sim {\cal{P}}(a)$, $a>0$ :
$$n P_n^a(X) = (\ X - a\ )( P_{n}^a)'(X) - a\sum_{k=2}^{n}\frac{ (\ P_n^a\ )^{(k)} }{ k!}(X)$$
\subsection*{ii) Seconde équation différentielle}
\begin{thm}
Les polynômes de Wick vérifient une équation différentielle donnée par :
$$nW_n(X)  - \dps\sum_{k=1}^{n}  \frac{kXm_{k-1}-nm_{k} }{k!}W_n^{(k)}(X) = 0.$$
\end{thm}
\subsection*{\underline{\textbf{Démonstration}}}
\begin{enumerate}
\item[]
\small
C'est exactement la même chose que précédemment.
\end{enumerate} \normalsize$\triangle$
\subsection*{Remarque}
Cette équation différentielle n'est pas une réécriture de la première
(confère le premier exemple). On a ainsi trouver deux équations différentielles vérifiées
par les polynômes de Wick.
\subsection*{Exemples : \bf{EQUATIONS DIFFERENTIELLES 2}}
1) $X\sim {\cal{N}}(0,1)$ :
$$ nH_n(X) + n\sum_{k=1}^{[n/2]}\frac{H_n^{(2k)}(X) }{ 2^kk!} - X\sum_{k=0}^{[(n-1)/2]}
\frac{H_n^{(2k+1)}(X) }{ 2^kk! } =0.$$
2) $X \sim \exp(\lambda)$ $(\lambda>0)$ :
$$n E_{n}^{\lambda}(X) - \sum_{k=1}^n(\ E^{\lambda}_n\ )^{(k)} \biggl(\ \frac{ \lambda X-n}
{ \lambda^k } \ \biggl) =0$$
3) $X\sim\Gamma(a,b)$ :
$$n\Gamma_{n}^{ab}(X) - \sum_{k=1}^n\frac{ a(a+1)\cdot\ldots\cdot(a+k-2) }{ k!b^{k-1} }
(\ kX - \frac{n}{b}(a+k-1)\ ) (\ \Gamma_{n}^{ab}(X)\ )^{(k)}(X)=0$$
4) $X\sim \Gamma( 1/2 , 1/2 )$ :
$$n\Gamma_{n}^{ \frac{1}{2}\frac{1}{2} }(X)  - \sum_{k=1}^n\frac{ kX -2n(k-\frac{1}{2} ) }
{ k2^{k-1} } { 2k-2\choose k-1}\
(\  \Gamma_{n}^{\frac{1}{2}\frac{1}{2} }\ )^{(k)}(X) = 0.$$
5) $X\sim\Gamma\Bigl(\alpha{a_1\choose b_1},\beta{a_2\choose b_2}\Bigl) $ :\\
\hspace*{4cm}$nQ_n(X) = \dps\sum_{k=0}^nQ_{n}^{(K)}(X)\ \cdot\Biggl\{$\\
\hspace*{0,8cm}$\sum_{p=0}^{k-1}\Biggl[{ a_2+p-1\choose p}{ a_2+k-p-2\choose k-p-1}
\biggl( \frac{ \alpha }{ b_1}
\biggl)^p\biggl(\frac{\beta}{b_2}\biggl)^{k-1-p}\cdot\biggl( X -n\Bigl( \frac{a_2+k-p-1}{ k-p }
\Bigl)\ \biggl)\Biggl]
-n{ a_1+k-1\choose k}\biggl(\frac{\alpha}{b_1}\biggl)^k\Biggl\},$\\
6) $ X\sim {\cal{P}}(a)$ :
$$ nP_n^a(X) = e^{-a}\sum_{k=1}^n\biggl(\ kX (aD)^{k-1}\ (e^a) - n(aD)^{k}\ (e^a)\ \biggl)
\frac{ (\ P_{n}^a\ )^{(k)}(X) }{ k! }.$$
\section{Quelques lois en probabilités}
\subsection{Loi NORMALE, ${\cal{N}}(0,1)$}
$\begin{array}{lcl}
\mbox{- Densité }                                   & : &
                   f(x) = \dps\frac{ e^{\frac{-x^2}{2}} }{ \sqrt{2\pi} },\\
                                                                &\\
\mbox{- transformée de Laplace }          & : &                      \phi_X(t) = e^{\frac{t^2}{2} }.
\end{array}$
\subsection{Loi GAMMA de paramètres $a>0$ et $b>0$, $\Gamma(a,b)$ }
$\begin{array}{lcl}
\mbox{- Densité }                                 & : &
                       f(x)  = \dps1_{]0,\infty[}(x)\dps\frac{b^a}
{\Gamma(a) } e^{-bx} x^{a-1}\mbox{, }a>0\mbox{, } b>0,\\
                                                              &\\
\mbox{- transformée de Laplace  }       & : &
                      \phi_X(t) = (\dps\frac{b}{b-t} ) ^a\mbox{ pour } |t|<b.
\end{array}$
\subsection{Combinaison linéaire de Lois GAMMA, $\Gamma\Bigl(\alpha{a_1\choose b_1},
\beta{a_2\choose b_2}\Bigl) $  }
Les paramètres étant tous des réels strictement positifs, c'est la loi d'une
 \va Z pouvant s'écrire $Z = \alpha X +\beta Y$ où X et Y sont des \vas
 indépendantes de loi
respective $\Gamma(a_1,b_1)$ et $\Gamma(a_2, b_2)$.\\
Notons $f_{1\alpha}$ la densité de $\alpha X$ (qui suit une loi $\Gamma(a_1, \frac{b_1}
{\alpha} )$) et $f_{2\beta}$ la densité de $\beta Y$ ( qui suit une
loi $\Gamma(a_2,\frac{b_2}{\beta})$). \\
$\begin{array}{lcl}
\mbox{- Densité }                                   & : &
                         \dps f(x)  =\dps f_{1\alpha}*f_{2\beta}(x),\\
                                                              &\\
\mbox{- transformée de Laplace  }         & : &
                           \dps \phi_X(t) = \biggl( \dps\frac{b_1}{b_1-\alpha t}\biggl)^{a_1}
\biggl(\dps \frac{b_2}{b_2-\beta
t}\biggl)^{a_2}\mbox{ pour } |t|<b_1\land b_2.
\end{array}$
\subsection{Loi  EXPONENTIELLE de paramètre $\lambda$, $\E (\lambda)$ }
C'est une loi $\Gamma(1,\lambda)$.\\
$\begin{array}{lcl}
\mbox{- Densité }                                      & : &
             \dps f(x) = \lambda e^{-\lambda x} 1_{]0,\infty[}(x),\\
                                                                  &\\
\mbox{- transformée de Laplace }              & : &
              \dps \phi_X(t) =\frac{\lambda}{ \lambda-t}\mbox{ pour } |t|<\lambda.
\end{array}$
\subsection{Loi de POISSON de paramètre $\lambda$, ${\cal{P}}(\lambda$)}
$\begin{array}{lcl}
\mbox{-  Loi (discrète) }                    & : &
                 \dps     P(X =k ) = e^{-\lambda}\frac{\lambda^{k} }{ k!},\\
                                                                  &\\
\mbox{- transformée de Laplace }      & : &
                  \dps   \phi_X(t) = \exp( \ \lambda(e^{t}-1)\ ).
\end{array}$
\subsection{Loi BINOMIALE de paramètres N,p,q : ${\cal{B}} (N;p,q)$}
C'est la loi d'une \va sur $\{0,\ldots ,N\}$, avec $p+q =1$ ($p>0,q>0,N\in\N$) :\\
$\begin{array}{lcl}
\mbox{- loi  }                                               & : &
           \dps   P(X = k) = {N \choose k}p^kq^{N-k},\\
%
%
\mbox{- transformée de Laplace  }               & : &            \dps  \phi_X(t) = ( e^tp + q)^N.
\end{array}$
\section{Quelques rappels sur les polynômes de LAGUERRE généralisés,
$ L_n^{\alpha}(x)$ }
\underline{Intervalle d'étude} :
$$ [0,+\infty) .$$
\underline{Standardisation} :
$$[x^n] L_n^{\alpha}(x)= \frac{(-1)^n}{n!} .$$
\underline{Fonction poids} :
$$ w(x) = x^{\alpha}e^{-x}\mbox{ , }\alpha>-1 .$$
\underline{Norme} :
$$ h_n =\frac{ \Gamma(n+\alpha+1) }{ n!}  .$$
\underline{Formulation explicite} :
$$L_n^{\alpha}(x) = \sum_{m=0}^n(-1)^m{n+\alpha\choose n-m}\frac{ 1}{m!}x^m. $$
\underline{Formule de récurrence }:
$$\left \{
\begin{array}{lcl}
 (n+1)L_{n+1}^{\alpha}(x)  & = &[(2n+\alpha+1)-x]L_n^{\alpha}(x)
-(n+\alpha)L_{n-1}^{\alpha}(x) \\
 L_0^{\alpha}(x)  & = &  1\\
  L_1^{\alpha}(x) & = &1+\alpha-x  .
\end{array}
\right.$$
\underline{Equation différentielle } :
$$xy''+(\alpha+1-x)y'+ny=0\mbox{ , }y=L_n^{\alpha}(x)  .$$
\underline{Formule de Rodriguez} :
$$ L_n^{\alpha}(x) = \frac{ 1 }{ n!x^{\alpha}e^{-x} } \frac{ d^n  }{  dx^n  } \{x^{n+\alpha}e^{-x}\}.$$
\underline{Fonction génératrice} :
$$ F_{\alpha}(t,x) =(1-z)^{-\alpha-1}\exp(\frac{xz}{ z-1} )  =
\dps  \sum_{n\geq0}L_n^{\alpha}(x)z^n\mbox{ , }|z|<1 .$$
\section{Une relation avec les polynômes de Laguerre}
On va noter $\Gamma_{a,b}(t,x)$ la fonction génératrice des polynômes de Wick,
$\Gamma^{a,b}_n(x)$, associée à une \va suivant une loi $\Gamma(a,b)$.
\begin{pp}
$$\Gamma_{a,b}(t,x) = F_{a-1}(\frac{t}{t-b}, bx).$$
\end{pp}
\subsection*{\underline{\textbf{Démonstration}}}
\begin{enumerate}
\item[]
\small
C'est une simple vérification.
\end{enumerate} \normalsize$\triangle$
\begin{pp}
$$\left \{ \begin{array}{lcl}
                                             \Gamma_0^{a,b}(x) & = & L_0^{a-1}(bx)\\
                                             \Gamma_n^{a,b}(x) & = &\dps \frac{n!}{b^n}\dps
 \sum_{k=0}^{n-1}(-1)^{n-k}{ n-1\choose k}L_{n-k}^{a-1}(bx)\mbox{, }n\geq1.
                 \end{array}
    \right.$$
\end{pp}
\subsection*{\underline{\textbf{Démonstration}}}
\begin{enumerate}
\item[]
\small
$ F_{a-1}(\frac{t}{t-b}, bx) = \dps\sum_{k\geq0}L_{n}^{a-1}(bx)\biggl (\ \frac{t}{t-b} \ \biggl)^n$.\\
Pour $n\geq1$ $(\ \frac{1}{1-x} \ )^n = \dps\sum_{k\geq0}{k+n-1 \choose k}x^k$ donc :\\
$ F_{a-1}(\frac{t}{t-b}, bx) = \dps\sum_{k\geq1}\biggl[\ L_{n}^{a-1}(bx)(\ \frac{ -t }{ b }\ )^n
\sum_{k\geq0}{k+n-1 \choose k}(\ \frac{ t }{ b }\ )^n \ \biggl]+
L_{0}^{a-1}(bx)$\\
$ F_{a-1}(\frac{t}{t-b}, bx) = \dps\sum_{m\geq1}\dps\sum_{k=0}^{m-1}(-1)^{m-k}L_{m-k}^{a-1}
(bx){m-1 \choose k}(\ \frac{ t }{ b }\ )^m+ L_{0}^{a-1}(bx)$\\
Ou encore :
$$
\sum_{n\geq0} \Gamma_n^{a,b}(x)\frac{t^n}{n!} =  \dps\sum_{m\geq1}\dps
\sum_{k=0}^{m-1}\frac{ (-1)^{m-k} }{ b^n } L_{m-k}^{a-1}(bx){m-1 \choose k} m! \frac{
t^m }{ m! }+ L_{0}^{a-1}(bx)$$
d'où le résultat. %
\end{enumerate} \normalsize$\triangle$
\section{Les premiers polynômes de Wick}
1) $X\sim{\cal{N}}(0,1)$ :  $X^{ \widetilde{n} } = H_n(X) =\dps n!\sum_{k=0}^{[n/2]}
 \frac{ (-1)^kX^{n-2k} }{ 2^kk!(n-2k)! }$
$$\begin{array}{lcllcl}
                  H_0(x)        & = &        1             &      H_3(x)     & = &     x^3-3x      \\
                  H_1(x)        & = &        x             &      H_4(x)     & = &     x^4-6x^2+3    \\
                  H_2(x)        & = &         x^2-1      &      H_5(x)    & = &       x^5-10x^3+15x.

\end{array}$$
2) $X\sim$exp($\lambda$) :  $X^{ \widetilde{n} } = E^{\lambda}_n(X) = X^n - \dps
\frac{n}{\lambda} X^{n-1}$
$$\begin{array}{lcllcl}
                  E^{\lambda}_0(x)         &=&     1                                         &         E^{\lambda}_3(x)      &=&     x^3-\dps\frac{3}{\lambda}x^2\\
                    &\\
                  E^{\lambda}_1(x)         &=&     x-\dps\frac{1}{\lambda}      &        E^{\lambda}_3(x)      &=&     x^4-\dps\frac{4}{\lambda}x^3\\
                    &\\
                  E^{\lambda}_2(x)        &=&      x^2-\dps\frac{2}{\lambda}    &       E^{\lambda}_3(x)       &=&    x^5-\dps\frac{5}{\lambda}x^4.

\end{array}$$
3) $X\sim\Gamma(a,b)$ :  $X^{ \widetilde{n} } = \Gamma^{ab}_n(X) =\dps n!
\sum_{k=0}^n{a \choose n-k} \frac{X^kÊ}{Êk! }\Bigl( \frac{ -1 }{b} \Bigl)^{n-k}$
$$\begin{array}{lcl}
                \Gamma^{ab}_0(x)           &=&    1                                                                                                                                                        \\
                  &\\
                 \Gamma^{ab}_1(x)          &=&     x-\dps\frac{a}{b}                                                                                                                                        \\
                   &\\
                  \Gamma^{ab}_2(x)         &=&      x^2-2 \dps\frac{a}{b}x +\dps\frac{a(a-1)}{b^2}                                                                                                   \\
                     &\\
                  \Gamma^{ab}_3(x)         &=&  x^3-3 \dps\frac{a}{b}x^2+3\dps\frac{a(a-1)}{b^2}x-\dps\frac{a(a-1)(a-2)}{b^3}\\
                     &\\
                 \Gamma^{ab}_4(x)         &=& x^4-4\dps\frac{a}{b}x^3+6\dps\frac{a(a-1)}{b^2}x^2-4\dps\frac{a(a-1)(a-2)}{b^3}x  +   \dps\frac{a(a-1)(a-2)(a-3)}{b^4}         \\
                   &\\
                 \Gamma^{ab}_5(x)         &=&x^5-5\dps\frac{a}{b}x^4+10\dps\frac{a(a-1)}{b^2}x^3-10\dps\frac{a(a-1)(a-2)}{b^3}x^2  +   5\dps\frac{a(a-1)(a-2)(a-3)}{b^4} x\\
                                                      &  &   -\dps\frac{a(a-1)(a-2)(a-3)(a-4)}{b^5}.
\end{array}$$
4) $ X\sim \Gamma(1/2,1/2)$ : $X^{ \widetilde{n} } =\Gamma_n^{\frac{1}{2} \frac{1}{2}}(X) = X^n - \dps\frac{ n! }{ 2^{n-1} }\dps\sum_{k=0}^{n-1}\frac{ 2^k }{ n-k } { 2(n-k-1)
\choose n-k-1} \frac{X^k}{ k!  }$
$$\begin{array}{lcllcl}
    \Gamma_0^{\frac{1}{2} \frac{1}{2}}(x)           &= &        1                                                &
    \Gamma_3^{\frac{1}{2} \frac{1}{2}}(x)           &= &        x^3-3x^2-3x-3                            \\
                                                                           &\\
    \Gamma_1^{\frac{1}{2} \frac{1}{2}}(x)           &=&        x-1                                                   &
\Gamma_4^{\frac{1}{2} \frac{1}{2}}(x)            &= &       x^4-4x^3-6x^2-12x-15                   \\
                                                                           &\\
     \Gamma_2^{\frac{1}{2} \frac{1}{2}}(x)          &=&        x^2-2x-1                                            &
     \Gamma_5^{\frac{1}{2} \frac{1}{2}}(x)          &= &        x^5-5x^4-10x^3-30x^2-75x-105.

\end{array}$$
5) $X\sim {\cal{P}}(a)$ : $X^{ \widetilde{n} } = P^a_n(X) = \dps\sum_{k=0}^n{ n \choose k} X^{n-k} \Bigl[\ e^{a} \bigl( aD \bigl) ^k\cdot \bigl( e^{-a}\bigl) \Bigl]$
$$\begin{array}{lcl}
P^a_0(x)      &=&    1                                                                                                                                                            \\
P^a_1(x)      &=&      x-a                                                                                                                                                        \\
P^a_2(x)      &=&      x^2-2ax+(a^2-a)                                                                                                                                      \\
P^a_3(x)      &=&      x^3-3ax^2+3(a^2-a)x +(-a^3+3a-a)                                                                                                           \\
P^a_4(x)      &=&      x^4-4ax^3+6(a^2-a)x^2 +4(-a^3+3a-a)x +(a^4-6a^3+7a^2-a)                                                                       \\
P^a_5(x)      &=&      x^5-5ax^4+10(a^2-a)x^3 +10(-a^3+3a-a)x^2 +5(a^4-6a^3+7a^2-a)x                                                              \\
                   &   &    +( -a^5+10a^4-25a^3+15a^2-a).
\end{array}$$

\newpage
%
\chapter{Densité de \vas discrètes}

\section{Introduction}
Les fonctions de Rademacher, $r_0 = 1$, $n\geq1$
$r_n = \dps\sum_{k=0}^{2^n-1}  (-1)^k
1_{  ]\frac{ k }{Ê2^n } ,  \frac{ k+1 }{2^n } ]  }$, forment un système
souvent utilisé et dont les propriétés sont bien connues. Ce sont notamment
des \vas  de $\Bigl( [0,1], B([0,1]), \lambda\Bigl)$ à valeurs dans
\{-1,1\}, centrées, indépendantes.\\
On sait, si on pose ${\cal{F}}_n = \sigma(r_1,\ldots,r_n)$, que
$B([0,1]) = \dps\bigvee_n {\cal{F}}_n$ et construire à partir des
$(r_n)_n$ un sous-espace
vectoriel (sev) dense de $L^2\Bigl( [0,1], B([0,1]), \lambda\Bigl)$.\\
Cet \ev  est celui engendré par les produits $\dps\prod_{i=1}^{n}r_{k_i}$ où
$1\leq k_1<\cdots<k_n$ , $n\in\N$. (On appelle ce système,
"système de Walsh" pour J.L.Walsh).\\
Avant de généraliser ce résultat, voyons deux propositions faciles.
\begin{pp}
Etant donnée une suite  de \vas indépendantes centrées (vaic) du second ordre, on ne peut pas construire (sans en prendre les produits) un sev dense de $L^2$.
\end{pp}
\subsection*{\underline{\textbf{Démonstration}}}
\begin{enumerate}
\item[]
\small
soit $X_n : \om\to(E,{\cal{E}})$, comme $E\bigl(\ (X_0X_1)X_j\ \bigl) = 0$
$\forall j\geq0$ on a trouvé une \va $\not\equiv 0$ orthogonale à
 l'\ev engendré par les
$X_j, j\in\N$, il n'est donc pas dense dans $L^2\om$.
\end{enumerate} \normalsize$\triangle$
\subsection*{Remarque}
Une suite de vaic réduites fournit un système orthonormal ce qui pourrait
nous pousser à essayer de montrer que c'est une base de $L^2\om$. Cette
proposition nous en dissuade immédiatement.
\begin{pp}
Soient X : ($\Omega,{\cal{F}}_{X}, P) \to \Bigl( \R,B(\R), P_X\Bigl)$,
$(Q_n)_n$ une famille de polynômes de $L^2(\R)$ formant un sev dense.\\
Alors $\Bigl( Q_n(X)\Bigl)_n,$ forme un sev dense de $L^2(\Omega)$.
\end{pp}
\subsection*{\underline{\textbf{Démonstration}}}
\begin{enumerate}
\item[]
\small
$f\in L^2(\Omega)$  telle que $E\Bigl( fQ_n(X)\Bigl) =0$ $\forall n$.
f est ${\cal{F}}_{X}$-mesurable donc $f=g\circ X$ $(g : \R\to\R$).
La condition s'écrit :
$$\left(\int_{\R}g(x)Q_n(x)dP_X(x) = 0\  \forall
n\right)
\Rightarrow
\Bigl( g\equiv0\Bigl)
\Rightarrow
\Bigl( f\equiv0\Bigl).$$
\end{enumerate} \normalsize$\triangle$
\subsection*{Remarque}
Cette proposition nous dit qu'un moyen d'investigation, dans la recherche
de \sevs denses, est d'étudier des familles de polynômes dans les espaces
image.
%
%
%
%
%
%
%
%
%
\section{Suite de \vas prenant 2 valeurs}
On se donne une suite de vaic réduites (vaicr) $(X_n)_{n\geq1}$ :
$\om\to\{ a_1 , a_2 \}$
où ${\cal{F}} =\dps\bigvee_n{\cal{F}}_n$ avec  ${\cal{F}}_n
=\sigma\bigl(X_1,\ldots,X_n\bigl)$. Posons $\left\{\begin{array}{lcl}
     H_0^n &  =  &   [1]\\
     H_k^n   &  =  &\bigl [ X_{i_1}\cdots X_{i_k} \  :
 \ 1 \leq  i_1<\cdots <i_k\leq n \bigl ] \ Ê1\leq k \leq n
 \end{array}\right.$\\
$ {\cal{H}}_n = \dps\bigcup_{k\geq0}H_k^n$\\
$A_k^j = \bigl(X_j\bigl)^{-1}(a_k)$\\
$A_k = A_{k_1}^1\cap\cdots\cap A_{k_n}^n$ avec $k=(k_1,\ldots,k_n)\in
\{1,2\}^n$.\\
Les $H_k^n$ sont les \evs engendrés par tous les monônes de degré k qu'on
peut obtenir avec les n premiers termes de la suite. (Il est assez facile
de voir que pour tout j, $X_j^k$ s'exprime avec 1 et $X_j$ pour tout $k\geq2$.)
\begin{thm}
$\dps\bigcup_{n\geq0}{\cal{H}}_n$  est dense dans $L^2\om$.
\end{thm}
\begin{lm}
$$H_k^n\bot H_{k'}^n\ \ \ k\not=k'.$$
\end{lm}
\subsection*{\underline{\textbf{Démonstration du lemme 4.2.1.}}}
\begin{enumerate}
\item[]
\small
Posons $E_{kk'} = E \Bigl(X_{i_1}\cdots X_{i_k}X_{j_1}\cdots X_{j_{k'}}\Bigl)$, supposons $k<k'$, il existe une permutation $\sigma$ de
$\Sigma_{k'}$ telle que : \\
$$E_{kk'} = E \Bigl(X_{i_1}\cdots X_{i_k}X_{\sigma(j_1)}\cdots X_{\sigma(j_{k})}\Bigl)
E\Bigl(   X_{\sigma(j_{k+1})}   \Bigl)    \cdots    E\Bigl(X_{\sigma(j_{k'})}\Bigl)$$
Les \vas étant centrées, $E_{kk'} =0$.
\end{enumerate} \normalsize$\triangle$
\begin{lm}
$$\forall n\geq 1\  \ L^2(\Omega,\F_n, P) = {\cal{H}}_n.$$
\end{lm}
\subsection*{\underline{\textbf{Démonstration du lemme 4.2.2.}}}
\begin{enumerate}
\item[]
\small
D'une part  $\F_n =\sigma (X_1,\ldots,X_n) =
\sigma\bigl(A_k, \ k\in\{1,2\}^n\bigl)$, par indépendance des \vas
aucun des $A_k$ n'est vide, les systèmes $\Bigl\{A_1^j,A_2^j\Bigl\}$
$j\in\{1,\ldots,n\}$ sont des partitions de
$\Omega$. Il en est de même du système $\Bigl\{\ A_k, \ k\in \{1,2\}
\Bigl\}$ (on le voit facilement par récurrence).\\
$\F_n$ est donc engendrée par $2^n$ éléments qui forment une partition de
$\Omega$ et , par suite,\\
dim$\Bigl(\L^2(\Omega,\F_n,P)\Bigl)  = 2^n$  (engendré par les $1_{A_k})$.\\
D'autre part les éléments de $H_k^n$ sont linéairement indépendants, en
effet :\\
soit
$\dps\sum_{ 1\leq i_1<\cdots<i_k\leq n}
\lambda_{i_1\cdots i_k}X_{i_1}\cdots X_{i_k} = 0$ en multipliant par
$X_{j_1}\cdots X_{j_k}$ puis en prenant
l'espérance on trouve $\lambda_{j_1\cdots j_k}=0$.\\
Chaque $H_k^n$ possède ${ n\choose k}$ éléments (qui sont libres), de part
l'orthogonalité, ${\cal{H}}_n$ possède $\dps\sum_{k=0}^n{ n\choose k} = 2^n$
éléments.\\
Finalement, comme ${\cal{H}}_n \subset L^2(\Omega,\F_n, P)$ le lemme est
démontré. %
\end{enumerate} \normalsize$\triangle$
%
%
%
%
%
%
%
%
%
%
%
%
%
%
\subsection*{\underline{\textbf{Démonstration du théorème 4.2.1.}}}
\begin{enumerate}
\item[]
\small
Prenons f dans $L^2\om$ orthogonale ˆ $\bigcup{\cal{H}}_n$. Par le
lemme 4.2.2., f est orthogonale à chaque $L^2(\Omega,\F_n,P)$, $n\geq1$, donc
pour tout F
dans $\dps\bigcup_{n\geq1}{\cal{F}}_n $, $E \Bigl( 1_Ff\Bigl) = 0$,
mais $\Bigl\{\  F,\ E \bigl( 1_Ff\bigl) = 0\ \Bigl\}$ est une classe
monotone contenant
algèbre $\dps\bigcup_{n\geq1}{\cal{F}}_n $ donc, par le théorème de la
classe monotone, $E(1_Ff)= 0$ pour tout F dans $\F$ et par suite $f\equiv 0$ d'où le
résultat.
\end{enumerate} \normalsize$\triangle$
\section{Suite de \vas prenant un nombre fini de valeurs}
Soit $X_n : \om \to \bigl\{ a_1,\ldots,a_N\bigl\}$, $n\geq1$ une suite de
vaic ayant des moments de tous  ordres, $\F =\dps\bigvee_n\F_n$ où $\F_n =
\sigma ( X_1,\ldots X_n)$.
\begin{pp}
Pour tout $j\geq1$ fixé dans $\N$, le système
$\Bigl\{1,X_j,\ldots,X_j^{N-1}\Bigl\}$ est une base de
$L^2(\Omega,\F_{X_j},P). $
\end{pp}
\subsection*{Remarque}
Cela veut dire que dans un ensemble générateur on ne peut pas avoir de
puissances supérieures à N-1 de chacune des \vas.
\subsection*{\underline{\textbf{Démonstration}}}
\begin{enumerate}
\item[]
\small
Posons comme précédemment $A_k^j = \bigl( X_j\bigl)^{-1}(a_k)$\ \
$k\in\{1,\ldots ,N\}^n$.\\
Alors $X_j = a_11_{A_1^j}+\cdots+a_N1_{A_N^j}$.\\
D'une part $\Bigl\{1_{A_1^j},\ldots , 1_{A_N^j}\Bigl\}$ est une base de
$L^2(\Omega,\F_{X_j},P)$.\\
D'autre part si $X^0_j,\ldots,X^{N-1}_j$ est un système libre de
$L^2(\Omega,\F_{X_j},P)$ alors c'est une base.\\
Formons : $\dps\sum_{k=0}^{N-1}\lambda_kX_j^k=0$\\
$\Leftrightarrow \dps\sum_{k=0}^{N-1}
\lambda_k\bigl(a_11_{A_1^j}+\cdots+a_N1_{A_N^j}\bigl)=0
\Leftrightarrow\dps\sum_{l=0}^{N-1}1_{A_l^j}\bigl( \lambda_0+
\lambda_1a^1_l+\cdots+\lambda_{N-1}a_l^{N-1}\bigl)=0$\\
$\Leftrightarrow\left(\begin{array}{ccccc}
 1           & a_1 &  a_1^2 & \cdots  & a_1^{N-1}\\
     \vdots  &       &            &              & \vdots    \\
  1          &a_N  & a_N^2  & \cdots  & a_N^{N-1}
                       \end{array}\right)
  \left(\begin{array}{c}
    \lambda_0\\
      \vdots\\
        \lambda_{N-1}
   \end{array}\right)
=0\Leftrightarrow \lambda_0=\cdots=\lambda_{N-1}=0$
(déterminant de Vandermonde).
\end{enumerate} \normalsize$\triangle$
\subsection*{Remarque}
On a envie de faire la même construction que précédemment en posant
$H_k^n$ pour l'ensemble des monômes de degré k ce qu'on peut écrire grace
à la proposition 4.3.1. : \\
%
%
%
%
%
%
%
%
%
$H_k^n =\Bigl[\ X_{i_1}^{ \alpha_{i_1} } \cdots X_{i_n}^{ \alpha_{i_n} }\
\ 1\leq \alpha_{i_1}<\cdots< \alpha_{i_n}\leq n\ , \
\alpha_{i_1}+\cdots+\alpha_{i_n}
= k\ ,Ê\ \alpha_{i_{\nu} }\in\{0,\ldots ,N-1\} \ \Bigl]$\\
et en posant ${\cal{H}}_n = \dps\bigcup_{k=0}^{n(N-1)}H_k^n$.\\
Cette décomposition n'a plus de sens pour $N\geq3$ car les $H_k^n$ ne sont
plus orthogonaux.
\subsection*{Exemple, N=3}
n=1 $\left\{\begin{array}{lcl}
                                  H_0^1 & = & [1]\\
                                  H_1^1 & = & [X_1]\\
                                  H_2^1 & = & [X_1^2]
         \end{array}\right.$\\
$E(X_1\cdot X_1^2)=E(X_1^3)$\\
n=2  $\left\{\begin{array}{lcl}
                                  H_0^2 & = & [1]\\
                                  H_1^2 & = & [X_1,X_2]\\
                                  H_2^2 & = & [X_1^2,X_2^2,X_1X_2]\\
                                  H_3^2 & = & [X_1^2X_2,X_1X_2^2]\\
                                  H_4^2 & = & [X_1^2X_2^2]
\end{array}\right.$\\
$E(X_1\cdot X_1^2)=E(X_1^3)$; $E(X_1X_2\cdot X_1^2X_2)=E(X_1^3)E(X_2^2)$ \\
En prenant $X = -1 \cdot1_{]0,1]}+3\cdot1_{]1,2]}-2\cdot1_{]2,3]}$,
$P=\frac{\lambda}{3}$, E($X^3$) = 6
\begin{thm}
Posons :$$\left\{ \begin{array}{lcl}
{\cal{H}}_0 & = &[1]\\
{\cal{H}}_n & = & \Bigl[ X_1^{\alpha_1}\cdots X_n^{\alpha_n}\ , \
(\alpha_1,\ldots\alpha_n)\in\{0,\ldots N-1\}^n\ \Bigl] \ n\geq1
\end{array}\right.$$.\\
Alors $\dps\bigcup_{n\geq0}{\cal{H}}_n$ est dense dans $ L^2\om$.
\end{thm}
\begin{lm}
Le système $\Bigl\{ X_1^{\alpha_1}\cdots X_n^{\alpha_n}\ , \ (\alpha_1,\ldots,\alpha_n)\in\{0,\ldots ,N-1\}^n\ \Bigl\}$ est libre.
\end{lm}
\subsection*{\underline{\textbf{Démonstration du lemme 4.3.1.}}}
\begin{enumerate}
\item[]
\small
Formons $\dps \sum_{  (\alpha_1,\ldots,\alpha_n)
\atop { \in \{0,\ldots,N-1\}^n  }  }
\lambda_{\alpha_1\cdots\alpha_n} X_1^{\alpha_1}\cdots X_n^{\alpha_n} =
0$\\
Le système $\Bigl\{1_{A_k}  \   ,  k\in\{1,\ldots ,N\}^n\ \Bigl\}$ forme
un base de $L^2(\Omega,\F_n,P)$ (c'est le même raisonnement que dans le cas
de deux valeurs, confère lemme 4.2.2.).\\
On peut écrire : $X_k^{\alpha_k} = \dps\sum_{\nu_k =1}^N
a_{\nu_k}^{\alpha_k} 1_{A_{\nu_k}^k}$.\\
Si $0=(\alpha_1,\ldots,\alpha_n)$, $1=X_1^{\alpha_1}\cdots X_n^{\alpha_n} =
\dps \sum_{ (\nu_1,\ldots,\nu_n)  \atop{\in\{1,\ldots,N\}^n  }}
1_{A^1_{\nu_1}\cap\cdots\cap A^n_{\nu_n} } $ (partition).\\
Supposons maintenant que $0\not=
 (\alpha_1,\cdots,\alpha_n) \in\{0,\ldots,N-1\}^n$ :\\
$\exists\  (i_1,\ldots,i_k)\in\{1,\ldots,N-1\}^k$, $i_1<\ldots<i_k$ :\\
$$X_1^{\alpha_1}\cdots X_n^{\alpha_n}  =X_{i_1}^{\alpha_{i_1}}\cdots
X_{i_k}^{\alpha_{i_k}}  = \dps\sum_{ (\nu_1,\ldots,\nu_k)
\atop{\in\{1,\ldots,N\}^k } }
a^{\alpha_{i_1}}_{\nu_1}\cdots a^{\alpha_{i_k}}_{\nu_k}
 1_{ A^{i_1}_{\nu_1}\cap\cdots\cap A^{i_k}_{\nu_k} }$$
%
%
%
%
%
%
%
%
%
mais $1_{A^{i_1}_{\nu_1}\cap\cdots\cap A^{i_k}_{\nu_k}  }=
\dps\sum_{ (\nu_{k+1},\ldots,\nu_n)  \atop{\in\{1,\ldots,N\}^{n-k}}}
1_{  ( A^{i_1}_{\nu_1}\cap\cdots\cap A^{i_k}_{\nu_k}  ) \cap
( A^{i_{k+1}}_{\nu_{k+1}}\cap\cdots\cap A^{i_n}_{\nu_n}  )  }$
 (partition) \\
où $ i_{k+1}<\cdots<i_n$ et $ \{  i_{k+1},\ldots,i_n\} =
 \{1,\ldots, N\}- \{i_{1},\ldots,i_k\}$\\
On peut écrire : \\
$X_{i_1}^{\alpha_{i_1}}\ldots X_{i_k}^{\alpha_{i_k}}  =
\dps\sum_{ (\nu_1,\ldots,\nu_k)  \atop{\in\{1,\ldots,N\}^k}}
a^{\alpha_{i_1}}_{\nu_1}\cdots a^{\alpha_{i_k}}_{\nu_k}
\dps\sum_{ (\nu_{k+1},\ldots,\nu_n)  \atop{\in\{1,\ldots,N\}^{n-k}}}
1_{  A^{i_1}_{\nu_1}\cap\cdots\cap A^{i_n}_{\nu_n}  } $\\
\hphantom{$X_{i_1}^{\alpha_{i_1}}\cdots X_{i_k}^{\alpha_{i_k}}  $}
$= \dps\sum_{ (\nu_1,\ldots,\nu_n)  \atop{\in\{1,\ldots,N\}^n }}a^{\alpha_{i_1}}_{\nu_1}\cdots a^{\alpha_{i_k}}_{\nu_k}
1_{  A^{i_1}_{\nu_1}\cap\cdots\cap A^{i_n}_{\nu_n}  } $\\
$\exists\  \sigma \in \Sigma_n$ (permutations d'ordre n) telle que
$ A^{i_1}_{\nu_1}\cap\cdots\cap A^{i_n}_{\nu_n} =   A^{1}_{\nu_{\sigma(1)  }}\cap\cdots\cap A^{n}_{\nu_{\sigma(n) }}$\\
$X_{i_1}^{\alpha_{i_1}}\cdots X_{i_k}^{\alpha_{i_k}}  = \dps\sum_{ (\nu_1,\ldots,\nu_n)  \atop{\in\{1,\ldots,N\}^n }}
a^{\alpha_{i_1}}_{\nu_1}\cdots a^{\alpha_{i_k}}_{\nu_k}  1_{   A^{1}_{\nu_{\sigma(1) }  }\cap\cdots\cap A^{n}_{\nu_{\sigma(n) }}    }$\\
\hphantom{$X_{i_1}^{\alpha_{i_1}}\cdots X_{i_k}^{\alpha_{i_k}}  $}
$ =\dps\sum_{ (\nu_1,\ldots,\nu_n)  \atop{\in\{1,\ldots,N\}^n }}
a^{\alpha_{i_1}}_{\nu_{i_1}}\cdots a^{\alpha_{i_k}}_{\nu_{i_k}}
1_{A^1_{\nu_1}\cap\cdots\cap A^n_{\nu_n}  }$.\\
Comme $\alpha_{i_{k+1}}=\cdots=\alpha_{i_n} = 0$ on peut écrire :  \\
$X_1^{\alpha_1}\cdots X_n^{\alpha_n}  = \dps\sum_{ (\nu_1,\ldots,\nu_n)
\atop{\in\{1,\ldots,N\}^n }}
a^{\alpha_{i_1}}_{\nu_1}\cdots a^{\alpha_{i_n}}_{\nu_n}
1_{A^1_{\nu_1}\cap\cdots\cap A^n_{\nu_n}  }$.\\
Donc : \\
$\dps \sum_{  (\alpha_1,\ldots,\alpha_n)   \atop { \in \{0,\ldots,N-1\}^n  }
 } \lambda_{\alpha_1\cdots\alpha_n} X_1^{\alpha_1}\cdots X_n^{\alpha_n} =0
\Leftrightarrow
\dps \sum_{ (\alpha_1,\ldots,\alpha_n) \atop{\in\{0,\ldots,N-1\}^n  }}
\lambda_{\alpha_1\cdots\alpha_n}
\dps\sum_{ (\nu_1,\ldots,\nu_n)  \atop{\in\{1,\ldots,N\}^n }}
a^{\alpha_{1}}_{\nu_1}\cdots a^{\alpha_n}_{\nu_n}
1_{A^1_{\nu_1}\cap\cdots\cap
A^n_{\nu_n}  } = 0$.\\
$\Leftrightarrow
\dps\sum_{ (\nu_1,\ldots,\nu_n)  \atop{\in\{1,\ldots,N\}^n }}
\Biggl(  \dps \sum_{ (\alpha_1,\ldots,\alpha_n) \atop{\in\{0,\ldots,N-1\}^n
 }}\lambda_{\alpha_1\cdots\alpha_n} a^{\alpha_{1}}_{\nu_1}\cdots
a^{\alpha_{n}}_{\nu_n}    \Biggl)  1_{A^1_{\nu_1}\cap\cdots
\cap A^n_{\nu_n}  } = 0$.\\
Comme le système $\Bigl\{ 1_{A_k}\ ,\ k\in\{1,\ldots,N\}^n\ \Bigl\}$ est
libre on a :
$$ \forall (\nu_1,\ldots,\nu_n)\in \{1,\ldots,N\}^n\ ,\
\dps \sum_{ (\alpha_1,\ldots,\alpha_n) \atop{\in\{0,\ldots,N-1\}^n  }}
\lambda_{\alpha_1\cdots\alpha_n}
a^{\alpha_{1}}_{\nu_1}\cdots a^{\alpha_{n}}_{\nu_n}   = 0.$$
Regardons le cas n=2 : $\forall  (\alpha_1,\alpha_2) \in\{0,\ldots,N-1\}^2
                                      \dps\sum_{ (\nu_1,\nu_2)  \atop{\in\{1,\ldots,N\}^2 }}    \lambda_{\alpha_1\alpha_2}a^{\alpha_{1}}_{\nu_1}
a^{\alpha_{2}}_{\nu_2} = 0$

En utilisant l'ordre lexicographique pour ranger le vecteur colonne $ \Bigl(  \lambda_{\alpha_1\alpha_2}  \Bigl)_{(\alpha_1\alpha_2)  }$, on
obtient l'écriture matricielle suivante :\\
$$\left(\begin{array}{cccc}
A_N & a_1A_N &\cdots & a_1^{N-1}A_N\\
A_N & a_2A_N &\cdots & a_2^{N-1}A_N\\
\vdots&             &            &  \vdots   \\
\vdots&             &            &  \vdots   \\
A_N & a_NA_N&  \cdots & a_N^{N-1}A_N
\end{array}\right)
\left(\begin{array}{c}
\lambda_{00}\\
\vdots\\
\lambda_{01}\\
\vdots\\
\lambda_{N-1N-1}
\end{array}\right) = 0.$$
Où on a posé $A_N =\left(\begin{array}{cccc}
1& a_1 &\cdots & a_1^{N-1}\\
\vdots&             &            &  \vdots   \\
1 & a_N&  \cdots & a_N^{N-1}
\end{array}\right)$\\
%
%
%
%
%
%
%
%
%
%
Ce qu'on peut réécrire :
$$\left(\begin{array}{ccc}
A_N  &                              &  \\
        &  0 \ddots 0             &     \\
       &                               &    A_N
\end{array}\right)
\left(\begin{array}{cccc}
I& a_1 I &\cdots & a_1^{N-1 }I\\
\vdots&             &            &  \vdots   \\
I & a_N I&  \cdots & a_N^{N-1} I
\end{array}\right)  = 0$$
où $I = Id_{N\mbox{\tiny{x}}N}$.\\
det $A_N\not=0$ (déterminant de Vandermonde), on peut appeler le déterminant
de la deuxième matrice "déterminant de Vandermonde généralisé"
(on le calcule
par la même méthode que le déterminant classique).\\
Prenons comme hypothèse de récurrence que le système s'écrit :\\
$\Bigl(\ M_n\ \Bigl) \Bigl(  \lambda_{\alpha_1\cdots\alpha_n} \Bigl) =
0$ avec det $M_n$ = 0 ( $M_n\in{\cal{M}}_{N^n\mbox{\tiny{x}}N^n}$).
Voyons le rang n+1.
On vérifie aisément qu'on a :  $\Bigl(\ M_{n+1}\ \Bigl) \Bigl(
 \lambda_{\alpha_1\cdots\alpha_{n+1} } \Bigl) = 0$ avec\\
$$M_{n+1}= \left(\begin{array}{cccc}
M_n& a_1M_n &\cdots & a_1^{N-1}M_n\\
\vdots&             &            &  \vdots   \\
M_n & a_N M_n&  \cdots & a_N^{N-1}M_n
\end{array}\right)$$
(en utilisant l'ordre lexicographique pour ranger le vecteur
$\lambda_{\alpha_1\cdots\alpha_{n+1}}$)\\
$M_{n+1} = \left(\begin{array}{ccc}
M_n  &                              &  \\
        &  0 \ddots 0             &     \\
       &                               &    M_n
\end{array}\right)
\left(\begin{array}{cccc}
I& a_1 I &\cdots & a_1^{N-1 }I\\
\vdots&             &            &  \vdots   \\
I & a_N I&  \cdots & a_N^{N-1} I
\end{array}\right)  = 0$
où $ I=  Id_{N^n\mbox{\tiny{x}}N^n}$.\\
On utilise l'hypothèse de récurrence pour déduire : det $M_{n+1}\not=0$. Par
 suite on peut conclure que $\lambda_{\alpha_1\cdots\alpha_{n+1}} = 0
$ $\forall(\alpha_1,\ldots,\alpha_{n+1})\in\{0,\ldots,N-1\}^{n+1}$, c'est à
dire que le système est libre, ce qui démontre le lemme 4.3.1.
\end{enumerate} \normalsize$\triangle$
\begin{lm}
$$\forall n\geq0\ L^2(\Omega,\F_n,P)  = {\cal{H}}_n.$$
\end{lm}
\subsection*{\underline{\textbf{Démonstration}}}
\begin{enumerate}
\item[]
\small
Comme ${\cal{F}}_n = \sigma\Bigl(\ A_k\ , \ k\in\{1,\ldots,N\}^n\ \Bigl)$,
$L^2(\Omega,\F_n,P)$ est engendré par $\Bigl\{ \ 1_{A_k}\ ,\ k
\in\{1,\ldots,N\}^n\
\Bigl\}$ (puisque $\{A_k,\  k\}$ est une partition de $\Omega$).\\
Donc dim$\Bigl(\ L^2(\Omega,\F_n,P)  \ \Bigl) = N^n$ (aucun des $A_k$ n'est
vide
par indépendance des \vas; comme aucun des $A_{k_j}^j$ n'est de probabilité
nulle, $p(A_k) \not= 0$ $\forall k \in\{1,\ldots,N\}^n)$.\\
On voit aisément que : card$\Bigl\{ X_1^{\alpha_1}\cdots X_n^{\alpha_n}\ ,\ (\alpha_1,\ldots,\alpha_n)\in\{0,\ldots, N-1\}^n\Bigl\} = N^n$\\
(on a une bijection
$\{0,\ldots,N-1\}^n \to {\cal{H}}_n$
$(\alpha_1,\ldots,\alpha_n)\to X_1^{\alpha_1}\cdots X_n^{\alpha_n}$).\\
Le lemme 4.3.1. nous permet d'en déduire : dim(${\cal{H}}_n) =N^n$\\
Ainsi $\left\{\begin{array}{l}
{\cal{H}}_n\subset L^2(\Omega,\F_n,P)  \\
{\cal{H}}_n  \mbox{\ et\ } L^2(\Omega,\F_n,P) \mbox{\ ont le même nombre de
générateurs\ }
\end{array}\right.$\\
 donc  ${\cal{H}}_n = L^2(\Omega,\F_n,P)  $.%
\end{enumerate} \normalsize$\triangle$
\subsection*{\underline{\textbf{Démonstration du théorème 4.3.1.}}}
\begin{enumerate}
\item[]
\small
Prenons f dans $L^2\om$ orthogonale à $\bigcup{\cal{H}}_n$. Par le lemme
4.3.2.
f est orthogonale à chaque $L^2(\Omega,\F_n,P)$, $n\geq1$, donc pour tout F
dans $\dps\bigcup_{n\geq1}{\cal{F}}_n $, $E \Bigl( 1_Ff\Bigl) = 0$, mais
$\Bigl\{\  F,\ E \bigl( 1_Ff\bigl) = 0\ \Bigl\}$ est une classe monotone
contenant
algèbre, $\dps\bigcup_{n\geq1}{\cal{F}}_n $ donc, par le théorème de la
classe monotone, $E(1_Ff)= 0$ pour tout F dans $\F$ et par suite $f\equiv 0$ d'où le
résultat.
\end{enumerate} \normalsize$\triangle$
%
%
%
%
%
%
%
%
%
%
%
%
%
%
\newpage
\chapter{Densité de produit de \vas continues}
\section{Densité}
\subsection{ Notations, premiers résultats}
Le but de ce paragraphe est d'étudier l'espace engendré par les produits
des termes d'une suite de \vas.\\
Pour cela, donnons nous une suite de variables aléatoires réelles,
centrée réduites
$$\Bigl( X_n\Bigl)_{n\in\N^*}$$
définies sur un espace probabilisé
$$ \Bigl( E,{\cal E},\nu \Bigl)$$
${\cal E}$ étant la tribu engendrée par la suite considérée
$${\cal E} :={\cal E}_{\infty}:=\sigma\Bigl( X_k|k\in\N^* \Bigl)$$
et E un espace vectoriel normé.\\
Nous noterons :
$${\cal E}_n:=\sigma\Bigl( X_k|k=1,\ldots,n \Bigl)$$
$${\cal E}_{X_k}:=\sigma\Bigl( X_k\Bigl)\mbox{\ où } k\in\N^*.$$
Il s'agit donc d'étudier l'espace des polynômes
$\R\Bigl[(X_n)_{n\in\N^*}\Bigl]$. Plus précisément, notre étude va
porter sur le $\R$-espace vectoriel engendré par les produits de la
forme :
$$X_1^{\alpha_1}\cdot\ldots\cdot X_n^{\alpha_n}\;\;\Bigl( n\in\N^*\Bigl)
\Bigl((\alpha_1,\ldots, \alpha_n)\in(\N^*)^n\Bigl).$$
Notons :
$${\cal X}_n:=Vect_{\R}\Bigl\{ X_1^{\alpha_1}\cdot\ldots\cdot
X_n^{\alpha_n} |(\alpha_1,\ldots, \alpha_n)\in(\N^*)^n \Bigl\}:={\cal
X}_n$$
$${\cal X}:=\bigcup_{n\in\N^*}{\cal X}_n$$
\begin{thm}
Supposons que
pour tout $n\in\N^*$, ${\cal X}_n$ soit dense dans $ L^2\Een,$\\
alors
$${\cal X}\mathop\subset_{dense}\limits
L^2\Ee.$$
\end{thm}
\subsection*{Remarque}
Les raisonnements donnant des résultats de densité s'inspirent souvent
du théorème que nous rappelons ici dans sa version
hilbertienne. Sa démonstration utilise le théorème
de Hahn-Banach sous sa deuxième forme géométrique (confère  {\bf [5]} p.7).
\begin{thm}
Soient H un espace de Hilbert et G un \sev de H.\\
Si pour tout $ h\in H$

$$
   \Bigl(  (\forall g\in G)(<h,g>=0)  \Bigl)
   \Longrightarrow
   \Bigl(h=0\Bigl)
$$
où $<,>$ désigne le produit scalaire de H, alors $\overline{G}=H$, ce
qu'on peut aussi écrire :
$$G^{\bot}=\{ 0 \}.$$
\end{thm}
\subsection*{\underline{\textbf{Démonstration}}}
\begin{enumerate}
\item[]
\small
Soit $f\in L^2\Ee\cap {\cal X}^{\bot}$.\\
La suite $\Bigl( E(f|{\cal E}_n) \Bigl)_{n\in\N^*}$ est une martingale
de carré intégrable (càd qui vérifie $\sup_nEM_n^2<\infty$) relativement
à la filtration
$\Bigl({\cal E}_n\Bigl)_{n\in\N^*}$, elle converge dans
$L^2\Ee$ vers $E(f|\dps\bigvee_n{\cal E}_n)=f$  (puisque
${\cal E}=\dps\bigvee_n{\cal E}_n$),
confère par exemple {\bf[6]} p.49 et p.62.\\
Fixons $n\in\N^*$, comme $f\in {\cal X}^{\bot}$ on a :
$$\Biggl(\forall (\alpha_1,\ldots,\alpha_n)\in\N^n\Biggl)
\Biggl(
E\Bigl(E(f|{\cal E}_n)X_1^{\alpha_1}\cdot\ldots\cdot
X_n^{\alpha_n}\Bigl)= E\Bigl(fX_1^{\alpha_1}\cdot\ldots\cdot
X_n^{\alpha_n} \Bigl)= 0
\Biggl) .$$
En particulier, n étant fixé dans $N^*$,
$$E(f|{\cal E}_n)\in L^2\Een\cap {\cal X}^{\bot},$$
et par densité de ${\cal X}_n$ dans $L^2\Een$
$$E(f|{\cal E}_n)=0.$$
ce qui entraîne  $f=0$ et par suite
$${\cal X}^{\bot} =\{0\}$$
i.e. ${\cal X}$ est dense dans $L^2\Ee$.\\
\end{enumerate} \normalsize$\triangle$
\subsection{Quelques conditions suffisantes}
\subsection*{i) Existence de moments exponentiels}
Nous allons voir que l'hypothèse d'existence de moments exponentiels
est une condition  suffisante de densité
(confère  {\bf [7]} chap. 2).
\begin{de}
On dira qu'une \va U admet des moments exponentiels si :
$$\Bigl(\exists t_0>0\Bigl)\Bigl(\forall |t|<t_0\Bigl)
\Bigl(E( e^{tU})<\infty\Bigl).$$
Nous noterons ${\cal E}xp(U)$ cette hypothèse.
\end{de}

\hspace*{-0,6cm}La proposition suivante donne le lien entre la condition $\Ex$ et la
densité des polynômes, commençons par rappeler une propriété de telles
\vas.
\begin{lm}
Toute \va X vérifiant  $\Ex$ admet des moments de tous ordres :
$$\Bigl(
\forall \alpha\in\N
\Bigl)
\Bigl(
E( |  X  |^{\alpha}   )   <  \infty
\Bigl)
$$
\end{lm}
\begin{pp}
Soit X une \va réelle définie sur $\Bigl( E,{\cal E}_X,\nu\Bigl)$
vérifiant
$\Ex$
alors
$${\cal X}:=Vect_{\R}\Bigl\{X^{\alpha} | \alpha\in\N\Bigl\}
\mathop{\subset}_{dense}\limits
L^2\Bigl( E,{\cal E}_X,\nu\Bigl).$$
\end{pp}
\subsection*{\underline{\textbf{Démonstration de la proposition 5.1.1.}}}
\begin{enumerate}
\item[]
\small
Soit $g\in L^2\Bigl(E,{\cal E}_X,\nu\Bigl)
\cap
{\cal X}^{\bot}$.\\
$\exists f\in L^2\Bigl( \R,B(\R),\nu_{ X}\Bigl)$
telle que $g=f\circ X$, où $\nu_{ X}$ est la loi image de X par $\nu$.\\
Posons :
$$\varphi(t):= E^{\nu}\Bigl(e^{itX}g\Bigl)$$
comme $g\in L^2\Bigl(E,{\cal E}_X,\nu\Bigl)$, $\varphi(t)$ est bien
définie, pour tout $t\in\R$. En passant dans l'espace image on a :
$$ \varphi(t) = \int_{\R} e^{itx}f(x)d\nu_X(x).$$
$\bullet$ \underline{Régularité de $\varphi$}\\
Soit $\alpha\in\N$
$$\begin{array}{rcl}
%
\varphi_{\alpha}(t)    &  :=      &
\dps\int_{\R}
|x^{\alpha}| \cdot |e^{itx}f(x)| d\nu_X(x)\\
%
                                 &  \mathop{\leq}_{CS}\limits       &
||f||_2
\Biggl(
 \dps\int_{\R}  |x^{2\alpha}| d\nu_X(x)
\Biggl)^{\frac{1}{2}}\\
%
                                 &\leq   &+\infty.
\end{array}$$
La dernière inégalité étant vraie par le lemme 5.1.1. On a noté
CS pour exprimer qu'on utilise l'inégalité de Cauchy-Schwartz. Le
théorème de dérivation de Lebesgue assure l'existence de toutes les
dérivées de $\varphi$.\\
$\bullet$ \underline{Existence d'un ouvert où $\varphi$ et toutes ses
dérivées sont nulles}\\
f étant orthogonale aux polynômes, on va développer $e^{itx}$
en série entière de façon à faire apparaître les monômes
$x^{\alpha}$ ($\alpha\in\N$).\\
Il faut qu'on puisse inverser $\int$ et $\sum$, vérifions le.\\
Regardons :
$$\begin{array}{rcl}
\psi (t)    &:= &
\dps\int_{\R}
\sum_{k\in\N} |\frac{(itx)^k}{k!}|f(x)|  d\nu_X(x)\\
            & = &
\dps\int_{\R}  e^{|tx|}|f(x)|  d\nu_X(x)\\
%
            & \mathop{\leq}_{ \mbox{CS} }\limits &
|| f ||_2
\Biggl(
\dps\int_{\R}  e^{2|tx|}  d\nu_X(x)
\Biggl)^{ \frac{1}{2} }\\
%
%
            & < & +\infty\  \mbox{ pour  }
\dps|t|<\frac{t_0}{2}\\
\end{array}$$
On peut intégrer terme à terme la série
$\dps\sum_{k\in\N}\frac{(itx)^k}{k!}$  sur l'intervalle
$D:= ]-\frac{t_0}{2},+\frac{t_0}{2}[$ et
$$
\varphi(t)
= \sum_{k\in\N}\frac{(it)^k}{k!}
            \int_{\R}  x^kf(x) d\nu_X(x)
= \sum_{k\in\N}\frac{(it)^k}{k!}
            E^{\nu}\Bigl( X^kg\Bigl)
$$
de la même façon, pour $j\in\N^*$ :
$$
\varphi^{(j)}(t)
= \sum_{k\in\N}\frac{t^k}{k!}
            \int_{\R}  (ix)^{k+j}f(x) d\nu_X(x)
= \sum_{k\in\N}\frac{t^ki^{k+j}}{k!}
            E^{\nu}\Bigl( X^{k+j}g\Bigl)
$$
Comme $g\in{\cal X}^{\bot}$, $\varphi^{(j)}_{|_D}=0$ $(j\in\N)$.\\
$\bullet$ \underline{ Prolongement }\\
Posons :
\begin{enumerate}
\item[-]
$D_0:=D$
\item[-]
$D_p:=] -\frac{t_0}{2} + p\frac{t_0}{2},\frac{t_0}{2} + p\frac{t_0}{2}[$
$p\in\N^*$
\end{enumerate}
On a $\varphi_{|_{D_0}}=0$, c'est le point précédent.\\
Pour prolonger ce résultat à $\R$ faisons un raisonnement par
récurrence sur p.\\
Soient $j\in\N$ et $t\in D_{p+1}$, $\exists t_p\in D_p$ tel que $|t-t_p|<\frac{t_0}{2}$,
$$\varphi^{(j)}(t) = \dps\int_{\R}(ix)^j  e^{i(t-t_p)x}e^{it_px}f(x)
d\nu_X(x).$$
$\left.
\begin{array}{cc}
\bullet  & |(ix)^je^{it_px}f| =|(x)^jf|\\
\bullet  & |t-t_p|<t_0
\end{array}\right\}$ $\Longrightarrow$
$\left\{\begin{array}{ll}
                 & \mbox{on peut intégrer terme à terme}\\
                  & \mbox{la série } \dps\sum_{k\in\N}
                             \frac{(i(t-t_p)x))^k}{k!}(ix)^je^{it_px}f(x).
\end{array}\right.$\\
Ce qui nous donne
$$\varphi^{(j)}(t) =\sum_{k\in\N}
\frac{(t-t_p)^k}{k!}
\varphi^{(k+j)}
(t_p)=0$$
et donc $\varphi^{(j)}_{|_{\dps\cup_{p\in\N}D_p}}=0$.\\
D'autre part
$$\Bigl(\forall t\in\R\Bigl)
\Bigl(
(\exists p\in\N)(|t|<\frac{t_0}{2}+p\frac{t_0}{2}
=\frac{p+1}{2}t_0)\Bigl)$$
puisque $\dps\lim_{p\rightarrow+\infty}
(\frac{p+1}{2}t_0)) =+\infty$\\
donc $\dps\bigcup_{p\in\N}D_p= \R$ et finalement, pour tout$j\in\N$
$\varphi^{(j)}\equiv 0$.\\
$\bullet$ \underline{Conclusion}\\
Par injectivité de la transformée de Fourier, on déduit que $f=0$ puis
$g=0$,\\
il s'ensuit :
$$ \overline{{\cal X} } =L^2\Bigl(E,{\cal E}_X,\nu\Bigl).$$
\end{enumerate} \normalsize$\triangle$
\subsection*{Remarque}
Dans {\bf [9]} pX on pourra trouver une autre démonstration de ce
résultat. Dans le cas de variables aléatoires ayant une densité par
rapport à la mesure de Lebesgue, cette proposition relève de la théorie
des polynômes. Dans {\bf [4]} p133, on trouvera une démonstration basée
sur le principe du prolongement analytique; une généralisation à
$\R^n$ de cette méthode impose l'utilisation de la théorie des
fonctions holomorphes de plusieurs variables. Nous avons choisi
une démonstration plus constructive, comme nous le verrons plus
loin, qui se généralise assez facilement.
\begin{de}
On dira qu'une suite de \vas $\Bigl( U_k\Bigl)_{k\in \N}$ admet des
moments exponentiels si :
$$\Bigl(\forall k\in\N\Bigl)\Bigl(\exists \tau_k>0\Bigl)
\Bigl(\forall |t|<\tau_k\Bigl)\Bigl(E( e^{tU_k})<\infty\Bigl).$$
Nous noterons ${\cal E}xp\Bigl( (U_k)_{k\in \N}\Bigl)$ cette
hypothèse.
\end{de}
\subsection*{ii) Deux corollaires au théorème 5.1.1.}
\begin{co}
On suppose de plus que les \vas $(X_n)_{n\in\N^*}$ sont
indépendantes.\\
Si
$$\Biggl(
\forall k\in\N^*
\Biggl)
\Biggl(
Vect_{\R}\Biggl\{X_k^{\alpha}|\alpha\in\N\Biggl\}
\mathop\subset_{dense}\limits
L^2\Bigl(E,{\cal E}_{X_k},\nu\Bigl)
\Biggl)$$
alors,
$$
\bigcup_{n\in\N^*}
Vect_{\R}
\Biggl\{
X_1^{\alpha_1}\cdot\ldots\cdot X_n^{\alpha_n}  |
(\alpha_1,\ldots,\alpha_n)\in\N^n
\Biggl\}
\mathop\subset_{dense}\limits
L^2\Bigl(E,{\cal E},\nu\Bigl).
$$
\end{co}
\subsection*{\underline{\textbf{Démonstration}}}
\begin{enumerate}
\item[]
\small
soit $\nu_k$ (resp. $\nu_X$)  la loi de $X_k$
$\Bigl($resp. $X:=(X_1,\ldots,X_n)\Bigl)$ alors $\nu_X
=\nu_1\otimes\cdots\otimes\nu_n$.\\
Dans les espaces images
$$\begin{array}{ccl}
\Ee     &  \longrightarrow    &  (\R^n,B(\R^n),\nu_X)\\
\omega  &   \longmapsto       &
\Bigl(X_1(\omega),\ldots,X_n(\omega)\Bigl):=(x_1,\ldots,x_n).
\end{array}$$
L'indépendance nous assure l'isomorphisme :
$$L^2\Bigl(\R^n,B(\R^n),\nu_X\Bigl)\simeq \bigotimes_{k=1}^n
L^2\Bigl(\R,B(\R),\nu_k\Bigl)$$
ce qui fournit la solution :\\
$\Bigl\{x_k^{\alpha} | \alpha\in\N\Bigl\}$
est un système générateur de $L^2\Bigl(\R,B(\R),\nu_k\Bigl)$,
$\Bigl\{ \dps\bigotimes_{k=1}^nx_k^{\alpha_k}|
(\alpha_1,\ldots,\alpha_n)\in\N^n\Bigl\}$
est donc générateur de
$\displaystyle\bigotimes_{k=1}^nL^2\Bigl(\R,B(\R),\nu_k\Bigl)$.
Par l'isomorphisme, on remonte à l'espace source et  utilise
le théorème 5.1.1. pour conclure.\\
\end{enumerate} \normalsize$\triangle$
\subsection*{Remarque}
Dans le corollaire suivant, on ne suppose plus les  \vas
indépendantes, ce qui impose une hypothèse plus forte à chacune
d'entre elles.
\begin{co}
Avec les notations du théorème 5.1.1., on suppose que la suite de
\vas
$\Bigl(X_k\Bigl)_{k\in\N^*}$ vérifie l'hypothèse $\Exk$;\\
alors :
$$
\bigcup_{n\in\N^*}
Vect_{\R}\Biggl\{
X_1^{\alpha_1}\cdot\ldots\cdot X_n^{\alpha_n}|
\Bigl((\alpha_1,\ldots,\alpha_n)\in\N^n\Bigl)
\Biggl\}
\mathop\subset_{dense}\limits
L^2\Ee.$$
\end{co}
\subsection*{\underline{\textbf{Démonstration}}}
\begin{enumerate}
\item[]
\small
Elle va consister en une généralisation de la démonstration de la
proposition 5.1.1.\\
Fixons $n\in \N^*$ :
$$\Bigl( X_1,\ldots,X_n\Bigl) \ : \Een\longrightarrow
\Bigl( \R^n,B(\R^n)
\Bigl)
$$
Soit $g\in L^2\Een\cap{\cal X}_n^{\bot}$\\
$\exists f\in L^2\Bigl( \R^n,B(\R^n),\nu_{ X}\Bigl)$
telle que $g=f\circ ( X_1,\ldots,X_n\Bigl)$.\\
Notons :
$$\varphi(t):= E^{\nu}\Bigl(e^{i<t,X>}g\Bigl).$$
Comme $g\in L^2\Een$, $\varphi(t)$ est bien définie pour tout
$t\in\R^n$. En passant dans l'espace image on a :
$$ \varphi(t) = \int_{\R^n} e^{i<t,x>}f(x)d\nu_X(x).$$
$\bullet$ \underline{Régularité de $\varphi$}\\
Soit $(\alpha_1,\ldots,\alpha_n)\in\N^n$
$$\begin{array}{rcl}
%
\varphi_{(\alpha_1,\ldots,\alpha_n)}(t)    &  :=      &
\dps\int_{\R^n}
|\prod_{k=1}^{n}x_k^{\alpha_k}| \cdot |e^{i<t,x>}f(x)| d\nu_X(x)\\
%
                                 &
\mathop{\leq}_{ \mbox{CS} }\limits       & ||f||_2
\Biggl(
 \dps\int_{\R^n}  |\prod_{k=1}^{n}x_k^{2\alpha_k}| d\nu_X(x)
\Biggl)^{\frac{1}{2}}\\
%
                                 &
 \mathop{\leq}_{\stackrel{\mbox{lois marginales}}{\mbox{et  CSG}}}
\limits     &
||f||_2
\Biggl(
 \dps\prod_{k=1}^{n}  \int_{\R}  |x_k^{2n\alpha_k}| d\nu_k(x)
\Biggl)^{\frac{1}{2n}}\\
%
                                 &\leq   &+\infty.
\end{array}$$
La dernière inégalité est vraie par l'existence de moments
exponentiels qui assure l'intégrabilité de tout polynôme. On a noté
CSG pour exprimer qu'on utilise l'inégalité de Cauchy-Schwartz
Généralisée. Le théorème de dérivation de Lebesgue assure l'existence
de toutes les dérivées de $\varphi$.\\
$\bullet$ \underline{Existence d'un ouvert où $\varphi$ et toutes ses
dérivées sont  nulles}\\
f étant orthogonale aux polynômes, on va développer $e^{i<t,x>}$
en série entière de façon à faire apparaître des produits
$\dps\prod_{k=1}^{n}x_k^{\alpha_k}$.\\
Comme
$$\begin{array}{rcl}
\psi (t)    &:= &
\dps\int_{\R^n}
\sum_{k\in\N} |\frac{(i<t,x>)^k}{k!}|f(x)|  d\nu_X(x)\\
            & = &
\dps\int_{\R^n}  e^{|<t,x>|}|f(x)|  d\nu_X(x)\\
             &  \mathop{\leq}_{CS}\limits      &
|| f ||_2 \Biggl(
\dps\int_{\R^n}
               \mbox{\LARGE{e}}^{2\dps\sum_{k=0}^n|t_kx_k|} d\nu_X(x)
\Biggl)^{\frac{1}{2}}\\
             &
\mathop{\leq}_{\stackrel{\mbox{lois marginales}}{\mbox{et  CSG}}}
\limits     &  || f||_2
\Biggl(
\dps\prod_{k=1}^{n}
\int_{\R}
        \mbox{\Large{e}}^{2n|t_kx_k|} d\nu_k(x)
\Biggl)^{ \frac{1}{2n} }\\
%
            & < & +\infty\  \mbox{ pour  }
\dps|t_k|<\frac{\tau_k}{2n}.
\end{array}$$
Sur $D:= \dps\prod_{k=1}^{n}
]-\frac{\tau_k}{2n},+\frac{\tau_k}{2n}[$ on peut intégrer terme à
terme la série $\dps\sum_{k\in\N}\frac{(i<t,x>)^k}{k!}$.
$$\begin{array}{lcl}
%
\varphi(t)   &  =   &  \dps\sum_{k\in\N}\frac{i^k}{k!}
\int_{\R^n}  <t,x>^kf(x) d\nu_X(x)\\
%
<t,x>^k    &  =  &
\Biggl(\dps\sum_{j=1}^k t_jx_j\Biggl)^k =\sum_{k_1+\cdots+k_n=k}
\frac{k!}{k_1!\cdots k_n!}
\prod_{j=1}^{n}(t_jx_j)^{k_j}\\
%
\varphi(t)     &  =   &
\dps\sum_{k\in\N}\frac{i^k}{k!}
\sum_{k_1+\cdots+k_n=k}
\frac{k!}{k_1!\cdots k_n!}
\prod_{j=1}^{n}(t_j)^{k_j}
\int_{\R^n} \prod_{j=1}^{n}(x_j)^{k_j}f(x) d\nu_X(x)\\
%
\varphi(t)     &  =   &
\dps\sum_{k\in\N}\frac{i^k}{k!}
\sum_{k_1+\cdots+k_n=k}
\frac{k!}{k_1!\cdots k_n!}
\prod_{j=1}^{n}(t_j)^{k_j}
E^{\nu}\Biggl(\prod_{j=1}^{n}(X_j)^{k_j}g \Biggl).
\end{array}$$
De la même façon on montre que, pour $p=p_1+\cdots+p_n$
$\Bigl((p_1,\ldots,p_n)\in\N^n\Bigl)$ :
$$\Biggl(
\frac{\partial^p\varphi}
{\partial t_1^{p_1}\cdot\ldots\cdot\partial t_n^{p_n}}
\Biggl)(t) =
\sum_{k\in\N}\frac{i^{k+p}}{k!}
\sum_{k_1+\cdots+k_n=k}
\frac{k!}{k_1!\cdots k_n!}
\prod_{j=1}^{n}(t_j)^{k_j}
E^{\nu}\Biggl(\prod_{j=1}^{n}(X_j)^{k_j+p_j}g \Biggl).$$
Comme $g\in\Bigl({\cal X}_n\Bigl)^{\bot}$, $\varphi_{|_D}\equiv0$ et
pour tout $p=p_1+\cdots+p_n$ $\Bigl((p_1,\ldots,p_n)\in\N^n\Bigl)$,
$\Biggl(\dps\frac{\partial^p\varphi}
{\partial t_1^{p_1}\cdot\ldots\cdot\partial t_n^{p_n}}
\Biggl)_{|_D}\equiv0.$\\
$\bullet$ \underline{ Prolongement }\\
Commençons par réduire l'ensemble D pour le rendre symétrique et plus
facile à manipuler.\\
Posons :
\begin{enumerate}
\item[-]
$\tau_0:=
\inf\Bigl\{
\frac{\tau_1}{2n},\ldots,\frac{\tau_n}{2n}
\Bigl\}$
\item[-]
$|t|_{\infty} := \dps\sup_{k=1,\ldots,n}|t_k|$
\item[-]
$D_p:=\Bigl\{
t\in\R^n  \ | \  |t|_{\infty} <  \frac{\tau_0}{2n} + p\frac{\tau_0}{2n}$\
$p\in\N\Bigl\}$ pour $p=0$ on a $D_0=D$.
\end{enumerate}
On a $\frac{\partial^k\varphi}
{\partial x_1^{k_1}\cdot\ldots\cdot\partial x_n^{k_n}}_{|_{D_0}}=0$, c'est le point précédent.\\
Un raisonnement par récurrence sur p permet de prolonger ce résultat à
$\;\R^n$ .\\
Soient des entiers naturels $j,\;j_1,\ldots,j_n$ tels que $j_1\cdots j_n=j$ et
$t\in D_{p+1}$, $\exists t_p\in D_p$ tel que $|t-t_p|<\frac{\tau_0}{2n}$, on peut écrire :
$$\Biggl(\frac{\partial^j\varphi}
{\partial x_1^{j_1}\cdot\ldots\cdot\partial x_n^{j_n}}\Biggl)(t) =
\dps\int_{\R^n}\prod_{r=1}^{n}(ix_r)^{j_r}  e^{i<t-t_p,x>}e^{i<t_p,x>}f(x)
d\nu_X(x).$$
$\left.
\begin{array}{cc}
\bullet  & |\dps\prod_{r=1}^{n}(ix_r)^{j_r}e^{i<t_p,x>}f| =
|\dps\prod_{r=1}^{n}(x_r)^{j_r}f|\\
\bullet  & |t-t_p|<\frac{\tau_0}{2n}
\end{array}\right\}$
$\Longrightarrow
\left\{\begin{array}{ll}
   & \mbox{on peut intégrer terme à terme}\\
                  & \mbox{la série }\\& \dps\sum_{k\in\N}
 \dps\frac{\dps(i<t-t_p,x>)^k\prod_{r=1}^{n}(ix_r)^{j_r}}{k!}.
\end{array}\right.$\\
Ce qui nous donne
$$\Biggl(\frac{\partial^j\varphi}
{\partial x_1^{j_1}\cdot\ldots\cdot\partial x_n^{j_n}}\Biggl)(t)=\sum_{k\in\N}
\sum_{k_1+\cdots +k_n = k} \prod_{j=1}^n \frac{(t_j-(t_p)_j)^k}{k_j!}
\Biggl(
\frac{\partial^{k+j}\varphi}
{\partial x_1^{k_1+j_1}\cdot\ldots\cdot\partial x_n^{k_n+j_n}}
\Biggl)
(t_p)=0$$
et par suite $\dps\Biggl(\frac{\partial^j\varphi}
{\partial x_1^{j_1}\cdot\ldots\cdot\partial x_n^{j_n}}
\Biggl)_{\left\vert\dps\bigcup_{p\in\N}D_p\right.}=0$.\\
D'autre part
$$\Bigl(\forall t\in\R^n\Bigl)
\Bigl(
(\exists p\in\N)(|t|_{\infty}<\frac{\tau_0}{2n}+p\frac{\tau_0}{2n}
=\frac{1+p}{2n}\tau_0)\Bigl)$$
puisque $\dps\lim_{p\rightarrow+\infty}
\bigl(\frac{1+p}{2n}\tau_0\bigl) =+\infty$\\
donc $\dps\bigcup_{p\in\N}D_p= \R^n$ d'où finalement
$\dps\frac{\partial^j\varphi}
{\partial x_1^{j_1}\cdot\ldots\cdot\partial x_n^{j_n}}\equiv 0$.\\
$\bullet$ \underline{Conclusion}\\
Par injectivité de la transformée de Fourier, on déduit que $f=0$ puis
$g=0$,\\
il s'ensuit :
$$ \overline{{\cal X}_n } =L^2\Een.$$
Le théorème 5.1.1. nous permet alors de conclure.\\
\end{enumerate} \normalsize$\triangle$
\subsection*{Remarque}
Sous les hypothèses du corollaire 5.1.2., les \vas étant indépendantes,
le résultat découle immédiatement du corollaire 5.1.1. et de la proposition
5.1.1.
\section{Application}
\subsection{Explicitation du modèle}
Considérons :
\begin{enumerate}
\item[-]
$L^2\Ee$ donné en 1
\item[-]
un espace de Hilbert H, séparable, muni d'un produit scalaire $< , >$
\item[-]
$\Bigl( e_ k\Bigl)_{k\in \N^*}$ une base hilbertienne de H.
\end{enumerate}
Notre modèle va être la donnée de ces deux espaces de Hilbert,
le premier muni de sa base hilbertienne, le second de la suite de \vas
$\Bigl( X_ k\Bigl)_{k\in \N^*}$ et d'une application $\phi$ les liant
l'un à l'autre de la façon suivante :
$$\begin{array}{cccl}
\phi :  &  H    &  \longrightarrow    &  L^2\Ee\\
      & \dps h=\sum_{k\in \N^*}<h,e_k>e_k &   \longmapsto       &
\dps\sum_{k\in \N^*}<h,e_k>X_k.
\end{array}$$
Je renvois à {\bf [7]} chap. 1 pour avoir plus  détails.
\subsection{Etude de la convergence dans tous les $L^p$}
\subsection*{i)  Convergence exponentielle}
On se place dans les conditions d'application du corollaire 5.1.2. On
va avoir besoin d'une hypothèse plus forte, il nous faudra un
intervalle non trivial sur lequel toutes les \vas admettent des
moments exponentiels.
\begin{pp} (Convergence exponentielle)\\
Supposons que $\dps\delta:=\inf_{k\in\N} \tau_k>0$ et soit
$\Bigl( h_ k\Bigl)_{k\in \N^*}$ une suite d'éléments de H,\\
alors :
$$\Biggl(  h_n\mathop{\longrightarrow}^{H}\limits 0\Biggl)
\Longrightarrow
\Biggl( \Bigl(
\forall |t|<\frac {\delta}{4}
\Bigl) \
\Bigl(
e^{|t\phi(h_n)|}
\mathop{\longrightarrow}^{L^1(\nu)}\limits 1
\Bigl)\Biggl).$$
\end{pp}
\hspace*{-0,6cm}Pour démontrer cette proposition on va borner
uniformément (en n),$\Bigl(e^{|t\phi(h_n)|}\Bigl)^2$, ceci  découlera
du résultat suivant :
\begin{lm}(Proposition 2.2 {\bf [7]})\\
Soient $\Bigl( a_ k\Bigl)_{k\in \N}\in l^2(\N)$ et $\Bigl( V_
k\Bigl)_{k\in\N}$
une suite de \vas indépendantes centrées, réduites (vaicr).  On
suppose
${\cal E}xp\Bigl( (V_k)_{k\in \N}\Bigl)$.\\
Si $\dps\tau:=\inf\frac{\tau_k}{|a_k|}>0$\\
alors $\dps U:=\sum_na_nV_n$
vérifie ${\cal E}(U)$ sur l'intervalle
$]-\frac{\tau}{2},+\frac{\tau}{2}[$ et
$E e^{tU}\leq e^{\frac{9}{8} t^2\sum a_k^2}$.
\end{lm}
\subsection*{\underline{\textbf{Démonstration de la proposition 5.2.1.}}}
\begin{enumerate}
\item[]
\small
\underline{Idée} : utiliser le théorème de convergence de Vitali\\
$\bullet$ \underline{Convergence en probabilité}\\
Soit $\eta>0$
$$\begin{array}{cclc}
\nu\biggl( [e^{|t\phi(h_n)|}>1+\eta] \biggl)   &   =   &
\nu\biggl(\Bigl[|t\phi(h_n)|>\ln(1+\eta) \Bigl]\biggl)     &
\\
                                              & \leq  &
\dps\frac {t^2  ||h_n||^2}{ \ln(1+\eta) }     &  B. Tchebychef
\end{array}$$
La convergence de $\Bigl( h_k\Bigl)_{k\in \N}$ donne alors le résultat.
\\
$\bullet$ \underline{Equi-intégrabilité} : \\
$$0\leq  \dps e^{2|t\phi(h_n)|} \leq    \dps e^{2t\phi(h_n)}
+e^{-2t\phi(h_n)}.$$
Posons $a^n_k:=<h_n,e_k>$, alors $|a_n^k|<||h_n||$.\\
Comme $h_n\mathop{\longrightarrow}^{H}\limits0$ :
$$(\exists N\in\N)(\forall n>N)(\forall k\in \N, |a_n^k|<1)$$
pour $n>N$, $\dps\inf\frac{\tau_k}{|a_k^n|}>\tau$.\\
Le lemme 5.2.1. nous montre alors que
$E e^{2t\phi(h_n)}\leq e^{\frac{9}{8} (2t)^2|| h ||^2}$
$(t\in]-\frac{\tau}{4},+\frac{\tau}{4}[$), d'où l'encadrement\\
suivant :
$$0\leq E  e^{2|t\phi(h_n)|}\leq2 e^{\frac{9}{2} t^2|| h_n||^2}
\leq
2e^{\frac{9}{2} t^2}\  \
(t\in]-\frac{\delta}{4},+\frac{\delta}{4}[).$$
On a ainsi majoré uniformément
 $\sup_n E\Bigl(e^{|t\phi(h_n)|}\Bigl)^2$
ce qui implique l'équi-intégrabilité. Le théorème
de convergence de Vitali  donne le résultat, à savoir, la
convergence dans $L^1(\nu)$ :
$$ E e^{|t\phi(h_n)|} \longrightarrow 1.$$
\end{enumerate} \normalsize$\triangle$
\subsection*{ii) Convergence dans tous les $L^p$}
\begin{pp} (Convergence dans tous les $L^p$)\\
Sous les hypothèses d'application de la proposition 5.2.1. on a :
$$\Bigl( \forall p\geq 1\Bigl)\Bigl( \phi(h_n)
\mathop{\longrightarrow}^{L^p(\nu)}\limits 0 \Bigl).$$
\end{pp}
\subsection*{Remarque}
Dans l'ensemble des suites ayant des moments exponentiels cette
proposition nous offre une réciproque  à la précédente. En
effet, on peut affirmer en particulier que
$\phi(h_n)\mathop{
\longrightarrow}^{L^2(\nu)}\limits 0$. La
réciproque découle alors de l'isométrie :
$$||h_n||_{H} = ||\phi(h_n)||_{L^2(\nu)}.$$
\subsection*{\underline{\textbf{Démonstration}}}
\begin{enumerate}
\item[]
\small
Fixons $n\in\N$ et $t\in]-\frac{\tau}{4},+\frac{\tau}{4}[$. On a le
développement en série entière suivant :
$$ E e^{|t\phi(h_n)|} = \sum_{k\in\N} \frac{ |t|^k}{
k!}E|\phi(h_n)|^k$$
(série à termes positifs, à valeurs dans $\;\overline{\R^+}$).
On veut prendre la limite des deux membres quand n tend vers $+\infty$.
Vérifions qu'on peut inverser
$\mathop{\lim}_{n\rightarrow\infty}\limits$ et $\dps \sum_{k\in\N}$.\\
Pour cela fixons $t_0\in]0,+\frac{\tau}{4}[$,
$$\Bigl( \exists N\in\N\Bigl)
\Bigl( \forall n\in\N, n>N\Bigl)\Longrightarrow
\Bigl( ||h_n||\leq1\Big).
$$
Prenons $n>N$, $|t|<t_0$ et $k\geq1$ un entier naturel. Nous devons
montrer que la série converge uniformément sur $\N$ :
$$\begin{array}{ccccc}
|\phi(h_n)|^k  & \leq  &
\dps\frac{k!}{t_0^k}  e^{t_0|\phi(h_n)|}  &  \leq  &
\dps\frac{k!}{t_0^k} \Bigl(e^{t_0\phi(h_n)} +  e^{-t_0\phi(h_n)}\Bigl)\\
E|\phi(h_n)|^k &  \leq  &
\dps\frac{k!}{t_0^k} 2 e^{\frac{9}{8}t_0^2||h_n||^2}  &
\mbox{par le lemme 2.1}\\
              &  \leq    &
\dps\frac{k!}{t_0^k} 2 e^{\frac{9}{8}t_0^2}  &
\mbox{car } n\geq N.
\end{array}$$
Mais
$$\Biggl(
\frac{k!}{t_0^k}2 e^{\frac{9}{8}t_0^2}
\Biggl)\Biggl(
\frac{|t|^k}{k!}
\Biggl) =
2 e^{\frac{9}{8}t_0^2}\Biggl(\frac{|t|}{t_0}\Biggl)^k$$
est le terme général d'une série géométrique convergente, ainsi
$$\sum_{k\in\N} E|\phi(h_n)|^k\frac{|t|^k}{k!}$$
converge normalement donc uniformément sur $\N$; par suite
$$1 =\mathop{\lim}_{n\rightarrow+\infty}\limits
Ee^{|t\phi(h_n)|} =
1 +\sum_{k\in\N^*}
\mathop{\lim}_{n\rightarrow+\infty}\limits
E|\phi(h_n)|^k\frac{|t|^k}{k!}$$
et
$$ 0 = \sum_{k\in\N^*} \mathop{\lim}_{n\rightarrow+\infty}\limits
E|\phi(h_n)|^k\frac{|t|^k}{k!}$$
d'où pour tout $k\in\N^*$ :
$$
\mathop{\lim}_{n\rightarrow+\infty}\limits
E|\phi(h_n)|^k = 0
$$
$\bullet$ \underline{Généralisation} \\
Comme nous travaillons avec des probabilités
$\left\{
\begin{array}{cccc}
L^p       & \subset   &L^q        &  p>q\\
||\ ||_q  & \leq      &|| \ ||_p.  &
\end{array}
\right.$\\
D'autre part
$$\Bigl(\forall p>1\Bigl)
\Bigl( [p]\leq p<[p]+ 1\Bigl)\mbox{  (partie entière)}.$$
Pour toute suite $\Bigl( U_ k\Bigl)_{k\in \N}$ de $L^{[p]+1}$ on a :
$$||U_n||_p\leq||U_n||_{[p]+1}$$
Ce qui nous démontre la proposition 5.2.1 pour tout $p\geq1$.\\
\end{enumerate} \normalsize$\triangle$
\subsection{Une application}
Tous ces préliminaires nous permettent d'énoncer le résultat qui nous
intéresse dans cette partie. Nous conservons les hypothèses précédentes
et utiliserons un espace de Hilbert particulier. On va prendre
$H:=L^2\Bigl(\R,B(\R),\mu\Bigl)$ où $\mu$ est une mesure bornée.\\
Introduisons quelques notations :
\begin{itemize}
\item
$W_0 := \R$
\item
$W_n := Vect_{\R}\Bigl\{W_{t_1}\cdot\ldots\cdot W_{t_n}  |
(t_1,\ldots,t_n)\in\R^n\Bigl\}$ $(n\geq1)$
\item
$W := \dps\bigcup_{n\in\N} W_n$
\end{itemize}
\begin{thm}
Dans les conditions d'application des propositions 5.2.1. et 5.2.2.
$$
\bigcup_{n\in\N} W_n\;
\mathop{\subset}_{dense}\limits
L^2\Bigl(\R,B(\R),\mu\Bigl)
$$
\end{thm}
\subsection*{\underline{\textbf{Démonstration}}}
\begin{enumerate}
\item[]
\small
Pour tout $k\in\N^*$ $e_k\in L^2\Bigl(\R,B(\R),\mu)$, il existe une
suite de fonctions simples $\Bigl( e_k^j\Bigl)_{j\in\N}$ convergeant
dans $L^2\Bigl(\mu\Bigl)$ (par densité) vers $e_k$.\\
C'est à dire, pour tout k fixé dans $\N^*$, il existe :
\begin{itemize}
\item
$\Bigl(N(j,k)\Bigl)_{j\in\N}$ suite  strictement croissante d'éléments
de $\N$,
\item
$\biggl( j\in\N\biggl)
\biggl(
\alpha_i(j,k)\in\R\;\Bigl(i=1,\ldots,N(j,k)\Bigl)
\biggl)$
\item
$t_0(j,k)<\cdots<t_N(j,k)$ des réels
\end{itemize}
tels que :
$$||
\sum_{i=0}^{N(j,k)}
\alpha_i (j,k)
1_{]t_i(j,k),t_{i+1}(j,k)[}
-
e_k\;
||_2
\mathop{\longrightarrow}_{j\rightarrow+\infty}\limits0.
$$
Pour simplifier on omettra de préciser le couple (j,k), on notera
$\Bigl( e_k^j\Bigl)_{j\in\N}$ cette suite et $\Bigl(
X_k^j\Bigl)_{j\in\N}$ son image par $\phi$.
\begin{enumerate}
%
\item
\underline{$\Bigl( j\in\N         \Bigl)
          \Bigl( \forall k\in\N \Bigl)
          \Bigl( X_k^j\in W_1   \Bigl)$}\\
Il suffit de voir que :
$$\phi (e_k^j) = \dps\sum_{i=0}^N
                 \alpha_i
                 \Bigl(  W_{t_{i+1}}-W_{t_{i}}  \Bigl).$$
%
\item
\underline{$\Bigl( \forall p\geq 1\Bigl)
           \Bigl(
X_k^j\mathop{\rightarrow}^{L^p}\limits X_k
           \Bigl)$}\\
Cela résulte de la proposition 2.1.
%
\item
\underline{$
\Bigl(  \forall n\in\N^* \Bigl)
\Bigl(  \forall(k_1,\ldots,k_n)\in\N^n\Bigl)
\Bigl(  X_{k_1}\cdot\ldots\cdot X_{k_n}\in\overline{W}\Bigl)$}\\
Par récurrence sur n.\\
n=1 :
$$\Bigl(  \forall p\geq1\Bigl)
  \Bigl(  X_k^j\mathop{\rightarrow}^{L^p}\limits X_k\Bigl)\;
\mbox{(point 1)}.$$
HR : fixons $n\in\N^*$ et supposons que  :
$$\Bigl(  \forall s\in\{1,\ldots,n\}\Bigl)
\Bigl(  \forall p\geq1\Bigl)
\Bigl(  \forall(k_1,\ldots,k_s)\in\N^s\Bigl)
\Bigl(  X_{k_1}^j\cdot\ldots\cdot X_{k_s}^j
\mathop{\longrightarrow}^{L^p}\limits
X_{k_1}\cdot\ldots\cdot X_{k_s}
\Bigl)$$
Prenons maintenant $(k_1,\ldots,k_{n+1})\in(\N^*)^{n+1}$. Alors :\\
$||
X_{k_1}^j\cdot\ldots\cdot X_{k_{n+1}}^j -
                        X_{k_1}\cdot\ldots\cdot X_{k_{n+1}}||_2$\\
$\leq
||
\Bigl(X_{k_1}^j\cdot\ldots\cdot X_{k_{n}}^j -
     X_{k_1}\cdot\ldots\cdot X_{k_{n}}\Bigl)X_{k_{n+1}}
||_2 +
||X_{k_1}\cdot\ldots\cdot X_{k_{n}}
\Bigl(X_{k_{n+1}}^j-X_{k_{n+1}}\Bigl)
||_2$\\
$\leq
||
X_{k_1}^j\cdot\ldots\cdot X_{k_{n}}^j -
     X_{k_1}\cdot\ldots\cdot X_{k_{n}}||_4
||  X_{k_{n+1}}  ||_4 +
||X_{k_1}\cdot\ldots\cdot X_{k_{n}}
||_4
||X_{k_{n+1}}^j-X_{k_{n+1}}
||_4$.\\
Or $\Bigl(
X_{k_{n+1}}^j
\mathop{\rightarrow}^{L^4}\limits
X_{k_{n+1}}
\Bigl)$ donc
$$\Bigl( \exists j_0\in\N \Bigl)
\Bigl(  \forall j\geq j_0\Bigl)
\Bigl(
||X_{k_{n+1}}^j||_4
\leq
||X_{k_{n+1}}||_4  +1
\Bigl)$$
en utilisant l'hypothèse de récurrence (HR) et le point 1), on a :
$$
X_{k_1}\cdot\ldots\cdot X_{k_{n+1}} \in
\overline{W}_{n+1}
\subset
\overline{W}
$$
et par suite
$$
\Bigl( \forall n\in\N^* \Bigl)
\Bigl(  \forall (k_1,\ldots,k_n)\in\N^n \Bigl)
\Bigl(  X_{k_1}\cdot\ldots\cdot X_{k_n}\in\overline{W}\Bigl)$$
finalement ${\cal X}\subset \overline{W}\subset L^2$ donc
$${\cal X} = \overline{W}= L^2\Ee.$$
\end{enumerate}
\end{enumerate} \normalsize$\triangle$
%
%
%
%
%
%
%
%
%
%
%
%
%
%
\part{CONSTRUCTION ET ÉTUDE D'UNE INTÉGRALE}
\chapter{Construction d'une intégrale}
\section{Résultats techniques}
\subsection{Modèle}
Nous utiliserons les données suivantes :
\begin{itemize}
\item[i)]
\begin{itemize}
\item E un espace vectoriel normé
\item ${\cal E}$ une tribu sur E
\item $({\mathcal E}_n)_{n\in\N* }$ une filtration de ${\cal E}$
\item $\nu$ une probabilité sur (E,${\cal E}$)
\item $(X_k)_{k\in\N^*}$ une suite de vaicr $\R$ sur $(E,{\cal E}, \nu)$ pour laquelle il existe, pour tout $p>1$, une contante $K_p$ telle que :
$$\sup_{k\in\N^*} E\vert X_k\vert^p<K_p$$
\end{itemize}
\item[ii)]
\begin{itemize}
\item $\mu$ une probabilité diffuse sur $({\R} , B(\R))$
\item $(e_j)_{j\in\N^*}$ une base hilbertienne de H:=$L^2(E,{\cal E}, \nu)$
\end{itemize}
\item[iii)]
\begin{itemize}
\item $\Phi$ :
$H\to L^2 \Bigl(   \R, B(\R), \mu   \Bigl)$
$h \mapsto\dps \sum_{k\in\N*} \langle  h , e_k  \rangle X_k$.
Nous noterons $\Phi(h)_n$ la somme partielle $\dps\sum_{k=1}^{n}\langle h , e_k\rangle X_k$
\end{itemize}
\end{itemize}
\subsection*{{Remarques} }
\begin{itemize}
\item
Pour $n\geq1$ la famille
$\Bigl(e_{j_1}\otimes\cdots\otimes e_{j_n}\Bigl)_{(j_1,\ldots,j_n)
\in\bigl(\N^*\bigl)^n}$
est naturellement une base hilbertienne de $H^{\otimes ^n}$
\item
Pour toute fonction $f\in H^{\otimes ^n}$ on a
$f=\displaystyle  \sum_{(j_1,\ldots,j_n)\in(N^*)^n}
\langle  f , e_{j_1}\otimes\cdots\otimes e_{j_n} \rangle
e_{j_1} \otimes \cdots \otimes e_{j_n}$
\item
Dans le cas où les \vas $\bigl(X_k\bigl)_{k\in\N^*}$ sont
gaussiennes on sait, confère par exemple {\bf{[7]}}, qu'il existe un
mouvement brownien $\beta$ pour lequel :
$$\Phi(h)=\int_0^th(s)d\beta_s$$
Nous généraliserons cette idée dans la suite en construisant une
intégrale dans le cas où les \vas $\bigl(X_k\bigl)_{k\in\N^*}$ ne
sont pas gaussiennes.
\end{itemize}
\begin{pp}
$\Phi$ est une isométrie pour les normes usuelles
\end{pp}
\subsection{Lemmes}
\begin{lm}
Considérons les données suivantes :
\begin{enumerate}
\item 
$\Bigl(M_n\Bigl)_{n\in\N^*}$
une martingale sur
$\Bigl( E, {\cal E}, ( {\cal E}_n )_{n\in\N^*}, \nu \Bigl)$
\item 
$\Bigl( Y_n \Bigl)_{n\in\N^*}$
un processus prévisible sur
$\Bigl( E, {\cal E}, ( {\cal E}_n )_{n\in\N^*}, \nu \Bigl)$
\end{enumerate}
Si l'on suppose les deux hypothèses suivantes satisfaites :
\begin{enumerate}
\item
$\displaystyle \sup_{k\geq2} E \Bigl( M_k-M_{k-1} \Bigl)^2<\infty$
\item
$\displaystyle\sum_{k\geq1} E Y_k^2<\infty$
\end{enumerate}
alors la martingale transformée :
$$\Bigl( Y\ast M\Bigl)_n := \sum_{k=2}^{n} Y_k \Bigl( M_k-M_{k-1} \Bigl)$$
est de carré intégrable.
\end{lm}
\subsection*{\underline{\textbf{Démonstration}}}
\begin{enumerate}
\item[]
\small
Choisissons un entier naturel n supérieur ou égal à deux, on a :
$$E \mid \Bigl( Y\ast M\Bigl)_n   \mid^2 \leq
\sup_{k\geq2} E \biggl(   M_k-M_{k-1}   \biggl)^2\sum_{k=0}^{n}  E Y_k^2$$
\end{enumerate} \normalsize$\triangle$
\begin{lm}
Pour un élément f de $H \otimes  H$ les deux suites
\begin{itemize}
\item[$\bullet$]
$\Biggl(
\dps\sum_{j=1}^{N}
\langle     f  ,    e_j   \otimes  e_j     \rangle   (X_j^2-1)
\Biggl)_{N\geq1}
=:\biggl(
\varphi_N^{(2)} (f)
\biggl)_{N\geq1}$
\item[$\bullet$]
$\Biggl(
\dps\sum_{j=2}^{N}
\Bigl(
\sum_{k=1}^{j-1}
\langle     f  ,    e_j   \otimes  e_k     \rangle  X_k   \Bigl)   X_j
\Biggl)_{N\geq2}
=:\biggl(
\varphi_N^{(1,1)} (f)
\biggl)_{N\geq2}$
\end{itemize}
convergent dans $L^2 \bigl( \R, B(\R),\mu\bigl)$
\end{lm}
\subsection*{Notation}
Nous noterons $\varphi^{(2)}(f)$ et $\varphi^{(1,1)}(f)$ leur limite respective.
\subsection*{\underline{\textbf{Démonstration}}}
\begin{enumerate}
\small
 \item[]
On utilise le lemme 6.1.1.
\newline
Pour la première suite on pose :
\begin{enumerate}
\item[$\bullet$]
$
\left\{
\begin{array}{lcl}
Y_k  & = & \langle f , e_k \otimes e_k \rangle\\
M_k  & = & \dps\sum_{j=1}^{k} (X_j^2-1)
\end{array}\right.
\left\vert
\begin{array}{lcl}
E(M_k-M_{k-1})^2 & \leq & K_4+1\\
&&\\
\dps\sum_{k\geq1}  EY_k^2  & \leq &\parallel   f  \parallel^2_{H\otimes H}
\end{array}\right.
$
\newline
pour la seconde suite on pose :
\item[$\bullet$]
$
\left\{
\begin{array}{lcl}
Y_k  & = &  \dps\sum_{j=1}^{k-1} \langle f , e_k \otimes e_j \rangle X_j\\
M_k  & = & \dps\sum_{j=1}^{k} X_j
\end{array}\right.
\left\vert
\begin{array}{lll}
E(M_k-M_{k-1})^2 &=1&\\
&&\\
\dps\sum_{k\geq2}  EY_k^2& \leq
\dps\sum_{1\leq j< k}\langle f , e_j\otimes e_j \rangle
& \leq \parallel   f  \parallel^2_{H\otimes H}
\end{array}\right.
$
\end{enumerate} \end{enumerate} \normalsize $\triangle$
\begin{lm}
$\varphi^{(2)}$ et $\varphi^{(1,1)}$ sont des applications linéaires continues de
$H\otimes H\to L^2\bigl( E, {\cal E}, \nu \bigl)$
\end {lm}
\subsection*{\underline{\textbf{Démonstration}}}
 \begin{enumerate} \small  \item[]
Soit $f\in H\otimes H$ on a, en utilisant, le lemme 6.1.2. :
\begin{enumerate}
\item
$E^{\nu}\vert \varphi^{(2)}(f)\vert^2 \leq (K_4+2K_2+1)\n f\n^2_{H\otimes H}$
\item
$E^{\nu}\vert \varphi^{(1,1)}(f)\vert^2 \leq\n f\n^2_{H\otimes H}$
\end{enumerate}
\end{enumerate}   \normalsize $\triangle$
\subsection*{Remarques}
\begin{enumerate}
\item[$\bullet$]
On a donc en particulier
$\left\{
\begin{array}{lcl}
\n\varphi^{(2)}\n & \leq & \sqrt{K_4+2K_2+1}\\
\n\varphi^{(1,1)}\n & \leq & 1\\
\end{array}\right.$
\item[$\bullet$]
Dans le cas où les \vas $X_k$ $(k\in\N^*)$ sont équidistribuées on a
$\n\varphi^{(2)}\n =M_4+2M_2+1$. Pour le voir il suffit de considérer
$U_1\in H$ et $U_1\not=0$ avec $U_2=\frac{U_1}{\n U_2 \n^2}$  (on peut également considérer $f=e_1\otimes e_2$).
\item[$\bullet$]
En prenant $f=e_2\otimes e_1$ on a
$E\vert \dps\sum_{n\geq2}\sum_{k=1}^{n-1}
\langle
f , e_n\otimes e_k
\rangle
X_k X_n
\vert^2 =
E\vert X_1 X_2\vert^2 =1$
donc $\n\varphi^{(1,1)}\n=1$
\end{enumerate}
\begin{lm}
Pour $f_1$, $f_2$ dans $H\otimes H$ on a :
\begin{enumerate}
\item[i)]
$E^{\nu}\Bigl( \varphi^{(1,1)}(f_1) \varphi^{(1,1)}(f_2)\Bigl) =
\dps\sum_{j\geq2}\sum_{k=1}^{j-1}
\langle
f_1 , e_j\otimes e_k
\rangle
\langle
f_2 , e_j\otimes e_k
\rangle$
\item[ii)]
$E^{\nu}\Bigl( \varphi^{(1,1)}(f_1) \varphi^{(2)}(f_2)\Bigl) =0$
\item[iii)]
$E^{\nu}\Bigl( \varphi^{(2)}(f_1) \varphi^{(2)}(f_2)\Bigl) =
\dps\sum_{j\geq1}
\langle
f_1 , e_j\otimes e_j\rangle
\langle
f_2 , e_j\otimes e_j\rangle$
\end{enumerate}
\end{lm}
\subsection*{\underline{\textbf{Démonstration}}}
\begin{enumerate}
\item[]
\small
Il suffit de polariser les résultats du lemme précédent
\end{enumerate} \normalsize$\triangle$
\begin{lm}
Si les \vas $X_k$ $(k\in\N^*)$ sont équidistribuées, on a pour
$f_1, f_2\in H\otimes H$ :
\newline
$E^{\nu}\Bigl( \varphi^{(2)}(f_1) \varphi^{(2)}(f_2)\Bigl) =
(M_4-2M_2+1)\dps\sum_{j\geq1}
\langle
f_1\otimes f_2 , e_j^{\otimes^{4}} \rangle$
\newline
$E^{\nu}\Bigl( \varphi^{(1,1)}(f_1) \varphi^{(1,1)}(f_2)\Bigl) =
\langle f_1, f_2 \rangle - \frac{1}{M_4-2M_2+1}E^{\nu}\Bigl( \varphi^{(2)}(f_1)
\varphi^{(2)}(f_2)\Bigl)$
\end{lm}
\subsection*{Remarque}
On vérifie assez facilement que
$$\n f\n^2_{H\otimes H} = \n \varphi^{(1,1)} (f)\n^2 +\frac{1}{M_4-2M_2+1}\n
 \varphi^{(2)}(f)\n^2$$
Ce résultat nous incitera dans la suite à considérer un opérateur $\Phi^{(2)}$ globalisant les deux opérateurs $\varphi^{(1,1)}$ et $\varphi^{(2)}$.\\
Nous garderons, cependant, ces deux derniers dans le but de simplifier les calculs.
\subsection*{\underline{\textbf{Démonstration}}}
\begin{enumerate}
\item[]
\small
On utilise la remarque précédente et la polarisation du produit scalaire.
\end{enumerate} \normalsize$\triangle$
\section{Construction d'une intégrale}
\subsection{Etude de produits}
Nous avons vu au chapitre précédent qu'on obtient un sous-espace dense de $L^2$ en
considérant les produits $X_{k_1}\cdot\ldots\cdot X_{k_n}$
$\bigl(  (k_1,\ldots,k_n)  \in(\N^*)^n, \ n\in\N^*\bigl)$,
sous l'hypothèse d'existence de moments exponentiels pour les variables
$X_k$.\\
Notre but, dans cette thèse, est l'obtention d'une décomposition en chaos.
Nous allons, dans un premier temps, considérer les chaos d'ordre 2 et pour cela
construire une intégrale.
\begin{pp}
Pour des fonctions $f$ et $g$ de H on a :
\begin{enumerate}
\item
$\Phi(h)\Phi(g)\in L^2$
\item
$\Phi(h)\Phi(g) = \varphi^{(2)}(h\otimes g) + \varphi^{(1,1)}(h\otimes g+g\otimes h) +
\langle h , g \rangle$ dans $L^2$
\end{enumerate} \end{pp}
\subsection*{Remarque}
Cette écriture nous poussera dans la suite à considérer des fonctions symétriques de $H\otimes H$ pour faire apparaître des propriétés particulières mais intéressantes.
\subsection*{\underline{\textbf{Démonstration}}}
\begin{enumerate}
\item[]
\small
$\dps\sum_{j=1}^{n}\sum_{k=1}^{j-1}\langle h\otimes g, e_j\otimes e_k\rangle X_j X_k =
\varphi^{(1,1)}_n(h\otimes g+g\otimes h)+\varphi^{(2)}_n(h\otimes g) +
\sum_{k=1}^{n}\langle h\otimes g, e_k\otimes e_k\rangle$
\newline
Il ne reste plus qu'à prendre la limite de chacun des termes dans $L^2$.
\end{enumerate} \normalsize$\triangle$
\begin{pp}
Si nous prenons une suite de partitions de $[0,1]$ : $0=t^n_0<\cdots<t^n_0=1$ dont le pas
 tend vers zéro, alors nous obtenons le résultat suivant :
\newline
$$\dps \sum_{k=0}^{n}\Phi \Bigl( h1_{]0,t^n_k]}\Bigl)\Phi \Bigl( g1_{]t^n_k,t^n_{k+1}]}\Bigl)
\stackrel{L^2}{\to}
\varphi^{(1,1)}(h\otimes g 1_C +g\otimes h 1_{\widetilde{C}}) +\varphi^{(2)}
 (h\otimes g 1_C)$$
où\\
$1_C =\{ (x,y)\in\R^2 \vert 0\leq x<y \leq 0\}$\\
et\\
$1_{\widetilde{C}}=\{ (x,y)\in\R^2 \vert 0\leq y<x \leq 0\}$.
\end{pp}
\subsection*{\underline{\textbf{Démonstration}}}
\begin{enumerate}
\item[]
\small
On utilise la proposition 6.2.1.:
\newline
$\dps\sum_{k=0}^{n}\Phi\Bigl(h1_{  ]t^n_k,t^n_{k+1} ]} \Bigl) =$\\
$\dps\sum_{k=0}^{n}
     \Bigl[
\varphi^{(1,1)}  \Bigl(    h\otimes g  1_{]0,t^n_k]}\otimes1_{  ]t^n_k,t^n_{k+1} ]} +
                                g\otimes h 1_{  ]t^n_k,t^n_{k+1} ]}\otimes1_{]0,t^n_k]}
                         \Bigl)
 +\varphi^{(2)}  \Bigl(    h\otimes g  1_{]0,t^n_k]}\otimes1_{  ]t^n_k,t^n_{k+1} ]}
                        \Bigl)$\\
$\hspace*{2cm}
+\langle h1_{]0,t^n_k]}, g1_{  ]t^n_k,t^n_{k+1} ]}\rangle
\Bigl]$
\newline
On remarque que $\langle h1_{]0,t^n_k]}, g1_{  ]t^n_k,t^n_{k+1} ]}\rangle=0$\newline
On va maintenant étudier les termes : \\
$\left\{
\begin{array}{lcl}
N^n_1   & := &
\left\Vert\dps\sum_{k=0}^{n}\varphi^{(1,1)}\biggl( h\otimes g  1_{]0,t^n_k]}\otimes1_{  ]t^n_k,t^n_{k+1} ]} \biggl)
-\varphi^{(1,1)}\biggl(h\otimes g1_C\biggl)\right\Vert_{L^2(\nu)}\\
N^n_2   & := &
\left\Vert
   \dps\sum_{k=0}^{n}
   \varphi^{(1,1)}\biggl(
                                 g\otimes h  1_{  ]t^n_k,t^n_{k+1}]} \otimes1_{]0,t^n_k]}
                      \biggl)
 -  \varphi^{(1,1)}\biggl(
                                  g\otimes h1_{\widetilde{C}}
                       \biggl)
\right\Vert_{L^2(\nu)}\\
N^n_3   & := &
\left\Vert\dps\sum_{k=0}^{n}
\varphi^{(2)}\biggl(
                              h\otimes g  1_{]0,t^n_k]}\otimes1_{  ]t^n_k,t^n_{k+1} ]}
                    \biggl)
-\varphi^{(2)}\biggl(   h\otimes g1_C      \biggl)       \right\Vert_{L^2(\nu)}
\end{array}\right.$
\newline
$\begin{array}{lcl}
N^n_1   & = &
\left\Vert\dps\sum_{k=0}^{n}
\varphi^{(1,1)}\biggl(
                                 h\otimes g  1_{]0,t^n_k]}\otimes1_{  ]t^n_k,t^n_{k+1} ]} -   h\otimes g1_C
                      \biggl)
\right\Vert_{L^2(\nu)}\\
          & \leq &  \left\Vert h\otimes g
\biggl(
         \dps\sum_{k=0}^{n}
 1_{]0,t^n_k]}\otimes1_{  ]t^n_k,t^n_{k+1} ]} - 1_C
\biggl)
\right\Vert_{L^2(\mu^{\otimes^2})}\\
\end{array}$
\newline
Nous avons :
\begin{enumerate}
\item[$\bullet$]
$\biggl( h\otimes g\biggl)^2
\biggl(
          \dps\sum_{k=0}^{n} 1_{]0,t^n_k]}\otimes1_{  ]t^n_k,t^n_{k+1} ]} - 1_C
\biggl)^2 \leq
4\biggl( h\otimes g\biggl)^2\in L^2(\mu^{\otimes^2})$
\item[$\bullet$]
$\biggl(
          \dps\sum_{k=0}^{n} 1_{]0,t^n_k]}\otimes1_{  ]t^n_k,t^n_{k+1} ]} - 1_C
\biggl)^2\dps\mathop{\stackrel{\mu^{\otimes^2}}{\longrightarrow}}\limits_{\mbox{\ ps}}0$
\end{enumerate}
$\Longrightarrow
\biggl( h\otimes g\biggl)^2
\biggl(
          \dps\sum_{k=0}^{n} 1_{]0,t^n_k]}\otimes1_{  ]t^n_k,t^n_{k+1} ]} - 1_C
\biggl)^2\dps\mathop{\stackrel{\mu^{\otimes^2}}{\longrightarrow}}\limits_{\mbox{\ ps}}0$
\newline
En appliquant le théorème de convergence dominée de Lebesgue on a : \\
$$
\left\Vert h\otimes g
\biggl(
         \dps\sum_{k=0}^{n}
 1_{]0,t^n_k]}\otimes1_{  ]t^n_k,t^n_{k+1} ]} - 1_C
\biggl)
\right\Vert_{L^2(\mu^{\otimes^2})}{\mathop{\longrightarrow}\limits_{n\to\infty}}0$$
d'où
$$\dps\sum_{k=0}^{n}\varphi^{(1,1)}\biggl( h\otimes g  1_{]0,t^n_k]}\otimes
1_{  ]t^n_k,t^n_{k+1} ]} \biggl)
\dps\mathop{  \stackrel{L^2(\nu)}{\longrightarrow} }\limits_{n\to\infty}
\varphi^{(1,1)}\biggl( h\otimes g1_C\biggl)
$$
Les termes $N_2^n$ et $N_3^n$ se traitent de la même façon.
\end{enumerate} \normalsize$\triangle$
\subsection{Définition de l'intégrale}
Nous avons maintenant tous les matériaux pour définir l'intégrale qui nous intéresse,
on la notera, pour $h$ et $g$ dans H :
$$\int \Phi(h)_sd\Phi(g)_s$$
où $\Phi(h)_s:=\Phi(h1_{]0,s]})$.
\begin{de}
$$\int \Phi(h)_sd\Phi(g)_s:=\varphi^{(1,1)}\bigl(h\otimes g1_C+
g\otimes h1_{\widetilde{C}}\bigl)+\varphi^{(2)}\bigl(h\otimes g1_C\bigl)$$
\end{de}
\subsection*{Remarque}
C'est une intégrale dans le sens suivant : \\
$\dps\int \Phi(h)_sd\Phi(g)_s =\dps\lim_{L^2}\sum_{k=0}^{n}\Phi(h)_{t_{k}^n}
\Bigl(\Phi(g)_{t_{k+1}^n}-\Phi(g)_{t_{k}^n}\Bigl)$\\
$=\dps\lim_{L^2}\sum_{k=0}^{n}\Phi\Bigl(h1_{]0,t_{k}^n]}\Bigl)
                                                 \Phi\Bigl(g1_{]t_{k}^n,t_{k+1}^n]}\Bigl)$
\section{Une autre présentation du même objet}
\subsection{Motivation}
La décomposition précédente est faite selon la base des polynômes $xy$ et $x^2-1$.
On pense soit aux polynômes de Wick,  on a recentré, soit aux polynômes
orthogonaux.
On sait associer à la mesure $\mu$ (confère {\bf{[2]}}) une famille de polynômes orthogonaux
$(P_n)_{n\in\N}$ où pour tout entier naturel n le degré du polynôme $P_n$ est exactement n.
\\
En considérant les produits $P_{i_1}(X_{j_1})\cdot\ldots\cdot P_{i_n}(X_{j_n})$, où les n-uples
$(i_1,\ldots,i_n)$ et $(j_1,\ldots,j_n)$ sont formés d'entiers naturels distincts, nous obtenons
une décomposition de $L^2(\nu)$ selon les chaos sur la famille de \vas $(X_k)_{k\in\N^*}$.
L'indépendance des \vas et l'orthogonalité des polynômes nous donne une famille orthogonale,
 mais plus normée; nous aurons l'occasion de résoudre ce problème dans la suite.\\
Nous allons réécrire cette décomposition selon ces polynômes orthogonaux de façon à rendre
les calculs qui suivent plus simples. Cette décomposition nous offre, de plus, une généralisation
plus aisée.
\subsection{Le modèle}
Nous allons utiliser deux opérateurs. Comme nous l'avons déjà remarqué, l'introduction de
fonctions symétriques simplifie l'écriture de notre intégrale. Pour garder toute la
généralité nous allons symétriser $f\in H\otimes H$ et considérer des opérateurs sur $H\circ H$. Pour $n\geq1$, $H^{\circ^n}$ est naturellement muni de la base orthogonale formée des $\bigl(e_{j_1}\circ\cdots\circ e_{j_n}\bigl)$ $\bigl( (j_1,\ldots,j_n) \in (\N^*)^n\bigl)$. Dans cette base toute fonction $f$ de $H^{\circ^n}$ admet la décomposition suivante :
$$f=\sum_{\mathop{(j_1,\ldots,j_n)}\limits_{\in\bigl(\N^*\bigl)^n}}
\frac{\langle f, e_{j_1}\circ\cdots\circ e_{j_n}\rangle}{\Vert e_{j_1}\circ\cdots\circ e_{j_n}\Vert}
 e_{j_1}\circ\cdots\circ e_{j_n}$$
De façon à retrouver les résultats que l'on connaît ({\bf[10]}) dans le cas  gaussien nous
allons introduire les polynômes de Hermites associés à la mesure de Lebesgue sur $\R$
pondérée par $e^{-\frac{1}{2}x^2}$ et adapter l'opérateur d'anihilation
introduit par P.A Meyer ({\bf[10]}).\\
\newline
Nous allons utiliser les notations suivantes :
\begin{enumerate}
\item[$\bullet$]
$H_n$ : le $n^{\mbox{ième}}$ polynôme de Hermite, avec $H_0=1$
\item[$\bullet$]
On rappelle la formule de Rodriguez des polynômes de Hermite :
 $$H_n(x)=(-1)^ne^{\frac{1}{2}x^2}\frac{d^n\bigl(e^{-\frac{1}{2}x^2} \bigl)}{dx^n}.$$
\item[$\bullet$]
Dans la base des polynômes de Hermite on écrira pour $n\in\N$ :
$$ X^n=\sum_{k=0}^{n}\Gamma_{n,k}H_k(X).$$
\item[$\bullet$]
Dans la base des polynômes associés à la mesure $\mu$ on écrira pour $n\in\N$ :
$$X^n=\sum_{k=0}^{n}\gamma_{n,k}P_k(X).$$
\item[$\bullet$]
Nous pouvons donner la table des premiers c\oe fficients $\Gamma_{n,k}$ et
$\gamma_{n,k}$ :
$$
\begin{array}{|*{6}{c|}}\hline
\Gamma_{n,k} & k=0&k=1 &k=2 &k=3 & k=4\\ \hline
n=0 & 1 &0 &0 &0 &0\\ \hline
n=1 & 0 &1 &0 &0 &0\\ \hline
n=2 & 1 &0 &1 &0 &0\\ \hline
n=3 & 0 &3 &0 &1 &0\\ \hline
n=4 & 3 &0 &6 &0 &1\\ \hline
\end{array}$$
$$\begin{array}{|*{5}{c|}}\hline
\gamma_{n,k} & k=0&k=1 &k=2 &k=3\\ \hline
n=0 & 1 &  0  & 0 & 0\\ \hline
n=1 & 0 &  1  & 0 & 0\\ \hline
n=2 & 1 & m_3 & 1 & 0\\ \hline
\end{array}$$\\
Les premiers polynômes :\\
Hermite$\left\{\begin {array}{lcl}
             H_0 & = & 1\\
             H_1 & = & X\\
             H_2 & = & X^2-1\\
             H_3 & = & X^3-3X
     \end{array}\right.$
Polynômes associés à $\mu$
$\left\{\begin {array}{lcl}
             P_0 & = & 1\\
             P_1 & = & X\\
             P_2 & = & X^2-m_3X-1\\
     \end{array}\right.$

\end{enumerate}
Nous pouvons maintenant définir nos deux opérateurs :
\begin{de}
Fixons
$\left\{\begin{array}{lll}
\bullet & n\in\N^*\\
\bullet & (r,k)\in\{1,\ldots,n\}^2\\
\bullet & (\alpha_1,\ldots,\alpha_r)\in\N^n\mbox{\ \ tels que \ }\alpha_1+\cdots+\alpha_r=n\\
\bullet &(j_1,\ldots,j_r)\in\N^n\mbox{\ \ tels que \ }1\leq j_1<\cdots<j_r
\end{array}\right.$\\
On définit les deux opérateurs qui suivent :
\begin{enumerate}
\item
$\Phi^{\circ^n} : \begin{array}{ccl}
H^{\circ^n} & \to          & L^2(\nu)\\
\mathop{\bigcirc}\limits_{i=1}^r e_{j_i}^{\circ^{\alpha_i}}
                 & \mapsto & P_{\alpha_1}(X_{j_1})\cdot\ldots\cdot P_{\alpha_r}(X_{j_r})
\end{array}$
\item
$a_k^n :\begin{array}{lcl}
      H^{\circ^n} & \to & H^{\circ^{n-k}}\\
\dps\mathop{\bigcirc}\limits_{i=1}^r e_{j_i}^{\circ^{\alpha_i}}  &\mapsto  &
\dps\sum_{\mathop{(k_1+\cdots+k_r=k)}\limits_{k_i\geq0}}
\prod_{(k_i\not=0)}\Bigl[\Bigl(
\gamma_{\alpha_i,\alpha_i-k_i} - \Gamma_{\alpha_i,\alpha_i-k_i}
\Bigl)1_{[\alpha_i\geq k_i]}\Bigl]
\dps\mathop{\bigcirc}\limits_{p=1}^re_{j_p}^{\alpha_p-k_p}
\end{array}$
\end{enumerate}
\end{de}
\subsection*{Remarques}
\begin{itemize}
\item
En passant sur l'espace $H^{\circ^n}$ nous n'avons plus une base
normale. Cela va nous obliger à transformer le produit scalaire
pour conserver l'isométrie que nous proposait $\Phi$.
\item
On peut également remarquer que pour n=1, $\Phi^{\circ^n}$ n'est
autre que $\Phi$.
\end{itemize}
\subsection{Première propriété du modèle généralisé}
Nous allons regarder ici des propriétés de nature topologique.
Nous commencerons naturellement par préciser les normes sur
chacune des espaces entrant en action.\vspace{0,5cm}
\begin{itemize}
\item Nous prendrons le produit scalaire usuel sur $L^2(\nu)$ ainsi que
sur $H^{\otimes^n}$.\vspace{0,5cm}
\item
On va considérer le produit scalaire $\langle\, , \, \rangle_A^{H^{\circ^n}}$
sur $H^{\otimes^n}$ défini par :
$$\langle\, , \, \rangle_A^{H^{\circ^n}}:=n!\langle A\cdot\, ,A\cdot\,
\rangle_{H^{\otimes^n}}.$$
$\begin{array}{lccc}
\mbox{Où\, }A: &H^{\circ^n} & \to      & H^{\circ^n}\\
   & \mathop{\bigcirc}\limits_{i=1}^{r} e_{j_i}^{\alpha_i}
                & \mapsto  &
 \sum A_{ { j_1,\ldots,j_r}\atop{\alpha_1,\ldots,\alpha_r}}
 \mathop{\bigcirc}\limits_{i=1}^{r}e_{j_i}^{\alpha_i}
\end{array}$\\
avec
${A}_{ { j_1,\ldots,j_r}\atop{\alpha_1,\ldots,\alpha_r}}:=
\frac{\dps\prod_{i=1}^r E\left\vert P_{\alpha_i}(X_{j_i})\right\vert^2}
     {\dps\prod_{i=1}^r\alpha_i!}$.
\newline
Cet opérateur apparaît de façon naturelle dans la recherche d'une
isométrie.
\end{itemize}
\subsection*{Notations}
\begin{itemize}
\item Posons
\begin{itemize}
\item
$C_{ {k_1,\ldots,k_r}\atop{ \alpha_1,\ldots,\alpha_r}}:=
\dps\prod_{k_i\not=0}\left[\left
(\gamma_{\alpha_i,\alpha_i-k_i}-\gamma_{\alpha_i,\alpha_i-k_i}
\right)1_{[\alpha_i\geq k_i]}\right]$
\item
$C_{(k,n)}:=
\mathop{\sup}\limits_{k_1+\cdots+k_r=k\atop{\alpha_1+\cdots+\alpha_r=n}}
C_{ { \alpha_1,\ldots,\alpha_r}\atop{k_1,\ldots,k_r}}$
\item
$A_{(k,n)}:=
\mathop{\sup}\limits_{k_1+\cdots+k_r=k\atop{\alpha_1+\cdots+\alpha_r=n}}
A_{ {k_1,\ldots,k_r}\atop{ \alpha_1,\ldots,\alpha_r}}$
\end{itemize}
\item
On prend le produit scalaire usuel sur $H^{\circ^n}$, à savoir :
$$\Vert\;\cdot\; \Vert_{H^{\circ^n}}:=n!\Vert\;\cdot\; \Vert_{H^{\otimes^n}}$$
\item
Les normes des opérateurs sont :
\begin{itemize}
\item
$\Vert T(f)\Vert_{A}:= \mathop{\sup}\limits_{\Vert f\Vert\not=0}
\dps\frac{\Vert T(f)\Vert_{L^2(\nu)}}
         {\Vert f \Vert_{A}}$
\item
$\Vert T(f)\Vert_{H^{\circ^n}}:=
\mathop{\sup}\limits_{\Vert f\Vert\not=0}
\dps\frac{\Vert T_k(f)\Vert_{H^{\circ^{n-k}}}}
         {\Vert f\Vert_{H^{\circ^{n-k}}}}$
\end{itemize}
\end{itemize}
\begin{pp}
On peut étendre par linéarité et continuité les opérateurs $\Phi^{\circ^n}$
et $a_k^n$ à toute fonction $f$ de $H^{\circ^n}$.\\
De plus :
\begin{itemize}
\item
$\Vert \Phi^{\circ^n}\Vert_A=1$
\item
$\Vert a_k^n\Vert_A\leq C_{(n,k)}A_{(n,k)}\Bigl(_{\;\;\;\;k}^{k+r+1}\Bigl)
\dps\sum_{k_1+\cdots+k_r=k}\frac{\prod(\alpha_i-k_i)!1_{[\alpha_i\geq k_i]}}
{\prod\alpha_i!}$
\end{itemize}
\end{pp}
\subsection*{Remarque}
Le premier point de cette proposition nous dit, entre autre, que  $\Phi^{\circ^n}$ est une
isométrie de $(H^{\circ^n},\Vert\;\Vert_A)$ dans $L^2(\nu)$.
\subsection*{\underline{\textbf{Démonstration}}}
\begin{enumerate}
\item[]
\small
\begin{itemize}
\item
$\Vert \Phi^{\circ^n}(f)\Vert_{L^2(\nu)}^2=
\dps\sum_{r=1}^n\sum_{1\leq j_1<\cdots<j_r}
\langle f , \mathop{\bigcirc}\limits_{i=1}^re_{j_i}^{\alpha_i}\rangle^2
\prod_{i=1}^rE\left(P_{\alpha_i}\bigl(X_{j_i}\bigl)^2\right)$
\begin{itemize}
\item
$\Vert f\Vert_{A}^2=\dps\sum_{r=1}^{n}
\langle f , \mathop{\bigcirc}\limits_{i=1}^re_{j_i}^{\alpha_i}\rangle^2
A_{{j_1,\ldots,j_r}\atop{\alpha_1,\ldots,\alpha_r}}
n!\Vert \mathop{\bigcirc}\limits_{i=1}^re_{j_i}^{\alpha_i}
\Vert_{H^{\otimes^n}}^2$
\item
$\Vert \mathop{\bigcirc}\limits_{i=1}^re_{j_i}^{\alpha_i}
\Vert_{H^{\otimes^n}}^2=\frac{1}{n!}\prod_{i=1}^n\alpha_i!$
(confère {\bf[10]})
\end{itemize}
d'où l'isométrie.
\item
$\left\Vert a_k^n\left(\mathop{\bigcirc}\limits_{i=1}^re_{j_i}^{\alpha_i}\right)
\right\Vert_{H^{\circ^n}}^2=
\dps\sum_{{k_1+\cdots+k_r=k}\atop{k_i\geq0}}
\left\vert
C_{{\alpha_1+\cdots+\alpha_r}\atop{k_1+\cdots+k_r}}\right\vert^2
\left\Vert
\mathop{\bigcirc}\limits_{i=1}^re_{j_i}^{\alpha_i-k_i}\right\Vert_A$
\begin{itemize}
\item[]
$\leq C_{(k,n)}A_{(k,n)}
\dps\sum_{{k_1+\cdots+k_r=k}\atop{k_i\geq0}}
\dps\prod_{i=1}^{r}
\left[(\alpha_i-k_i)1_{[\alpha_i\geq k_i]}\right]!
\left(\begin{array}{c}
k+r-1\\
k
\end{array}\right)$
\item[]
$\leq C_{(k,n)}A_{(k,n)}(n-k)!
\sharp\{(k_1,\ldots,k_r)\in\N\;\vert\;k_1+\cdots+k_r=k\}$
\item[]
comme
$\left\Vert
\mathop{\bigcirc}\limits_{i=1}^re_{j_i}^{\alpha_i-k_i}
\right\Vert_{H^{\circ^n}}=\prod \alpha_i$ on a le résultat.
\end{itemize}
\end{itemize}
\end{enumerate} \normalsize$\triangle$
\subsection{Réécriture de l'intégrale}
On va ici regarder ce qu'il se passe, après un rapide rappel de la
définition de l'intégrale, pour nos opérateurs dans le cas où n=2.
\begin{enumerate}
\item
On a défini $I(h,f):=\varphi^{(1,1)}\bigl(h\otimes g1_C+
g\otimes h1_{\widetilde{C}}\bigl)+
\bigl(h\otimes g1_C\bigl)$, son écriture sous forme de série étant :
\newline
$I(f,g) = \dps\sum_{j\geq2}\sum_{k=1}^{j-1}
\langle h\otimes g1_C+g\otimes h1_{\widetilde{C}},e_j\otimes e_k
\rangle X_kX_j+
\dps\sum_{j\geq1}
\langle        h\otimes g1_C ,e_j^{\otimes^2}    \rangle (X_j^2-1)$\\
\hphantom{I(f,g) } =
$2\dps\sum_{j\geq2}\sum_{k=1}^{j-1}
\langle \bigl(h\otimes g1_C\bigl)^{\circ},e_j\otimes e_k
\rangle X_kX_j+
\dps\sum_{j\geq1}
\langle        h\otimes g1_C ,e_j^{\otimes^2}    \rangle
(X_j^2-1)$\\
d'une part
$$\langle
\bigl(h\otimes g1_C\bigl)^{\circ},e_j\otimes e_k
\rangle =
\langle
\bigl(h\otimes g1_C\bigl)^{\circ},e_k\otimes e_j
\rangle$$
donc
$$2\langle
\bigl(h\otimes g1_C\bigl)^{\circ},e_j\otimes e_k
\rangle =
\langle
\bigl(h\otimes g1_C\bigl)^{\circ},e_k\circ e_j
\rangle$$
d'autre part
$$\langle
\bigl(h\otimes g1_C\bigl)^{\circ},e_j\otimes e_j
\rangle =
\langle
\bigl(h\otimes g1_C\bigl)^{\circ},e_j\circ e_j
\rangle$$
\item
soit $f\in H^{\circ^2}$ et regardons $\Phi^{\circ^2}(f)$ :\\
$\Phi^{\circ^2}(f)=
\dps\sum_{1\leq j\not=k}
\langle
f^{\circ}, e_j\circ e_k
\rangle P_1(X_j)P_1(X_k)+
\dps\sum_{j\geq1}
\langle
f^{\circ},e_j^{\circ^2}
\rangle P_2(X_j)$\\
$= 2\dps\sum_{j\geq2}\sum_{k=1}^{j-1}
\langle
f^{\circ}, e_j\circ e_k
\rangle X_jX_k+
\dps\sum_{j\geq1}
\langle
f^{\circ},e_j^{\circ^2}
\rangle \bigl(X_j^2-m_3X_j-1\bigl)$\\
\item
Soit $f\in H^{\circ^2}$ et regardons $a_1^2$ :\\
$a_1^2\bigl(e_j\circ e_k \bigl) =
\left\{\begin{array}{ll}
0 & \mbox{si\;} j\not=k\\
\gamma_{2,1} & \mbox{si\;} j=k
\end{array}\right.$
avec $\gamma_{2,1}=m_3$
\end{enumerate}
Nous pouvons rassembler tous ces résultats dans la proposition
suivante :
\begin{pp}
$$I(f,g) = \Bigl(
\Phi^{\circ^2}+\Phi\circ a_1^2
\Bigl)
\Bigl(\left[h\otimes g1_{C}\right]^{\circ^2}\Bigl)$$
\end{pp}
\subsection*{Remarque}
Dans le cas gaussien on retrouve la décomposition de P.A. Meyer :
$$I(f,g) =
\Phi^{\circ^2}
\Bigl(\left[h\otimes g1_{C}\right]^{\circ^2}\Bigl)$$
\section{Intégration par partie}
Nous allons établir une formulation du type formule de Itô pour
les produits $\Phi(h)\Phi(g)$ $(h,g\in H)$.
\begin{thm}
Soient $h,g\in H$, on a :
$$\Phi(h)\Phi(g)=\int\Phi(h)_sd\Phi(g)_s+\int\Phi(g)_sd\Phi(h)_s
+\left[\Phi(h),\Phi(g)\right]$$
avec $\left[\Phi(h),\Phi(g)\right]:=\left(\Phi^{\circ^2}+a_1^2\right)
(h\otimes g1_{\Delta}) + \langle h,g\rangle$
où $\Delta:=\left\{(x,x)\;\vert\;x\in\R\right\}$
\end{thm}
\subsection*{\underline{\textbf{Démonstration}}}
\begin{enumerate}
\item[]
\small
Posons $\Delta_{k,k+1}\Phi(h):=\Phi(h)_{t_{k+1}}-\Phi(h)_{t_{k}}$\\
$\dps\sum_{j=0}^n \Delta_{j,j+1}\Phi(h)
     \sum_{k=0}^n\Delta_{k,k+1}\Phi(g)
=   \dps\sum_{j=0}^n\Phi\left(h1_{]t_j,t_{j+1}]}\right)
     \sum_{k=0}^n\Phi\left(g1_{]t_k,t_{k+1}]}\right)\\
\begin{array}{cl}

= &\dps\sum_{j=0}^n\Phi\left(h1_{]t_j,t_{j+1}]}\right)
                   \Phi\left(g1_{]t_j,t_{j+1}]}\right)\\
&+  \dps\sum_{j=1}^n\sum_{k=0}^{j-1}\Phi\left(g1_{]t_k,t_{k+1}]}\right)
                   \Phi\left(h1_{]t_j,t_{j+1}]}\right)\\
&+  \dps\sum_{k=1}^n\sum_{j=0}^{k-1}\Phi\left(h1_{]t_j,t_{j+1}]}\right)
\Phi\left(g1_{]t_k,t_{k+1}]}\right).\\
\end{array}$
\newline
Appelons $S_i^n$ $(i\in\{1,2,3\})$ la $i^{\mbox{ième}}$ somme du
membre de droite.
\begin{enumerate}
\item
On utilisant les calculs faits dans la démonstration de la
proposition 6.2.2.
$$S_1^n = \left(\Phi^{\circ^2}+\Phi\circ a_1^2\right)
\left(
h\otimes g \sum_{j=0}^n 1_{]t_j,t_{j+1}]}\otimes1_{]t_j,t_{j+1}]}
\right)$$
$$+\sum_{j=0}^n \langle
h1_{]t_j,t_{j+1}]},g 1_{]t_j,t_{j+1}]}
\rangle.$$
$\bullet$ D'une part :
$\dps\sum_{j=0}^n1_{]t_j,t_{j+1}]}\otimes1_{]t_j,t_{j+1}]}
\stackrel{L^2(\mu^{\otimes^2})}{\longrightarrow}1_{\Delta}$\\
\hphantom{0,5cm} et
$\langle
h\otimes g 1_{\Delta},e_{k_1}\otimes e_{k_2}
\rangle = 0.$\\
$\bullet$ D'autre part
$\dps\sum_{j=0}^n\langle
h1_{]t_j,t_{j+1}]},g1_{]t_j,t_{j+1}]}
\rangle = \langle h , g \rangle$.\\
Donc $S_1^n
\dps\mathop{\stackrel{L^2(\nu)}{\longrightarrow}}
\limits_{n\to\infty}
\langle h , g \rangle$
\item
$$S_2^n = \left(\Phi^{\circ^2}+\Phi\circ a_1^2\right)
\left(
g\otimes h \sum_{j=1}^n\sum_{k=0}^{j-1}
1_{]t_k,t_{k+1}]}\otimes1_{]t_j,t_{j+1}]}
\right)$$
$$+\sum_{j=1}^n \sum_{k=0}^{j-1}
\langle
h1_{]t_k,t_{k+1}]},g 1_{]t_j,t_{j+1}]}
\rangle.$$
$\bullet$ D'une part :
$\dps\sum_{j=1}^n\sum_{k=0}^{j-1}
1_{]t_k,t_{k+1}]}\otimes1_{]t_j,t_{j+1}]}
=\sum_{j=0}^n
1_{]0,t_{j}]}\otimes1_{]t_j,t_{j+1}]}
\stackrel{L^2(\mu^{\otimes^2})}{\longrightarrow}1_{C}$.\\
$\bullet$ D'autre part
$\langle
h1_{]t_k,t_{k+1}]},g1_{]t_j,t_{j+1}]}
\rangle = o$ car $j\not=k$.\\
Donc $S_2^n
\dps\mathop{\stackrel{L^2(\nu)}{\longrightarrow}}
\limits_{n\to\infty}
\left(
\Phi^{\circ^2}+\Phi\circ a_1^2
\right)
\left(
g\otimes h1_C
\right)=:\int\Phi(g)_sd\Phi(h)_s$
\item
$$S_3^n = \left(\Phi^{\circ^2}+\Phi\circ a_1^2\right)
\left(
h\otimes g \sum_{k=1}^n\sum_{j=0}^{k-1}
1_{]t_j,t_{j+1}]}\otimes1_{]t_k,t_{k+1}]}
\right)$$
on procède sensiblement de la même façon que dans le point
précédent et on a :
$$S_3^n
\dps\mathop{\stackrel{L^2(\nu)}{\longrightarrow}}
\limits_{n\to\infty}
\int\Phi(h)_sd\Phi(g)_s$$
\end{enumerate}
\end{enumerate} \normalsize$\triangle$
\section{Norme}
\begin{pp}
$$\left\Vert
\int\Phi(h)_sd\Phi(g)_s
\right\Vert_{L^2(\nu)}^2 =
\left\Vert
h\otimes g1_C
\right\Vert_{L^2(\nu)}^2
+\sum_{j\geq0}\langle h\otimes g,e_j\otimes e_j \rangle^2
\left(EX_j^4-3\right)$$
\end{pp}
\subsection*{Remarque}
Dans le cas gaussien (équidistribuées, normales, centrées) on a $EX_j^4=3$
et donc :
$$\left\Vert
\int\Phi(h)_sd\Phi(g)_s
\right\Vert_{L^2(\nu)}^2 =
\left\Vert
h\otimes g1_C
\right\Vert_{L^2(\nu)}^2.$$
Ce qui est conforme au résultat concernant l'intégrale de Ito.
\subsection*{\underline{\textbf{Démonstration}}}
\begin{enumerate}
\item[]
\small
Notons $f:=h\otimes g1_C$ alors :\\
$I(h,g)=\varphi^{(1,1)}(f+\widetilde{f}) +
\varphi^{(2)}(f)$\\
$$\begin{array}{ll}
E\left(\left\vert I(h,g)\right\vert^2\right)= &
E\left(\left\vert\varphi^{(1,1)}(f)\right\vert^2\right)+
E\left(\left\vert\varphi^{(1,1)}(\widetilde{f})\right\vert^2\right)+
E\left(\left\vert\varphi^{(2)}(f)\right\vert^2\right)\\
& +2\left[
E\left(\varphi^{(1,1)}(f)\varphi^{(1,1)}(\widetilde{f})\right)+
E\left(\varphi^{(1,1)}(f)\varphi^{(2)}(f)\right)\right.\\
& \left.+E\left(\varphi^{(1,1)}(\widetilde{f})\varphi^{(2)}(f)\right)
\right].
\end{array}$$
Les lemmes 6.1.3. et 6.1.4. nous donnent :\\
$$\begin{array}{ll}
E\left(\left\vert I(h,g)\right\vert^2\right)= &
\dps\sum_{j\geq2}\sum_{k=1}^{j-1}
\langle f, e_j\otimes e_k \rangle^2 +
\sum_{j\geq2}\sum_{k=1}^{j-1}
\langle \widetilde{f}, e_j\otimes e_k \rangle^2\\
& \dps+\sum_{j\geq1}
\langle f, e_j\otimes e_j \rangle^2\left(EX_j^4-1\right)\\
& \dps+2\sum_{j\geq2}\sum_{k=1}^{j-1}
\langle f, e_j\otimes e_k \rangle
\langle \widetilde{f}, e_j\otimes e_k \rangle.
\end{array}$$
$\bullet$ D'une part
$\langle \widetilde{f},e_j\otimes e_k\rangle =
\langle f,e_k\otimes e_j\rangle.$\\
$\bullet$ D'autre part :
$\left\Vert f\right\Vert_{H^{\otimes^2}} =
\dps\sum_{(j,k)\in\atop\bigl(\N^*\bigl)^2} \langle f,e_j\otimes e_k\rangle^2$\\
\hphantom{$\bullet$ D'autre part : $\left\Vert f\right\Vert_{H^{\otimes^2}}$}
$=\dps\sum_{j\geq2}\sum_{k=1}^{j-1}
\langle f,e_j\otimes e_k\rangle^2
+\sum_{j\geq2}\sum_{k=1}^{j-1}
\langle f,e_k\otimes e_j\rangle^2$\\
\hspace*{4,5cm}
$+\dps\sum_{j\geq1}\langle f,e_j\otimes e_j\rangle^2.$\\
On a donc, en ajoutant et retranchant le terme
$\dps\sum_{j\geq1}\langle f,e_j\otimes e_j\rangle^2$ :\\
$E\left\vert I(h,g)\right\vert^2=
\left\Vert h\otimes g1_C\right\Vert^2_{L^2(\mu^{\otimes^2})}+
\dps\sum_{j\geq1}\sum_{j\geq1}
\langle f,e_j\otimes e_j\rangle^2
\left(EX_j^4-2\right)$\\
\hphantom{$E\left\vert I(h,g)\right\vert^2=$}
$+2\dps\sum_{j\geq2}\sum_{k=1}^{j-1}
\langle f,e_j\otimes e_k\rangle
\langle \widetilde{f},e_j\otimes e_k\rangle.$\\
$\bullet$ D'une part
$2\dps\sum_{j\geq2}\sum_{k=1}^{j-1}
\langle f,e_j\otimes e_k\rangle
\langle \widetilde{f},e_j\otimes e_k\rangle =
\sum_{j\not=k}\langle f,e_j\otimes e_k\rangle
\langle \widetilde{f},e_j\otimes e_k\rangle.$\\
$\bullet$ D'autre part
$\langle f,\widetilde{f}\rangle=
\dps\sum_{(j,k)\in\atop\bigl(\N^*\bigl)^2}
\langle f,e_j\otimes e_k\rangle
\langle \widetilde{f},e_j\otimes e_k\rangle.$\\
En ajoutant et retranchant
$\langle f,e_j\otimes e_j\rangle
\langle \widetilde{f},e_j\otimes e_j\rangle =
\dps\sum_{j\geq1}\langle f,e_j\otimes e_j\rangle^2$ on a donc :
$$E\left\vert I(h,g)\right\vert^2 =
\left\Vert h \otimes g1_C\right\Vert^2_{L^2(\mu^{\otimes^2})}
+ \sum_{j\geq1}
\langle f,e_j\otimes e_j\rangle\left(EX_j^4-3\right) +
\langle f ,\widetilde{f}\rangle.$$
Pour finir il suffit d'expliciter $\langle f ,\widetilde{f}\rangle$
et remarquer que $1_C\cdot1_{\widetilde{C}}= 0$
\end{enumerate} \normalsize$\triangle$
\section{Isométries partielles}
\subsection{Notations}
\begin{enumerate}
\item
Pour $f\in H\otimes H$ on notera :
$J^2(f):=\left(\Phi^{\circ^2}+\Phi\circ a_1^2\right)(f).$
\item
$S_2:=\left\{f\in H\otimes H\; \vert\;f(x,y)=f(y,x)\right\}$ le
sous-espace de $H\otimes H$ des fonctions symétriques.
\item On utilisera l'opérateur :
$$\begin{array}{llcl}
B : & H\otimes H  &\to  & H\otimes H\\
&e_j\otimes e_k &\mapsto &
              \left\{\begin{array}{ll}
              e_j\otimes e_k &j\not=k\\
              E(X_j^2-1)^2e_j\otimes e_j &j=k
              \end{array}\right.
\end{array}$$
\item
$H^{(1,1)}:=\mbox{Vect}\left\{e_j\otimes
e_j\;\vert\;j\in\N^*\right\}$\\
$H^{(2)}:=\mbox{Vect}\left\{e_j\otimes
e_k\;\vert\;j,k\in\N^*,\;j\not=k\right\}$
\end{enumerate}
\subsection{Isométries}
\begin{pp} Les applications
$$\varphi^{(1,1)} \;:\;\left(H^{(1,1)},\;\; \Vert\;\Vert_B\right)
\to L^2(\nu) $$
$$\varphi^{(2)} \;:\;\left(H^{(2)},\; \Vert\;\;\Vert_B\right)
\to L^2(\nu) $$
sont des isométries.
\end{pp}
\subsection*{\underline{\textbf{Démonstration}}}
\begin{enumerate}
\item[]
\small
$\bullet$ Soit $f\in H^{(1,1)}$ on a :
$$\left\Vert\varphi^{(1,1)}(f)\right\Vert_{L^2(\nu)}^2 =
\sum_{j\not=k}\langle f , e_j\otimes e_k\rangle^2 =
\left\Vert f\right\Vert_{B}$$
$\bullet$ Soit $f\in H^{(2)}$ on a :
$$\begin{array}{lcl}
\left\Vert\varphi^{(2)}(f)\right\Vert_{L^2(\nu)}^2 & = &
\dps\sum_{j\geq1}\langle f , e_j\otimes e_j\rangle^2E(X_j^2-1)^2\\
& = &\dps\sum_{j\geq1}\frac{\langle f , B(e_j\otimes e_j)\rangle^2}
                  {\Vert B(e_j\otimes e_j)\Vert_{H\otimes H} }\\
& = &\dps\sum_{j\geq1}\frac{\langle f , e_j\otimes e_j\rangle_B^2}
                  {\Vert e_j\otimes e_j\Vert_{B} }\\
& = &\left\Vert f\right\Vert_{B}
\end{array}$$
\end{enumerate} \normalsize$\triangle$
\subsection{Isométrie sur $S^2$}
Dans le cas gaussien $J_2=\Phi^{(\circ^2)}$ est une isométrie sur
$S^2$. Dans le cas général nous perdons l'isométrie mais
conservons une relation entre les normes de $L^2(\nu)$ et $S^2$ :
\begin{pp}
$$\Vert\;\;\Vert_{H\circ H}\;\mbox{et}\; \Vert\;\;\Vert_{L^2(\nu)}\circ J_2$$
sont des normes équivalentes sur $S^2$.
\end{pp}
\subsection*{\underline{\textbf{Démonstration}}}
\begin{enumerate}
\item[]
\small
Soit $f \in S^2$, des calculs déjà faits dans la proposition
6.5.1. nous donnent :
$$\left\Vert J_2(f)\right\vert^2_{L^2(\nu)} =
2\sum_{j\not=k}\langle f, e_j\circ e_k\rangle^2 +
\sum_{j\geq1}\langle f, e_j\circ e_j\rangle^2. $$
On pose :
$\left\{\begin{array}{lcl}
a &:= &\min\{E(X_j^2-1)^2\;(j\in\N^*),\;2\}\\
b &:= &\max\{E(X_j^2-1)^2\;(j\in\N^*),\;2\}\\
\end{array}\right.$\\
on trouve :
$$a\left\Vert f\right\Vert^2_{S^2}\leq
\left\Vert J_2(f)\right\Vert^2_{L^2(\nu)}\leq
b\left\Vert f\right\Vert^2_{S^2}$$
\end{enumerate} \normalsize$\triangle$
\subsection*{Remarque}
L'ensemble de ces calculs fournit l'opérateur qu'il faut utiliser
pour obtenir une isométrie. Les opérateurs A et B introduits
précedemment ne conviennent pas. On va donc en introduire un
troisième, C, sur $H\otimes H$ :
$$C(e_j\otimes e_k) := \left\{\begin{array}{ll}
2e_j\otimes e_k & j\not=k\\
E(X_j^2-1)^2e_j\otimes e_j & j=k
\end{array}\right.$$
\begin{pp}
$$J_2 : \left(S^2,\;\Vert\;\;\Vert_C\right)\to L^2(\nu)$$
est une isométrie.
\end{pp}
\subsection*{\underline{\textbf{Démonstration}}}
\begin{enumerate}
\item[]
\small
Il suffit de  remarquer que $\left\Vert e_j\circ e_k\right\Vert_C =
\left\{\begin{array}{ll}

2& j\not=k\\
E(X_j^2-1)^2 &j=k
\end{array}\right.$\\
et $\left\Vert e_j\circ e_k\right\Vert_{L^2(\nu)} =
\dps\sum_{(j,k)\in\atop\bigl(\N^*\bigl)^2}
\frac {\langle f, e_j\circ e_k\rangle_C^2}
      {\Vert e_j\circ e_k\Vert_C}$
\end{enumerate} \normalsize$\triangle$
\subsection{Vers une généralisation}
Pour des fonctions $h, g\in F$, nous pouvons définir un processus
$\left(Z_t\right)_{t>0}$ de la façon suivante :
$$Z_t:=\int\Phi(h)_sd\Phi\left(g 1_{]0,t]}\right)_s.$$
Nous noterons plutôt :
$Z_t:=\int_0^t\Phi(h)_sd\Phi\left(g \right)_s$.\\
Nous allons étudier l'existence de modifications continues puis la
variation quadratique de ce processus.\\
Au préalable, nous établirons des résultats techniques utiles pour
mener à bien ces études.
\chapter{Résultats techniques}
\section{Etude de $Z_t^2$}
Pour simplifier l'écriture des calculs nous allons expliciter les
calculs de $\left(J_2(f)\right)^2$ pour $f\in H\circ H$ et
noterons :
$\left\{\begin{array}{lcl}
a_{jk} & := &\langle f,e_j\circ e_k\rangle\\
a_j    & := &a_{jj}
\end{array}\right.$
\begin{lm} On a pour $f\in H\circ H$\\
$\begin{array}{ccl}
\left(J_2(f)\right)^2 & = &\\
\mbox{(ordre 0)} & &4\dps\sum_{j_2<j_1}a_{j_1j_2}^2 +
(\gamma_{4,0}-1)\sum_{j\geq1}a_j^2\\
\mbox{(ordre 1)}& +&
4\gamma_{3,0}\dps\sum_{j_2<j_1}
a_{j_1j_2}a_{j_1}P_1\left(X_{j_2}\right)\\
 & + &
4\gamma_{3,0}\dps\sum_{j_2<j_1}
a_{j_1j_2}a_{j_2}P_1\left(X_{j_1}\right)\\
 & + &
(\gamma_{4,1}-2m_3)\dps\sum_{j\geq1}a_{j}^2P_1\left(X_{j}\right)\\
\mbox{(ordre 2)}& +&
8\dps\sum_{j_3<j_2<j_1}
a_{j_1j_2}a_{j_1j_3}P_1\left(X_{j_2}\right)P_1\left(X_{j_3}\right)\\
& +&
8\dps\sum_{j_3<j_2<j_1}
a_{j_1j_2}a_{j_2j_3}P_1\left(X_{j_1}\right)P_1\left(X_{j_3}\right)\\
& +&
8\dps\sum_{j_3<j_2<j_1}
a_{j_2j_3}a_{j_1j_3}P_1\left(X_{j_1}\right)P_1\left(X_{j_2}\right)\\
& +&
4\dps\sum_{j_2<j_1}
a_{j_1j_2}^2\left(P_2\left(X_{j_1}\right)+P_2\left(X_{j_2}\right)\right)\\
& +&
\dps\sum_{j_2<j_1}
\left[4m_2^2a_{j_1j_2}^2 + 2m_3^2a_{j_1}a_{j_2}+4(\gamma_{3,1}-1)\left(a_{j_1j_2}a_{j_1}+a_{j_1j_2}a_{j_2}\right)\right]
P_1\left(X_{j_1}\right)P_1\left(X_{j_2}\right)\\
& +&
\left(\gamma_{4,2}-2\right)\dps\sum_{j\geq1}
a_{j}^2P_2\left(X_{j}\right)\\
\end{array}\\
\begin{array}{ccl}
\mbox{(ordre 3)}& +&
m_3\dps\sum_{j_3<j_2<j_1}\left[
8\left(a_{j_1j_2}a_{j_1j_3}+a_{j_1j_2}a_{j_2j_3}+a_{j_1j_3}a_{j_2j_3}
\right)+\right.\\
&&\left.4\left(a_{j_1}a_{j_2j_3}+a_{j_1j_3}a_{j_2}+a_{j_1j_2}a_{j_3}\right)\right]
P_1\left(X_{j_1}\right)P_1\left(X_{j_2}\right)P_1\left(X_{j_3}\right)\\
& +&
\dps\sum_{j_2<j_1}\left(
4m_3a_{j_1j_2}^2+2m_3a_{j_1}a_{j_2}+4\gamma_{3,2}a_{j_1j_2}a_{j_1}\right)
P_2\left(X_{j_1}\right)P_1\left(X_{j_2}\right)\\
& +&
\dps\sum_{j_2<j_1}\left(
4m_3a_{j_1j_2}^2+2m_3a_{j_1}a_{j_2}+4\gamma_{3,2}a_{j_1j_2}a_{j_2}\right)
P_1\left(X_{j_1}\right)P_2\left(X_{j_2}\right)\\
& +&
\gamma_{4,3}\dps\sum_{j\geq1}
a_{j}^2P_3\left(X_{j}\right)\\
\mbox{(ordre 4)}& +&
8\dps\sum_{j_4<j_3<j_2<j_1}\left(
a_{j_1j_2}a_{j_3j_4}+a_{j_1j_3}a_{j_2j_4}+a_{j_2j_3}a_{j_1j_4}\right)
\prod_{i=1}^4P_1\left(X_{j_i}\right)\\
& +&
\dps\sum_{j_3<j_2<j_1}\left(
8a_{j_1j_2}a_{j_1j_3}+4a_{j_1}a_{j_2j_3}\right)
P_2\left(X_{j_1}\right)P_1\left(X_{j_2}\right)P_1\left(X_{j_3}\right)\\
& +&
\dps\sum_{j_3<j_2<j_1}\left(
8a_{j_1j_2}a_{j_2j_3}+4a_{j_1j_3}a_{j_2}\right)
P_1\left(X_{j_1}\right)P_2\left(X_{j_2}\right)P_1\left(X_{j_3}\right)\\
& +&
\dps\sum_{j_3<j_2<j_1}\left(
8a_{j_1j_3}a_{j_2j_3}+4a_{j_1j_2}a_{j_3}\right)
P_1\left(X_{j_1}\right)P_1\left(X_{j_2}\right)P_2\left(X_{j_3}\right)\\
& +&
\dps\sum_{j_2<j_1}
4a_{j_1j_2}a_{j_1}
P_3\left(X_{j_1}\right)P_1\left(X_{j_2}\right)\\
& +&
\dps\sum_{j_2<j_1}
4a_{j_1j_2}a_{j_2}
P_1\left(X_{j_1}\right)P_3\left(X_{j_2}\right)\\
& +&
\dps\sum_{j_2<j_1}\left(
4a_{j_1j_2}^2+2a_{j_1}a_{j_2}\right)
P_2\left(X_{j_1}\right)P_2\left(X_{j_2}\right)\\
& +&
\dps\sum_{j\geq1}
a_{j}^2P_4\left(X_{j_1}\right)\\
\end{array}$
\end{lm}
\subsection*{\underline{\textbf{Démonstration}}}
\begin{enumerate}
\item[]
\small
On peut écrire :
$$\begin{array}{lcl}
\left(\Phi^{\circ^2}(f)+\phi( a_1^2f)\right)^2 & = &
4\left(\dps\sum_{j_2<j_1}
\langle f , e_{j_1}\circ e_{j_3}\rangle X_{j_1}X_{j_2}\right)^2\\
&&+\left(\dps\sum_{j\geq1}
\langle f , e_{j}\circ e_{j}\rangle (X_{j}^2-1)\right)^2\\
&&+\dps\sum_{j_2<j_1}
\langle f , e_{j_1}\circ e_{j_3}\rangle X_{j_1}X_{j_2}
\sum_{j_3\geq1}
\langle f , e_{j_3}\circ e_{j_3}\rangle (X_{j_3}^2-1)
\end{array}$$
On va procéder en quatre étapes, une pour le calcul de chacun des
termes du membre de droite et une dernière dans laquelle nous
sommerons ces trois résultats.
\subsubsection*{\underline{\bf Etape 1 :} calcul de $\left(\dps\sum_{j_2<j_1}
\langle f , e_{j_1}\circ e_{j_3}\rangle X_{j_1}X_{j_2}\right)^2=:M^2$}
$M^2=\dps\sum_{j_1\geq2}^n  \left( \sum_{k_1=1}^{j_1-1}
a_{j_1k_1}X_{k_1}\right)X_{j_1}
\sum_{j_2\geq2}^n  \left( \sum_{k_2=1}^{j_2-1}
a_{j_2k_2}X_{k_2}\right)X_{j_2}$\\
\hphantom{$M^2$}$ = \dps\sum_{0\leq k_1<j_1\atop0\leq k_2<j_2}
a_{j_1k_1}a_{j_2k_2}
X_{k_1}X_{j_1}X_{k_2}X_{j_2}$\\
On va réordonner la somme en ordonnant les paramètres.\\
On fixe $k_1<j_1$, on a cinq cas à voir :
\begin{enumerate}
\item $j_2<j_1$ \\
Il y a cinq sous cas
\begin{enumerate}
\item $j_2<k_2<j_1<k_1$ 
\item $j_2<k_2=j_1<k_1$  
\item $j_2<j_1<k_2<k_1$  
\item $j_2<j_1<k_2=k_1$   
\item $j_2<j_1<k_1<k_2$    
\end{enumerate}
\item $j_2=j_1$\\  
Il y a trois sous cas :
\begin{enumerate}
\item $j_1=j_2<k_2<k_1$  
\item $j_1=j_2<k_1=k_2$  
\item $j_1=j_2<k_1<k_2$  
\end {enumerate}
\item $j_1<j_2<k_1$\\  
Il y a trois sous cas :
\begin{enumerate}
\item $j_1<j_2<k_2<k_1$  
\item $j_1<j_2<k_2=k_1$  
\item $j_1<j_2<k_1<k_2$  
\end{enumerate}
\item $j_1<j_2=k_1$\\  
Il y a un sous cas :
\begin{enumerate}
\item[] $j_1<j_2=k_1<j_2$  
\end{enumerate}
\item $j_1<k_1<k_2$\\ 
Il y a un sous cas :
\begin{enumerate}
\item[] $j_1<k_1<j_2<k_2$  
\end{enumerate}
\end{enumerate}
Nous allons maintenant réécrire la somme en prenant chacune
de ces indexations et en  ordonnant les sommes par
ordre décroissant du nombre d'indices :\\
\newline
$$\begin{array}{lclc}
M^2 &= &\dps\sum_{0\leq j_2<k_2<j_1<k_1}
a_{j_1k_1}a_{j_2k_2}X_{k_1}X_{k_2}X_{j_1}X_{j_2}  &\mbox{(a-i)}\\
(4\; ind.)& + &\dps\sum_{0\leq j_2<j_1<k_2<k_1}
a_{j_1k_1}a_{j_2k_2}X_{k_1}X_{k_2}X_{j_1}X_{j_2}  &\mbox{(a-iii)}\\
& +&  \dps\sum_{0\leq j_2<j_1<k_1<k_2}
a_{j_1k_1}a_{j_2k_2}X_{k_1}X_{k_2}X_{j_1}X_{j_2}  &\mbox{(a-v)}\\
& +&  \dps\sum_{0\leq j_1<j_2<k_2<k_1}
a_{j_1k_1}a_{j_2k_2}X_{k_1}X_{k_2}X_{j_1}X_{j_2}  &\mbox{(c-i)}\\
& +&  \dps\sum_{0\leq j_1<j_2<k_1<k_2}
a_{j_1k_1}a_{j_2k_2}X_{k_1}X_{k_2}X_{j_1}X_{j_2}  &\mbox{(c-iii)}\\
& +&  \dps\sum_{0\leq j_1<k_1<j_2<k_2}
a_{j_1k_1}a_{j_2k_2}X_{k_1}X_{k_2}X_{j_1}X_{j_2}  &\mbox{(e)}\\
(3 \;ind.)& +&  \dps\sum_{0\leq j_2<k_2<k_1}
a_{k_1k_2}a_{j_2k_2}X_{k_2}^2X_{j_2}X_{j_1}  &\mbox{(a-ii)}\\
& +&  \dps\sum_{0\leq j_2<j_1<k_1}
a_{j_1k_1}a_{j_2k_1}X_{k_1}^2X_{j_1}X_{j_2}  &\mbox{(a-iv)}\\
& +&  \dps\sum_{0\leq j_1<k_2<k_1}
a_{j_1k_2}a_{k_1j_1}X_{j_1}^2X_{k_2}X_{k_1}  &\mbox{(b-i)}\\
& +&  \dps\sum_{0\leq j_1<k_1<k_2}
a_{j_1k_2}a_{k_1j_1}X_{j_1}^2X_{k_1}X_{k_2}  &\mbox{(b-iii)}\\
& +&  \dps\sum_{0\leq j_1<j_2<k_1}
a_{j_1k_1}a_{j_2k_2}X_{k_1}^2X_{j_1}X_{j_2}  &\mbox{(c-ii)}\\
& +&  \dps\sum_{0\leq j_1<j_2<k_2}
a_{j_1j_2}a_{j_2k_2}X_{j_2}^2X_{j_1}X_{k_2}  &\mbox{(d)}\\
(2\; ind.)& +&  \dps\sum_{0\leq k<j}
a_{jk}^2X_{j}^2X_{j}^2  &\mbox{(b-ii)}\\
\end{array}$$\\
\newline
On range les indices $j_4<j_3<j_2<j_1$, $j_3<j_2<j_1$ et
$j_2<j_1$, puis on factorise :
\newline
$$\begin{array}{lcl}
M^2 & = & \dps\sum_{j_4<j_3<j_2<j_1}
\left(a_{j_1j_2}a_{j_3j_4} + a_{j_1j_3}a_{j_2j_4}+a_{j_2j_3}a_{j_1j_4}+
a_{j_1j_4}a_{j_2j_3}+a_{j_2j_4}a_{j_1j_3}\right.\\
&&\left.\hspace{2cm}+a_{j_3j_4}a_{j_1j_2}
\right)X_{j_4}X_{j_3}X_{j_2}X_{j_1}\\
& +&\dps\sum_{j_3<j_2<j_1}
\left(a_{j_1j_2}a_{j_1j_3} + a_{j_1j_3}a_{j_1j_2}
\right)X_{j_3}X_{j_2}X_{j_1}^2\\
& +&\dps\sum_{j_3<j_2<j_1}
\left(a_{j_1j_2}a_{j_2j_3} + a_{j_2j_3}a_{j_1j_2}
\right)X_{j_3}X_{j_2}^2X_{j_1}\\
& +&\dps\sum_{j_3<j_2<j_1}
\left(a_{j_2j_3}a_{j_1j_3} + a_{j_1j_3}a_{j_2j_3}
\right)X_{j_3}^2X_{j_2}X_{j_1}\\
& +&\dps\sum_{j_2<j_1}
a_{j_1j_2}^2X_{j_3}X_{j_2}^2X_{j_1}^2\\
\end{array}$$
\newline
Il ne nous reste plus qu'à exprimer ces résultats dans la base des
polynômes $\left(P_n\right)_n$ associée à $\mu$.\\
\newline
On rappel les relations suivantes :
\newline
$$\left\{\begin{array}{lcl}
1     & = & P_0\\
X_j   & = & P_1(X_j)\\
X_j^2 & = & P_2(X_j)+m_3P_1(X_j)+1\\
X_j^n & = & \dps\sum_{k=0}^n\gamma{n,k}P_k(X_j)
\end{array}\right.$$
\newline
On remplace, on calcule et on trouve :\\
\newline
$$\left\vert\begin{array}{lcl}\hline
M^2 & = &\\
(ordre\; 0)&& \dps\sum_{j_2<j_1}
a_{j_1j_2}^2\\
(ordre \;1) & + &
m_3\dps\sum_{j_2<j_1}a_{j_1j_2}^2\left(P_1(X_{j_1})+P_1(X_{j_2})\right)\\
(ordre \;2) & + &
2\dps\sum_{j_3<j_2<j_1}a_{j_1j_2}a_{j_1j_3}P_1(X_{j_2})P_1(X_{j_3})\\
& + &
2\dps\sum_{j_3<j_2<j_1}a_{j_1j_2}a_{j_2j_3}P_1(X_{j_1})P_1(X_{j_3})\\
& + &
2\dps\sum_{j_3<j_2<j_1}a_{j_2j_3}a_{j_1j_3}P_1(X_{j_1})P_1(X_{j_2})\\
& + &
\dps\sum_{j_2<j_1}a_{j_1j_2}^2\left(P_2(X_{j_1})+P_2(X_{j_2})\right)\\
& + &
m_3^2\dps\sum_{j_2<j_1}a_{j_1j_2}^2P_1(X_{j_1})P_1(X_{j_2})\\
(ordre \;3)& + &
2m_3\dps\sum_{j_3<j_2<j_1}\left(
a_{j_1j_2}a_{j_1j_3}+a_{j_1j_2}a_{j_2j_3}\right.\\
&&\hspace{2cm}\left.+a_{j_2j_3}a_{j_1j_3}\right)
P_1(X_{j_1})P_1(X_{j_2})P_1(X_{j_3})\\
& + &
m_3\dps\sum_{j_2<j_1}a_{j_1j_2}^2
\left(P_2(X_{j_1})P_1(X_{j_2})+P_2(X_{j_2})P_1(X_{j_1})\right)\\
(ordre \;4)& + &
2\dps\sum_{j_4<j_3<j_2<j_1}\left(
a_{j_1j_2}a_{j_3j_4}+a_{j_1j_3}a_{j_2j_4}\right.\\
&&\hspace{2cm}\left.+a_{j_2j_3}a_{j_1j_4}\right)
P_1(X_{j_1})P_1(X_{j_2})P_1(X_{j_3})P_1(X_{j_4})\\
& + &
2\dps\sum_{j_3<j_2<j_1}a_{j_1j_2}a_{j_1j_3}
P_2(X_{j_1})P_1(X_{j_2})P_1(X_{j_3})\\
& + &
2\dps\sum_{j_3<j_2<j_1}a_{j_1j_2}a_{j_2j_3}
P_1(X_{j_1})P_2(X_{j_2})P_1(X_{j_3})\\
& + &
2\dps\sum_{j_3<j_2<j_1}a_{j_1j_3}a_{j_2j_3}
P_1(X_{j_1})P_1(X_{j_2})P_2(X_{j_3})\\
& + &
\dps\sum_{j_2<j_1}a_{j_1j_2}^2
P_2(X_{j_1})P_2(X_{j_2})\\\hline
\end{array}\right\vert$$
\subsubsection*{\underline{\bf Etape 2 :} calcul de $\left(\dps\sum_{j\geq1}
\langle f , e_{j}\circ e_{j}\rangle (X_{j}^2-1)\right)^2=:N^2$}
$$N^2 = \sum_{j\geq1}
a_j^2\left(X_j^2-1\right)^2 +
2\sum_{k<j}a_ka_j
\left(X_k^2X_j^2-X_j^2-X_k^2+1\right)$$
On regarde maintenant dans la base des polynômes :
$$\begin{array}{lcl}
N^2 & = & \dps\sum_{j\geq1}a_j^2
\Bigl[P_4(X_j)+\gamma_{4,3}P_3(X_j)+\gamma_{4,2}P_2(X_j)+\gamma_{4,1}P_1(X_j)
+\gamma_{4,0}\\
 & &\hspace{1,1 cm} -2P_2(X_j)-2m_3P_1(X_j)-2+1\Bigl]\\
& + & 2\dps\sum_{k<j}\Bigl[
\Bigl(P_2(X_k)+m_3P_1(X_k)+1\Bigl)\Bigl(P_2(X_j)+m_3P_1(X_j)+1\Bigl)
\\
& &\hspace{1,1 cm}-\Bigl(P_2(X_k)+m_3P_1(X_k)+1\Bigl)
-\Bigl(P_2(X_j)+m_3P_2(X_j)+1\Bigl)+1\Bigl]\\
\end{array}$$
On calcule et ordonne :

$$\left\vert\begin{array}{lcl}\hline
N^2 & = &\\
\mbox{(ordre 4)} & & \dps\sum_{j\geq1}a_j^2P_4(X_j) +
 2\dps\sum_{k<j}a_ka_jP_4(X_k)P_4(X_j)\\
\mbox{(ordre 3)} & + & \gamma_{4,3}\dps\sum_{j\geq1}a_j^2P_3(X_j) +
2m_3\dps\sum_{k<j}a_ka_j
\left(P_2(X_k)P_1(X_j)+P_1(X_k)P_2(X_j)\right)\\
\mbox{(ordre 2)} & + & (\gamma_{4,2}-2)\dps\sum_{j\geq1}a_j^2P_2(X_j) +
2m_3^2\dps\sum_{k<j}a_ka_j
P_1(X_k)P_1(X_j)\\
\mbox{(ordre 1)} & + & (\gamma_{4,1}-2m_3)\dps\sum_{j\geq1}a_j^2P_1(X_j)\\
\mbox{(ordre 0)} & + &
(\gamma_{4,0}-1)\dps\sum_{j\geq1}a_j^2\\
\hline
\end{array}\right\vert$$
\subsubsection*{\underline{\bf Etape 3 :} calcul de $\dps\sum_{j_2<j_1}
\langle f , e_{j_1}\circ e_{j_3}\rangle X_{j_1}X_{j_2}
\sum_{j_3\geq1}
\langle f , e_{j_3}\circ e_{j_3}\rangle (X_{j_3}^2-1)$}
$$\begin{array}{lclc}
MN & = & \dps\sum_{j_1\geq1}\sum_{j_2=0}^{j_1-1}
a_{j_1j_2}X_{j_2}X_{j_1}
\sum_{j_3\geq1}a_{j_3}(X_{j_3}^2-1)\\
 &= &\dps\sum_{j_1\geq1}
\Bigl(\sum_{j_2<j_1}a_{j_1j_2}X_{j_2}\Bigl)a_{j_1}(X_{j_1}^2-1) & j_1=j_3\\
& + &\dps\sum_{j_3\geq1}
\sum_{j_1<j_3}\sum_{j_2<j_1}
a_{j_1j_2}a_{j_3}X_{j_1}X_{j_2}(X_{j_3}^2-1) & j_1<j_3\\
& + &\dps\sum_{j_1\geq1}\left(
\sum_{j_2<j_1}
a_{j_1j_2}X_{j_2}
\sum_{j_3<j_1}a_{j_3}(X_{j_3}^2-1)\right)
                     X_{j_1} & j_3<j_1\\
&= & \dps\sum_{j_2<j_1}
a_{j_1}a_{j_1j_2}(X_{j_1}^3-X_{j_1})X_{j_2}&\\
&+& \dps\sum_{j_3<j_2<j_1}
a_{j_2j_3}(X_{j_1}^2-1)X_{j_2}X_{j_3}&\\
&+& \dps\sum_{j_1\geq1}\Bigl[X_{j_1}
\sum_{j_2<j_1}a_{j_1j_2}a_{j_2}X_{j_2}(X_{j_2}^2-1) & j_2=j_3\\
&+&
\dps\sum_{j_2<j_1}\sum_{j_3<j_2}
a_{j_1j_2}a_{j_3}X_{j_2}(X_{j_3}^2-1) & j_3<j_2\\
&+&
\dps\sum_{j_3<j_1}\sum_{j_2<j_3}
a_{j_1j_2}a_{j_3}X_{j_2}(X_{j_3}^2-1) \Bigl]& j_2<j_3\\
 &= &
\dps\sum_{j_3<j_2<j_1}
a_{j_1}a_{j_2j_3}(X_{j_1}^2-1)X_{j_2}X_{j_3}\\
& + &
\dps\sum_{j_3<j_2<j_1}
a_{j_1j_3}a_{j_2}X_{j_1}(X_{j_2}^2-1)X_{j_3}\\
&&\dps\sum_{j_3<j_2<j_1}
a_{j_1j_2}a_{j_3}X_{j_1}X_{j_2}(X_{j_3}^2-1)\\
& + &\dps\sum_{j_2<j_1}
a_{j_1j_2}a_{j_1}(X_{j_1}^3-X_{j_1})X_{j_2}\\
& + & \dps\sum_{j_2<j_1}
a_{j_1j_2}a_{j_2}X_{j_1}(X_{j_2}^3-X_{j_2})\\
\end{array}$$\\
On ordonne et on utilise la base des $(P_n)_n$ :
$$\begin{array}{lcl}
MN & = & \dps\sum_{j_2<j_1}
a_{j_1j_2}a_{j_1}
\Bigl(P_3(X_{j_1})+\gamma_{3,2}P_2(X_{j_1})+\gamma_{3,1}P_1(X_{j_1})
+\gamma{3,0}-P_1(X_{j_1})\Bigl)P_1(X_{j_2})\\
 & + &\dps\sum_{j_2<j_1}
a_{j_1j_2}a_{j_2}P_1(X_{j_1})
\Bigl(P_3(X_{j_2})+\gamma_{3,2}P_2(X_{j_2})+\gamma_{3,1}P_1(X_{j_2})
+\gamma{3,0}-P_1(X_{j_2})\Bigl)\\
& + &\dps\sum_{j_3<j_2<j_1}
a_{j_1}a_{j_2j_3}
\Bigl(P_2(X_{j_1})+m_3P_1(X_{j_1})\Bigl)
P_1(X_{j_2})P_1(X_{j_3})\\
& + &\dps\sum_{j_3<j_2<j_1}
a_{j_1j_3}a_{j_2}P_1(X_{j_1})
\Bigl(P_2(X_{j_2})+m_3P_1(X_{j_2})\Bigl)P_1(X_{j_3})\\
& + &\dps\sum_{j_3<j_2<j_1}
a_{j_1j_2}a_{j_3}P_1(X_{j_1})P_1(X_{j_2})
\Bigl(P_2(X_{j_3})+m_3P_1(X_{j_3})\Bigl)\\
\end{array}$$
On calcule et on réordonne :
$$\left\vert\begin{array}{lcl}\hline
MN & = &\\
\mbox{(ordre 4)}& & \dps\sum_{j_2<j_1}
a_{j_1j_2}a_{j_1}P_3(X_{j_1})P_1(X_{j_2})\\
& + & \dps\sum_{j_2<j_1}
a_{j_1j_2}a_{j_2}P_1(X_{j_1})P_3(X_{j_2})\\
& + & \dps\sum_{j_3<j_2<j_1}
a_{j_1}a_{j_2j_3}P_2(X_{j_1})P_1(X_{j_2})P_1(X_{j_3})\\
& + & \dps\sum_{j_3<j_2<j_1}
a_{j_1j_3}a_{j_2}P_1(X_{j_1})P_2(X_{j_2})P_1(X_{j_3})\\
& + & \dps\sum_{j_3<j_2<j_1}
a_{j_1j_2}a_{j_3}P_1(X_{j_1})P_1(X_{j_2})P_2(X_{j_3})\\
\mbox{(ordre 3)}& + & \gamma_{3,2}\dps\sum_{j_2<j_1}
a_{j_1j_2}a_{j_1}P_2(X_{j_1})P_1(X_{j_2})\\
& + & \gamma_{3,2}\dps\sum_{j_2<j_1}
a_{j_1j_2}a_{j_2}P_1(X_{j_1})P_2(X_{j_2})\\
& + & m_3\dps\sum_{j_3<j_2<j_1}
\bigl(a_{j_1}a_{j_2j_3}+a_{j_1j_3}a_{j_2}+a_{j_1j_2}a_{j_3}\bigl)
P_1(X_{j_1})P_1(X_{j_2})P_1(X_{j_3})\\
\mbox{(ordre 2)}& + & (\gamma_{3,1}-1)\dps\sum_{j_2<j_1}
\bigl(a_{j_1j_2}a_{j_1}+a_{j_1j_2}a_{j_2}\bigl)
P_1(X_{j_1})P_1(X_{j_2})\\
\mbox{(ordre 1)}& + & \gamma_{3,0}\dps\sum_{j_2<j_1}
\biggl(a_{j_1j_2}a_{j_1}P_1(X_{j_2})+
a_{j_1j_2}a_{j_2}P_1(X_{j_1})\biggl)\\ \hline
\end{array}\right\vert$$
\subsubsection*{\underline{\bf Etape 4 :} sommation}
On calcule simplement $4M^2+N^2+4MN$ et on a le résultat annoncé.
\end{enumerate} \normalsize$\triangle$
\section{identification}
L'expression trouvée n'est, pour le moins, pas très maniable. Nous
allons identifier chacun des termes, ordre par ordre, de façon à
pouvoir, par la suite, calculer $EM^4$ pour montrer l'existence de
modification continues et expliciter la variation
quadratique. Nous noterons
$O_k\biggl(\bigl[\bigl(\Phi^{\circ^2}+\Phi\circ a_1^2\bigl)(f)\bigl]^2\biggl)$
les termes d'ordre k de
$\bigl[\bigl(\Phi^{\circ^2}+\Phi\circ a_1^2\bigl)(f)\bigl]^2$
\subsection{Calculs préliminaires}
Dans cette section nous allons utiliser le développement en série
de $f\circ f$ et un opérateur introduit par P. A. Meyer
{\bf{[10]}}.
\begin{lm}
Soit $f\in H\circ H$ on a :
$$\left\vert\begin{array}{lcl}\hline
f\circ f & = & \dps\sum_{j_4<j_3<j_2<j_1}
8(a_{j_1j_2}a_{j_3j_4}+a_{j_1j_3}a_{j_2j_4}+a_{j_1j_4}a_{j_2j_3})
e_{j_1}\circ e_{j_2}\circ e_{j_3}\circ e_{j_4}\\
& + &
 \dps\sum_{j_3<j_2<j_1}
(8a_{j_1j_2}a_{j_1j_3}+4a_{j_1}a_{j_2j_3})
e_{j_1}^{\circ^2}\circ e_{j_2}\circ e_{j_3}\\
& + &
 \dps\sum_{j_3<j_2<j_1}
(8a_{j_1j_2}a_{j_2j_3}+4a_{j_2}a_{j_1j_3})
e_{j_1}\circ e_{j_2}^{\circ^2}\circ e_{j_3}\\
& + &
 \dps\sum_{j_3<j_2<j_1}
(8a_{j_1j_3}a_{j_2j_3}+4a_{j_1j_2}a_{j_3})
e_{j_1}\circ e_{j_2}\circ e_{j_3}^{\circ^2}\\
& + &
 \dps\sum_{j_2<j_1}
4a_{j_1}a_{j_1j_2}
e_{j_1}^{\circ^3}\circ e_{j_2}\\
& + &
 \dps\sum_{j_2<j_1}
4a_{j_1j_2}a_{j_2}
e_{j_1}\circ e_{j_2}^{\circ^3}\\
& + &
 \dps\sum_{j_2<j_1}
(4a_{j_1j_2}^2+2a_{j_1}a_{j_2})
e_{j_1}^{\circ^2}\circ e_{j_2}^{\circ^2}\\
& + &
 \dps\sum_{1\leq j}
a_{j}^2
e_{j}^{\circ^4}\\ \hline
\end{array}\right\vert$$
\end{lm}
\subsection*{\underline{\textbf{Démonstration}}}
\begin{enumerate}
\item[]
\small
$f=\dps\sum_{(j_1, j_2)}\langle f ,e_{j_1}\circ e_{j_2}\rangle
e_{j_1}\circ e_{j_2}$\\
$f\circ f = \dps\sum_{(j_1, j_2, j_3, j_4)}
\langle f ,e_{j_1}\circ e_{j_2}\rangle\langle f ,e_{j_3}\circ e_{j_4}\rangle
e_{j_1}\circ e_{j_2}e_{j_3}\circ e_{j_4}$\\
On va maintenant réécrire cette somme en fonction du nombre
d'indices égaux, cela revient à regarder les différentes écritures du nombre
quatre en somme d'entiers.\\
\fbox{4 indices égaux}\;(4 = 4)\\
\begin{center}
\fbox{$\dps\sum_{j\geq 1}\langle f ,e_{j}\circ e_{j}\rangle
e_{j}^{\circ^4}$}
\end{center}
\vspace*{1cm}
\fbox{3 indices égaux}\;(4 = 3+1)\\
Dans chacune des sommes on appellera $j_1$ l'indice distinct :\\
\underline{$j_1\not=j_2=j_3=j_4$}
$$\sum_{j_1\geq 1}\sum_{j_2\geq 1}
\langle f ,e_{j_2}\circ e_{j_1}\rangle
\langle f ,e_{j_1}\circ e_{j_1}\rangle
e_{j_1}^{\circ^3}\circ e_{j_2}$$
\underline{$j_2\not=j_1=j_3=j_4$}
$$\sum_{j_1\geq 1}\sum_{j_2\geq 1}
\langle f ,e_{j_1}\circ e_{j_2}\rangle
\langle f ,e_{j_1}\circ e_{j_1}\rangle
e_{j_1}^{\circ^3}\circ e_{j_2}$$
\underline{$j_3\not=j_1=j_2=j_4$}
$$\sum_{j_1\geq 1}\sum_{j_2\geq 1}
\langle f ,e_{j_1}\circ e_{j_1}\rangle
\langle f ,e_{j_2}\circ e_{j_1}\rangle
e_{j_1}^{\circ^3}\circ e_{j_2}$$
\underline{$j_4\not=j_1=j_2=j_3$}
$$\sum_{j_1\geq 1}\sum_{j_2\geq 1}
\langle f ,e_{j_1}\circ e_{j_1}\rangle
\langle f ,e_{j_1}\circ e_{j_2}\rangle
e_{j_1}^{\circ^3}\circ e_{j_2}$$
On fait la somme de ces trois termes et on trouve :
$$4\dps\sum_{j_1\geq 1}\sum_{j_2\geq 1}
\langle f ,e_{j_1}\circ e_{j_1}\rangle
\langle f ,e_{j_1}\circ e_{j_2}\rangle
 e_{j_1}^{\circ^3}\circ e_{j_2}$$
$$\begin{array}{clc}
=& 4\dps\sum_{j_1\geq1}\sum_{j_2=1}^{j_1-1}
\langle f ,e_{j_1}\circ e_{j_2}\rangle
\langle f ,e_{j_1}\circ e_{j_1}\rangle
e_{j_1}^{\circ^3}\circ e_{j_2} & j_1<j_2\\
+& 4\dps\sum_{j_2\geq1}\sum_{j_1=1}^{j_2-1}
\langle f ,e_{j_1}\circ e_{j_2}\rangle
\langle f ,e_{j_1}\circ e_{j_1}\rangle
e_{j_1}^{\circ^3}\circ e_{j_2} & j_2<j_1\\
+& 4\dps\sum_{j_1\geq1}\sum_{j_2=1}^{j_1-1}\bigl(
a_{j_1j_2}a_{j_2}e_{j_1}^{\circ^3}\circ e_{j_2}+
a_{j_1}a_{j_1j_2}e_{j_1}\circ e_{j_2}^{\circ^3}\bigl)\\
\end{array}$$\\
\fbox{2 indices égaux les autre distincts} \;(4=2+1+1)\\
$$\begin{array}{ll}
j_1=j_2\\
+ & \dps\sum_{j_1\not=j_2\not=j_3}
\langle f ,e_{j_1}\circ e_{j_1}\rangle
\langle f ,e_{j_2}\circ e_{j_3}\rangle
e_{j_1}^{\circ^2}\circ e_{j_2} \circ e_{j_3}\\
j_1=j_3\\
+ & \dps\sum_{j_1\not=j_2\not=j_3}
\langle f ,e_{j_1}\circ e_{j_2}\rangle
\langle f ,e_{j_1}\circ e_{j_3}\rangle
e_{j_1}^{\circ^2}\circ e_{j_2} \circ e_{j_3}\\
j_1=j_4\\
+ & \dps\sum_{j_1\not=j_2\not=j_3}
\langle f ,e_{j_1}\circ e_{j_3}\rangle
\langle f ,e_{j_2}\circ e_{j_1}\rangle
e_{j_1}^{\circ^2}\circ e_{j_2} \circ e_{j_3}\\
j_2=j_3\\
+ & \dps\sum_{j_1\not=j_2\not=j_3}
\langle f ,e_{j_2}\circ e_{j_1}\rangle
\langle f ,e_{j_1}\circ e_{j_3}\rangle
e_{j_1}^{\circ^2}\circ e_{j_2} \circ e_{j_3}\\
j_2=j_4\\
+ & \dps\sum_{j_1\not=j_2\not=j_3}
\langle f ,e_{j_2}\circ e_{j_1}\rangle
\langle f ,e_{j_3}\circ e_{j_1}\rangle
e_{j_1}^{\circ^2}\circ e_{j_2} \circ e_{j_3}\\
j_3=j_4\\
+ & \dps\sum_{j_1\not=j_2\not=j_3}
\langle f ,e_{j_2}\circ e_{j_3}\rangle
\langle f ,e_{j_1}\circ e_{j_1}\rangle
e_{j_1}^{\circ^2}\circ e_{j_2} \circ e_{j_3}\\
= & \dps2\sum_{j_1\not=j_2\not=j_3}
\langle f ,e_{j_1}\circ e_{j_1}\rangle
\langle f ,e_{j_2}\circ e_{j_3}\rangle
e_{j_1}^{\circ^2}\circ e_{j_2} \circ e_{j_3}\\
&+  \dps4\sum_{j_1\not=j_2\not=j_3}
\langle f ,e_{j_1}\circ e_{j_2}\rangle
\langle f ,e_{j_1}\circ e_{j_3}\rangle
e_{j_1}^{\circ^2}\circ e_{j_2} \circ e_{j_3}\\
\end{array}$$
Pour expliciter ces deux sommes on va considérer les différentes
permutations suivantes, qui constitue une partition de l'ensemble
$\{(j_1,j_2,j_3)\;\vert\;j_1\not=j_2\not=j_3\}$ :
$$\begin{array}{ccc}
j_1<j_2<j_3 & j_2<j_1<j_3 & j_2<j_3<j_1\\
j_1<j_3<j_2 & j_3<j_1<j_2 & j_3<j_2<j_1
\end{array}$$
\underline{Première somme}\\
après réindexaton on trouve :
$$\begin{array}{ll}
2\dps\sum_{j_3<j_2<j_1}\Bigl[&
\langle f ,e_{j_3}\circ e_{j_3}\rangle
\langle f ,e_{j_2}\circ e_{j_1}\rangle
e_{j_1}\circ e_{j_2} \circ e_{j_3}^{\circ^2}\\
& +
\langle f ,e_{j_2}\circ e_{j_2}\rangle
\langle f ,e_{j_1}\circ e_{j_3}\rangle
e_{j_1}\circ e_{j_2}^{\circ^2} \circ e_{j_3}\\
& +
\langle f ,e_{j_1}\circ e_{j_1}\rangle
\langle f ,e_{j_3}\circ e_{j_3}\rangle
e_{j_1}^{\circ^2}\circ e_{j_2} \circ e_{j_3}\Bigl]\\
\end{array}$$
\underline{deuxième somme}\\
après réindexaton on trouve :
$$\begin{array}{cl}
\dps\sum_{j_3<j_2<j_1}\Bigl[&
a_{j_1j_2}a_{j_1j_3}e_{j_1}^{\circ^2}\circ e_{j_2} \circ e_{j_3}+
a_{j_2j_1}a_{j_2j_3}e_{j_1}\circ e_{j_2} ^{\circ^2}\circ e_{j_3}\\
& +
a_{j_3j_1}a_{j_3j_2}e_{j_1}\circ e_{j_2} \circ e_{j_3}^{\circ^2}+
a_{j_1j_3}a_{j_1j_2}e_{j_1} ^{\circ^2}\circ e_{j_2}\circ e_{j_3}\\
& +
a_{j_2j_3}a_{j_2j_1}e_{j_1}\circ e_{j_2}^{\circ^2} \circ e_{j_3}+
a_{j_3j_2}a_{j_3j_1}e_{j_1} \circ e_{j_2}\circ e_{j_3}^{\circ^2}\Bigl]\\
= &\dps\sum_{j_3<j_2<j_1}
(a_{j_1j_2}a_{j_1j_3}+a_{j_1j_3}a_{j_1j_2})
e_{j_1}^{\circ^2}\circ e_{j_2} \circ e_{j_3}\\
& +\dps\sum_{j_3<j_2<j_1}
2a_{j_1j_2}a_{j_2j_3}
e_{j_1}\circ e_{j_2}^{\circ^2}\circ e_{j_3}\\
& +\dps\sum_{j_3<j_2<j_1}
2a_{j_1j_3}a_{j_2j_3}
e_{j_1}\circ e_{j_2}\circ e_{j_3}^{\circ^2}\\
\end{array}$$
On réunit ces deux résultats et on trouve :
$$\left\vert\begin{array}{cl}\hline
& \dps\sum_{j_3<j_2<j_1}\bigl(
4a_{j_1}a_{j_2j_3} + 8a_{j_1j_2}a_{j_1j_3}\bigl)
e_{j_1}^{\circ^2}\circ e_{j_2}\circ e_{j_3}\\
 + & \dps\sum_{j_3<j_2<j_1}\bigl(
4a_{j_2}a_{j_1j_3} + 8a_{j_1j_2}a_{j_2j_3}\bigl)
e_{j_1}\circ e_{j_2}^{\circ^2}\circ e_{j_3}\\
+ & \dps\sum_{j_3<j_2<j_1}\bigl(
4a_{j_1j_2}a_{j_3} + 8a_{j_1j_3}a_{j_2j_3}\bigl)
e_{j_1}\circ e_{j_2}\circ e_{j_3}^{\circ^2}\\ \hline
\end{array}\right\vert$$\\
\fbox{2 termes égaux $\not=$ deux termes égaus}\; (4=2+2)\\
On a quatre cas à considérer :
$$j_1=j_2\not=j_3=j_4$$
$$j_1=j_3\not=j_2=j_4$$
$$j_1=j_4\not=j_2=j_4$$
On a alors la somme :\\
$\dps\sum_{j_1\not=j_2}\Bigl(a_{j_1}a_{j_2}+2a_{j_1j_2}^2\Bigl)
e_{j_1}^{\circ^2}\circ e_{j_2}^{\circ^2}$
$$\left\vert\begin{array}{cc}\hline
= & \dps\sum_{j_1\geq1}\sum_{j_2<j_1}\Bigl(
a_{j_1}a_{j_2}+2a_{j_1j_2}^2\Bigl)
e_{j_1}^{\circ^2}\circ e_{j_2}^{\circ^2}\\ \hline
\end{array}\right\vert$$\\
\fbox{Tous les indices distincts}\; (4=1+1+1+1)\\
On doit considérer l'ensemble
$\{(j_1,j_2,j_3,j_4)\; \vert\;j_1\not=j_2\not=j_3\not=j_4\}$,
puisque nous étudions la somme :
$$\sum_{j_1\not=j_2\not=j_3\not=j_4}
a_{j_1j_2}a_{j_3j_4}e_{j_1}\circ e_{j_2}\circ e_{j_3}\circ e_{j_4}$$
Les 24 permutations ne sont pas à regarder puisque échanger $j_1$
et $j_2$ ou $j_3$ et $j_4$ n'a aucun effet. Il reste en fait 3
termes qui se répètent 8 fois :
$$\left\vert\begin{array}{c}\hline
\dps\sum_{j_4<j_3<j_2<j_1}8\Bigl(
a_{j_1j_2}a_{j_3j_4}+a_{j_1j_3}a_{j_2j_4}+a_{j_1j_2}a_{j_2j_3}\Bigl)
e_{j_1}\circ e_{j_2}\circ e_{j_3}\circ e_{j_4}\\ \hline
\end{array}\right\vert$$
Pour terminer il ne reste qu'à sommer tous ces résultats
\end{enumerate} \normalsize$\triangle$
\begin{de}
$$\bigl(f\mathop{\sim}\limits_{1}f\bigl)(s_1,s_3):=\int f(s_1,s_2) f(s_3,s_2)d\mu(s_2)$$
\end{de}
\begin{lm} Pour $f\in H\circ H$ :
$$\left\vert\begin{array}{lcl}\hline
\bigl(f\mathop{\sim}\limits_{1}f\bigl)= & &
\dps\sum_{j_3<j_2<j_1}
2\bigl(
a_{j_1j_2}a_{j_2j_3}e_{j_1}\circ e_{j_3}+
a_{j_1j_3}a_{j_2j_3}e_{j_1}\circ e_{j_2}+
a_{j_1j_2}a_{j_1j_3}e_{j_2}\circ e_{j_3}\bigl)\\
& + &
\dps\sum_{j_2<j_1}
2\bigl(
a_{j_1j_2}a_{j_2}+
a_{j_1j_2}a_{j_1}\bigl)e_{j_1}\circ e_{j_2}\\
& + &
\dps\sum_{j_2<j_1}
a_{j_1j_2}^2\bigl(e_{j_1}^{\circ^2}+ e_{j_2}^{\circ^2}\bigl)\\
& + &
\dps\sum_{j\geq1}
a_{j}^2e_{j}^{\circ^2}\\ \hline
\end{array}\right\vert$$
\end{lm}
\subsection*{\underline{\textbf{Démonstration}}}
\begin{enumerate}
\item[]
\small
$f\mathop{\sim}\limits_{1}f=\dps
\dps\sum_{(j_1,j_2,j_3,j_4)}a_{j_1j_2}a_{j_3j_4}
\Bigl(
\langle e_{j_1},e_{j_3}\rangle e_{j_2}\otimes e_{j_4}+
\langle e_{j_1},e_{j_4}\rangle e_{j_2}\otimes e_{j_3}+
\langle e_{j_2},e_{j_3}\rangle e_{j_1}\otimes e_{j_4}+
\langle e_{j_2},e_{j_4}\rangle e_{j_1}\otimes e_{j_3}\Bigl)$\\
$\begin{array}{clc}
&\dps\int f(s_1,s_2)f(s_2,s_3)d\mu(s_2)\\
= & \dps\sum_{(j_1,j_2)}
\langle f,e_{j_1}\circ e_{j_2}\rangle
\sum_{(j_3,j_4)}
\langle f,e_{j_3}\circ e_{j_4}\rangle
\int e_{j_1}\circ e_{j_2}(s_1,s_2)e_{j_3}\circ
e_{j_4}(s_2,s_3)d\mu(s_2)\\
= &\dps \frac{1}{4}\sum_{(j_1,j_2,j_3,j_4)}a_{j_1j_2}a_{j_3j_4}
\int
\Bigl(e_{j_1}\otimes e_{j_2}(s_1,s_2)+e_{j_2}\otimes e_{j_1}(s_1,s_2)
\Bigl)\\
 & \hspace{2cm}\cdot\Bigl(e_{j_3}\otimes e_{j_4}(s_2,s_3)+e_{j_4}\otimes e_{j_3}(s_2,s_3)
\Bigl)d\mu(s_2)\\
= &\dps \frac{1}{4}\sum_{(j_1,j_2,j_3,j_4)}a_{j_1j_2}a_{j_3j_4}
\Bigl(
e_{j_1}\otimes e_{j_4}\delta_{j_2j_3}+e_{j_1}\otimes e_{j_3}\delta_{j_2j_4}
+e_{j_2}\otimes e_{j_4}\delta_{j_1j_3}\\
& \hspace*{4cm}+e_{j_2}\otimes e_{j_3}\delta_{j_1j_4}
\Bigl)(s_1,s_3)\\
= &\dps \frac{1}{4}
\sum_{(j_1,j_2,j_3)}a_{j_1j_2}a_{j_2j_3}e_{j_1}\otimes e_{j_3}(s_1,s_3)+
\frac{1}{4}\sum_{(j_1,j_2,j_3)}a_{j_1j_2}a_{j_3j_2}e_{j_1}\otimes e_{j_3}(s_1,s_3)\\
&+\dps \frac{1}{4}
\sum_{(j_1,j_2,j_3)}a_{j_1j_2}a_{j_1j_3}e_{j_2}\otimes e_{j_3}(s_1,s_3)+
\frac{1}{4}\sum_{(j_1,j_2,j_3)}a_{j_3j_1}a_{j_3j_2}e_{j_2}\otimes e_{j_3}(s_1,s_3)\\
= &\dps \frac{1}{2}
\sum_{(j_1,j_2,j_3)}a_{j_1j_2}a_{j_2j_3}e_{j_1}\otimes e_{j_3}(s_1,s_3)+
\frac{1}{2}
\sum_{(j_1,j_2,j_3)}a_{j_1j_2}a_{j_1j_3}e_{j_2}\otimes e_{j_3}(s_1,s_3)\\
&   \mbox{(en échangeant $j_1$ et $j_2$ dans la deuxième somme)}\\
= &\dps\sum_{(j_1,j_2,j_3)}a_{j_1j_2}a_{j_2j_3}e_{j_1}\otimes e_{j_3}(s_1,s_3)\\
\fbox{1+1+1}&\dps\sum_{j_1\not=j_2\not=j_3}
a_{j_1j_2}a_{j_2j_3}e_{j_1}\otimes e_{j_3}(s_1,s_3)\\
& \begin{array}{ccccc}
j_1<j_2<j_3  & ; & j_2<j_1<j_3  & ; &  j_2<j_3<j_1\\
j_1<j_3<j_2  & ; & j_3<j_1<j_2  & ; &  j_3<j_2<j_1\\
\end{array}\\
& \mbox{(après réindexation)}\\
= & \dps\sum_{j_3<j_2<j_1}
\Bigl[
a_{j_3j_2}a_{j_2j_1}e_{j_3}\otimes e_{j_1}(s_1,s_3)+
a_{j_2j_3}a_{j_3j_1}e_{j_2}\otimes e_{j_1}(s_1,s_3)
\\
& +
a_{j_1j_3}a_{j_3j_2}e_{j_1}\otimes e_{j_2}(s_1,s_3)+
a_{j_3j_1}a_{j_1j_2}e_{j_3}\otimes e_{j_2}(s_1,s_3)\\
& +
a_{j_2j_1}a_{j_1j_3}e_{j_2}\otimes e_{j_3}(s_1,s_3)+
a_{j_1j_2}a_{j_2j_3}e_{j_1}\otimes e_{j_3}(s_1,s_3)
\Bigl]\\
= &2\dps\sum_{j_3<j_2<j_1}
\Bigl(
a_{j_1j_2}a_{j_2j_3}e_{j_1}\circ e_{j_3}(s_1,s_3)+
a_{j_1j_3}a_{j_2j_3}e_{j_1}\circ e_{j_2}(s_1,s_3)+\\
& \hspace{3cm}+a_{j_1j_2}a_{j_1j_3}e_{j_2}\circ e_{j_3}(s_1,s_3)
\Bigl)\\
\fbox{2+1} & \dps\sum_{j_2\not=j_1\atop j_2=j_3}
a_{j_1j_2}a_{j_2}e_{j_1}\otimes e_{j_2}(s_1,s_3)+
\sum_{j_2\not=j_1\atop j_1=j_2}
a_{j_1}a_{j_2j_3}e_{j_1}\otimes e_{j_2}(s_1,s_3)+\\
 & \hspace{4cm}+\dps\sum_{j_2\not=j_1\atop j_1=j_3}
a_{j_1j_2}a_{j_2j_3}e_{j_1}\otimes e_{j_1}(s_1,s_3)\\
= &2\dps\sum_{j_2<j_1}
a_{j_1j_2}a_{j_2}e_{j_1}\circ e_{j_2}(s_1,s_3)+
2\dps\sum_{j_2<j_1}
a_{j_1}a_{j_1j_2}e_{j_1}\circ e_{j_2}(s_1,s_3)+\\
 & \hspace{4cm}+2\dps\sum_{j_2<j_1}
a_{j_1j_2}^2\Bigl(e_{j_1}\circ e_{j_1}+e_{j_2}\circ
e_{j_2}\Bigl)(s_1,s_3)\\
\fbox{3}& \dps\sum_ja_j^2e_j\circ e_j
\end{array}$\\
Pour conclure il suffit de réunir tous ces résultats.
\end{enumerate} \normalsize$\triangle$
\subsection{Identification des termes d'ordre 4}
\begin{pp}
$$O_4\Bigl(\bigl[(\Phi^{\circ^2}+\Phi\circ a_1^2)(f)\bigl]^2\Bigl)=
\Phi^{\circ^4}(f\circ f)$$
\end{pp}
\subsection*{\underline{\textbf{Démonstration}}}
\begin{enumerate}
\item[]
\small
C'est une simple comparaison des termes en question.
\end{enumerate} \normalsize$\triangle$
\subsection{Identification des termes d'ordre 3}
\begin{lm}
$$\left\vert\begin{array}{ccl}\hline
a_1^4(f\circ f)
=&&\dps\sum_{j_3<j_2<j_1}\gamma_{2,1}\Bigl[
8\bigl(
a_{j_1j_2}a_{j_1j_3}+a_{j_1j_2}a_{j_2j_3}+a_{j_1j_3}a_{j_2j_3}
\bigl)\\
&&\hspace{1,5cm}+4\bigl(a_{j_1}a_{j_2j_3}+a_{j_2}a_{j_1j_3}+a_{j_3}a_{j_1j_2}
\bigl)\Bigl]e_{j_1}\circ e_{j_2}\circ e_{j_3}\\
&+&\dps\sum_{j_2<j_1}
\Bigl[4\gamma_{3,2}a_{j_1}a_{j_1j_2}+
\gamma_{2,1}(4a_{j_1j_2}^2+2a_{j_1}a_{j_2})
\Bigl]e_{j_1}^{\circ^2}\circ e_{j_2}\\
&+&\dps\sum_{j_2<j_1}
\Bigl[4\gamma_{3,2}a_{j_1j_2}a_{j_2}+
\gamma_{2,1}(4a_{j_1j_2}^2+2a_{j_1}a_{j_2})
\Bigl]e_{j_1}\circ e_{j_2}^{\circ^2}\\
&+&\dps\sum_{1\leq j}
\gamma_{4,3}a_{j}^2 e_{j}^{\circ^3}\\ \hline
\end{array}\right\vert$$
\end{lm}
\subsection*{\underline{\textbf{Démonstration}}}
\begin{enumerate}
\item[]
\small
Fixons le 4-uple $(j_1,j_2,j_3,j_4)$ tel que $j_4<j_3<j_2<j_1$ :\\
$$\begin{array}{cclcl}
a_1^4 & : & H^{\circ^4} & \to & H^{\circ^3}\\
&&e_{j_1}\circ e_{j_2}\circ e_{j_3}\circ e_{j_4}&
\mapsto & 0\\
&& e_{j_1}^{\circ^2}\circ e_{j_2}\circ e_{j_3} &
\mapsto & \gamma_{2,1} e_{j_1}\circ e_{j_2}\circ e_{j_3}\\
&& e_{j_1}^{\circ^3}\circ e_{j_2} &
\mapsto & \gamma_{3,2} e_{j_1}^{\circ^2}\circ e_{j_2}\\
&& e_{j_1}^{\circ^2}\circ e_{j_2}^{\circ^2} &
\mapsto & \gamma_{2,1} e_{j_1}^{\circ^2}\circ e_{j_2}+
\gamma_{2,1} e_{j_1}\circ e_{j_2}^{\circ^2}\\
&& e_{j_1}^{\circ^4} &
\mapsto & \gamma_{4,3} e_{j_1}^{\circ^3}\\
\end{array}$$
\end{enumerate} \normalsize$\triangle$
\begin{pp}
$$O_3\Bigl(\bigl[(\Phi^{\circ^2}+\Phi\circ a_1^2)(f)\bigl]^2\Bigl)=
\Phi^{\circ^3}\bigl(a_1^4(f\circ f)\Bigl)$$
\end{pp}
\subsection*{\underline{\textbf{Démonstration}}}
\begin{enumerate}
\item[]
\small
C'est une simple comparaison des termes en question.
\end{enumerate} \normalsize$\triangle$
\subsection{Identification des termes d'ordre 2}
\begin{pp}
$$O_2\Bigl(\bigl[(\Phi^{\circ^2}+\Phi\circ a_1^2)(f)\bigl]^2\Bigl)=
\Phi^{\circ^2}\Bigl(4f\mathop{\sim}\limits_{1}f +a_2^4(f\circ f)
\Bigl)$$
\end{pp}
\subsection*{\underline{\textbf{Démonstration}}}
\begin{enumerate}
\item[]
\small
$$\begin{array}{lcl}
&&O_2\Bigl(\bigl[(\Phi^{\circ^2}+\Phi\circ a_1^2)(f)\bigl]^2\Bigl)
-\Phi^{\circ^2}\bigl(4f\mathop{\sim}\limits_{1}f\bigl)\\
= & &\dps\sum_{j_2<j_1}\Bigl[4m_3^2a_{j_1j_2}^2+2m_3^2a_{j_1}a_{j_2}
+4(\gamma_{3,1}-3)(a_{j_1j_2}a_{j_1}+a_{j_1j_2}a_{j_2})
\Bigl]P_1(X_{j_1}) P_1(X_{j_2})\\
&&+\dps\sum_{1\leq j}(\gamma_{4,2}-6)a_jP_2(X_j)
\end{array}$$
Il nous reste à calculer $\Phi^{\circ^2}\bigl(a_2^4(f\circ
f)\bigl)$.\\
Fixons le 4-uple $(j_1,j_2,j_3,j_4)$ tel que $j_4<j_3<j_2<j_1$ :\\
$$\begin{array}{cclcl}
a_2^4 & : & H^{\circ^4} & \to & H^{\circ^3}\\
&&e_{j_1}\circ e_{j_2}\circ e_{j_3}\circ e_{j_4}&
\mapsto & 0\\
&& e_{j_1}^{\circ^2}\circ e_{j_2}\circ e_{j_3} &
\mapsto & 0\\
&& e_{j_1}^{\circ^3}\circ e_{j_2} &
\mapsto & (\gamma_{3,1} -\Gamma_{3,1})
e_{j_1}\circ e_{j_2}\\
&& e_{j_1}^{\circ^2}\circ e_{j_2}^{\circ^2} &
\mapsto & (\gamma_{2,0}- \Gamma_{2,0})^2
 e_{j_1}\circ e_{j_2}\\
&& e_{j_1}^{\circ^4} &
\mapsto & (\gamma_{4,2}- \Gamma_{4,2})e_{j_1}^{\circ^2}\\
\end{array}$$
\end{enumerate} \normalsize$\triangle$
\subsection{Identification des termes d'ordre 1}
\begin{lm}
$$a_3^4(f\circ f)=
4\gamma_{3,0}\Bigl(\sum_{j_2<j_1}a_{j_1}a_{j_1j_2}e_{j_2}+
\sum_{j_2<j_1}a_{j_1j_2}a_{j_2}e_{j_1}\Bigl)+
\gamma_{4,1}\sum_{j\geq1}a_j^2e_j$$
\end{lm}
\subsection*{\underline{\textbf{Démonstration}}}
\begin{enumerate}
\item[]
\small
Fixons le 4-uple $(j_1,j_2,j_3,j_4)$ tel que $j_4<j_3<j_2<j_1$ :\\
$$\begin{array}{cclcl}
a_3^4 & : & H^{\circ^4} & \to & H^{\circ^3}\\
&& e_{j_1}^{\circ^3}\circ e_{j_2} &
\mapsto & (\gamma_{3,1} -\Gamma_{3,1})
e_{j_1}\circ e_{j_2}\\
&& e_{j_1}^{\circ^2}\circ e_{j_2}^{\circ^2} &
\mapsto & 0\\
&& e_{j_1}^{\circ^4} &
\mapsto & (\gamma_{4,1}- \Gamma_{4,1})e_{j_1}\\
\end{array}$$
\end{enumerate} \normalsize$\triangle$
\begin{lm}
$$a_1^2(f\mathop{\sim}\limits_{1}f) =
m_3\dps\sum_{j_2<j_1}a_{j_1j_2}^2(e_{j_1}+e_{j_2}) +
m_3\dps\sum_{1\leq j}a_j^2e_j$$
\end{lm}
\subsection*{\underline{\textbf{Démonstration}}}
\begin{enumerate}
\item[]
\small
Fixons le couple $(j_1,j_2)$ tel que $j_2<j_2$ :\\
$$a_1^2(e_{j_1}\circ e_{j_2}) =0$$
$$a_1^2(e_{j_1}^{\circ^2}) =(\gamma_{2,1}-\Gamma_{2,1}) = m_3$$
\end{enumerate} \normalsize$\triangle$
\begin{pp}
$$O_1\Bigl(\bigl[(\Phi^{\circ^2}+\Phi\circ a_1^2)(f)\bigl]^2\Bigl)=
\Phi\bigl(
a_3^4(f\circ f) + 4 a_1^2(f\mathop{\sim}\limits_{1}f)-
6\bigl(a_1^2\circ \pi_1\bigl)(f\mathop{\sim}\limits_{1}f)
\Bigl)$$
où $\pi_1$ est la projection orthogonale sur $\{e_j\circ e_k\;\vert\;j=k\}$
\end{pp}
\subsection*{Remarque}
Nous noterons $\pi_2$ la projection sur  $\{e_j\circ e_k\;\vert\;j\not=k\}$
\subsection*{\underline{\textbf{Démonstration}}}
\begin{enumerate}
\item[]
\small
Pour cette identification il suffit de remarquer que : \\
$O_1\Bigl(\bigl[(\Phi^{\circ^2}+\Phi\circ a_1^2)(f)\bigl]^2\Bigl)-
\Phi\bigl(a_3^4(f\circ f) + 4 a_1^2(f\mathop{\sim}\limits_{1}f)\bigl)=
-6m_3\dps\sum_{1\leq j}a_j^2P_1(X_j)$\\
et : $\bigl(a_1^2\circ \pi_1\bigl)(f\mathop{\sim}\limits_{1}f)
=m_3\dps\sum_{1\leq j}a_j^2e_j$

\end{enumerate} \normalsize$\triangle$
\subsection{Identification des termes d'ordre 0}
\begin{pp}
$$O_0\Bigl(\bigl[(\Phi^{\circ^2}+\Phi\circ a_1^2)(f)\bigl]^2\Bigl)=
2\n f\n^2+a_4^4(f\circ f)$$
\end{pp}
\subsection*{\underline{\textbf{Démonstration}}}
\begin{enumerate}
\item[]
\small
On vérifie que $a_4^4(f\circ f)=
(\gamma_{4,0}-\Gamma_{4,0})\dps\sum_{1\leq j}a_j^2$ et
$\n f\n^2=2\dps\sum_{j_2<j_1}a_{j_1j_2}^2+\sum_{1\leq j}a_j^2$
$$\begin{array}{lcl}
O_0\Bigl(\bigl[(\Phi^{\circ^2}+\Phi\circ
a_1^2)(f)\bigl]^2\Bigl)&=&
4\dps\sum_{j_2<j_1}a_{j_1j_2}^2+ (\gamma_{4,0}-1)\sum_{1\leq j}a_j^2\\
& =&2\n f\n^2+(\gamma_{4,0}-3)\dps\sum_{1\leq j}a_j^2\\
& =&2\n f\n^2+a_4^4(f\circ f)\\
\end{array}$$
\end{enumerate} \normalsize$\triangle$
\subsection{Identification}
En réunissant les résultats des précédentes sections on a la
proposition suivante :
\begin{pp}
$$\left\vert\begin{array}{cll}\hline
&\bigl[(\Phi^{\circ^2}+\Phi\circ a_1^2)(f)\bigl]^2\\
 = & \Phi^{\circ^4}(f\circ f)&\mbox{(ordre 4)}\\
 +& \Phi^{\circ^3}\bigl(a_1^4(f\circ f)\Bigl) &\mbox{(ordre 3)}\\
 +&
\Phi^{\circ^2}\Bigl(4f\mathop{\sim}\limits_{1}f +a_2^4(f\circ f)
\Bigl) & \mbox{(ordre 2)}\\
 +& \Phi\bigl(
a_3^4(f\circ f) + 4 a_1^2(f\mathop{\sim}\limits_{1}f)-
6\bigl(a_1^2\circ \pi_1\bigl)(f\mathop{\sim}\limits_{1}f)
\Bigl) & \mbox{(ordre 1)}\\
+ & 2\n f\n^2+a_4^4(f\circ f) & \mbox{(ordre 0)}\\ \hline
\end{array}\right\vert$$
\end{pp}
\section{Etude d'une limite}
Soient $h_1,\;h_2\in H$ et $0=t_0<t_1<\cdots<t_n=t$ une suite de
partitions de $[0,t]$ dont le pas tend vers 0.\\
Nos noterons $f_k:=h_1\otimes \left(h_21_{]0,t_k]}\right)1_C$
\begin{lm}
Pour $\Psi=\Phi^{\circ^2}$ ou $\Psi=\Phi\circ a_1^4$ on a :
$$\mathop{\lim}\limits_{k\to\infty}E\left[
4\Psi
\bigl(\sum_{k}f_k\mathop{\sim}\limits_{1}f_k\bigl)-
\int_0^th_2^2\Psi
\Bigl(\bigl[h_11_{]0,.]}
\bigl]^{\circ^2}\Bigl)d\mu\right]^2=0$$
\end{lm}
\subsection*{\underline{\textbf{Démonstration}}}
\begin{enumerate}
\item[]
\small
Nous commençons par le cas $\Psi=\Phi^{\circ^2}$. Nous
Introduisons quelques notations :\\
\begin{itemize}
\item
$S_k:=\left[
4\Phi^{\circ^2}
\bigl(\dps\sum_{k}f_k\mathop{\sim}\limits_{1}f_k\bigl)-
\int_0^th_2^2\Phi^{\circ^2}
\Bigl(\bigl[h_11_{]0,.]}
\bigl]^{\circ^2}\Bigl)d\mu\right]^2$
\item
$h^k:=h1_{]0,t_k]}$
\end{itemize}
$$\begin{array}{clc}
&4f_k\mathop{\sim}\limits_{1}f_k(s_1,s_3)=\\
&\int h_1^{\otimes^2}(s_1,s_2)h_2^{k\otimes^2}(s_2,s_3)
1_C(s_1,s_2)1_C(s_2,s_3)d\mu(s_2)  &
\bigl(:=I_1^k(s_1,s_2)\bigl)\\
+&\int h_1^{\otimes^2}(s_1,s_3)h_2^{k\otimes^2}(s_2,s_2)
1_C(s_1,s_2)1_C(s_3,s_2)d\mu(s_2)  &
\bigl(:=I_2^k(s_1,s_2)\bigl)\\
+&\int h_1^{\otimes^2}(s_2,s_2)h_2^{k\otimes^2}(s_1,s_3)
1_C(s_2,s_1)1_C(s_2,s_3)d\mu(s_2)  &
\bigl(:=I_3^k(s_1,s_2)\bigl)\\
+&\int h_1^{\otimes^2}(s_2,s_3)h_2^{k\otimes^2}(s_1,s_2)
1_C(s_2,s_1)1_C(s_3,s_2)d\mu(s_2)  &
\bigl(:=I_4^k(s_1,s_2)\bigl)\\
\end{array}$$
L'inégalité de convexité nous permet d'écrire :
$$S_k\leq
4E\Bigl[\Phi^{\circ^2}\sum_{k}I^k_1\Bigl]^2
+4E\Bigl[\Phi^{\circ^2}\sum_{k}I^k_2-
\int_0^th^2_2\Phi^{\circ^2}\bigl[(h_11_{]0,.]})^{\circ^2}\bigl]d\mu
\Bigl]^2$$
$$+4E\Bigl[\Phi^{\circ^2}\sum_{k}I^k_3\Bigl]^2
+4E\Bigl[\Phi^{\circ^2}\sum_{k}I^k_4\Bigl]^2
$$
\fbox{Etape 1 : } Calcul de
$\mathop{\lim}\limits_{k\to\infty}E\Bigl[\Phi^{\circ^2}
\sum_{k}I^k_1\Bigl]^2$\\
$$\begin{array}{cl}
&E\Phi^{\circ^2}(I_1^k)^2 =\n
I_1^k\n^2_{L^2(\mu^{\otimes^2})}\hspace{2cm} \mbox{(isométrie)}\\
= &\dps \int\Bigl(\sum_{k}
\int h_1^{\otimes^2}(s_1,s_2)h_2^{k\otimes^2}(s_2,s_3)
1_C(s_1,s_2)1_C(s_2,s_3)d\mu(s_2)
\Bigl)^2d\mu^{\otimes^2}(s_1,s_3)\\
= &\dps\int\Bigl[\dps\sum_{k}
\int h_1^{\otimes^4}(s_1,s_2,s_1,t_2)h_2^{k\otimes^4}(s_2,s_3,t_2,s_3)
1_{C^4}(s_1,s_2,s_2,s_3,s_1,t_2,t_2,s_3)d\mu^{\otimes^2}(s_2,t_2)
\\
& +2\dps\sum_{k<j}
\int h_1^{\otimes^4}(s_1,s_2,s_1,t_2)
h_2^{k\otimes^2}(s_2,s_3)h_2^{j\otimes^2}(t_2,s_3)\\
&\hspace{1,5cm}1_{C^4}(s_1,s_2,s_2,s_3,s_1,t_2,t_2,s_3)d\mu^{\otimes^2}(s_2,t_2)
\Bigl]d\mu^{\otimes^2}(s_1,s_3)\\
\end{array}$$
\begin{itemize}
\item
D'une part $\dps\sum_{k}1_{]t_k,t_{k+1}]^4}(s_2,s_3,t_2,s_3)
\mathop{\longrightarrow}\limits_{k\to\infty}^{L^2(\mu^{\otimes^4})}
1_{[s_2=s_3=t_2]}$
\item
D'autre part :
$\dps\sum_{k<j}1_{]t_k,t_{k+1}]^2\times[t_j,t_{j+1}]^2}(s_2,s_3,t_2,s_3)=0$\\
car $(k<j)\Rightarrow(]t_k,t_{k+1}]\bigcap[t_j,t_{j+1}]=\emptyset)$
\end{itemize}
Le théorème de convergence dominée de Lebesgue nous donne :
$$\mathop{\lim}\limits_{k\to\infty}E\Bigl[\Phi^{\circ^2}
\sum_{k}I^k_i\Bigl]^2=\int h_1^{\otimes^2}(s_1,s_2)h_2^{k\otimes^2}(s_2,s_3)
1_{[s_2=s_3=t_2]}
d\mu^{\otimes^4}(s_1,s_2,t_2,s_3)=0$$
\fbox{Etape 2} Calcul de
$\mathop{\lim}\limits_{k\to\infty}E\Bigl[\Phi^{\circ^2}
\dps\sum_{k}I^k_i\Bigl]^2$ pour $i\in\{3,4\}$\\
On procède de la même façon et on a également :
$$\mathop{\lim}\limits_{k\to\infty}E\Bigl[\Phi^{\circ^2}
\sum_{k}I^k_3\Bigl]^2=\mathop{\lim}\limits_{k\to\infty}E\Bigl[\Phi^{\circ^2}
\sum_{k}I^k_4\Bigl]^2=0$$
\fbox{Etape 3 : } Calcul de $\mathop{\lim}\limits_{k\to\infty}
E\Bigl[\Phi^{\circ^2}\dps\sum_{k}I^k_2-
\int_0^th^2_2\Phi^{\circ^2}\bigl(\bigl[h_11_{]0,.]}\bigl]^{\circ^2}\bigl)d\mu
\Bigl]^2$\\
$$\begin{array}{lcl}
\dps\sum_{k}I_2^k(s_2,s_3) & = &
\dps\sum_{k}\int_{t_k}^{t_{k+1}}
h_2^2(s_2)h_1^{\otimes^2}(s_1,s_3)1_C(s_2,s_3)1_{\widetilde C}(s_2,s_3)
d\mu(s_2)\\
& = &\dps\int_0^th_2^2(s_2)h_1^{\otimes^2}(s_1,s_3)
1_C(s_2,s_3)1_{\widetilde C}(s_2,s_3)
d\mu(s_2)\\
& = &\dps\int_0^th_2^2(s_2)
\Bigl(h_11_{]0,s_2]}\otimes h_11_{]0,s_2]}
\Bigl)(s_1,s_3)d\mu(s_2)\\
\end{array}$$
Il s'agit donc de montrer que :
$$\Phi^{\circ^2}\Bigl(
\int_0^th_2^2(s_2)
h_11_{]0,s_2]}^{\otimes^2}d\mu(s_2)\Bigl)
\stackrel{L^2}{=}
\int_0^th_2^2(s_2)
\Phi^{\circ^2}\Bigl(\bigl(
h_11_{]0,s_2]}\bigl)^{\circ^2}\Bigl)d\mu(s_2)
$$
Pour cela on va décomposer $1_C$ en produit tensoriel :
$$1_C (s_1,s_2) \stackrel{L^2}{=}
\mathop{\lim}\limits_{k\to\infty}1_{]0,t_{k-1}]\times]t_k,t_{k+1}]}(s_1,s_2) $$
$$1_{\widetilde C} (s_3,s_2) \stackrel{L^2}{=}
\mathop{\lim}\limits_{k\to\infty}1_{]t_k,t_{k+1}]\times]0,t_{k-1}]}(s_3,s_2) $$
Nous allons noter $C_1^k:=1_{]0,t_{k-1}]}$ et $C_2^k:=1_{]t_k,t_{k+1}]}$
$$\begin{array}{cl}
&\dps\int_0^th_2^2(s_2)h_1^{\otimes^2}(s_1,s_3)
1_C(s_2,s_3)1_{\widetilde C}(s_2,s_3)
d\mu(s_2)\\
= &\dps
h_1^{\otimes^2}\int_0^th_2^2(s_2)
\mathop{\lim}\limits_{{k\to\infty}\atop{l\to\infty}}\dps\sum_k\sum_l
1_{C_1^k\times C_2^k}(.,s_2)1_{C_2^l\times
C_1^l}(s_2,.)d\mu(s_2)\\
& \mathop{=}\limits_{{\mbox{\scriptsize cvg}}\atop{\mbox{dominée}}}
\mathop{\lim}\limits_{{k\to\infty}\atop{l\to\infty}}\dps\sum_{k,l}
h_1^{\otimes^2}1_{C_1^k}1_{C_1^l}
\int_0^th_2^2(s_2)1_{C_2^k}(s_2)1_{C_2^l}(s_2)d\mu(s_2)
\end{array}$$
$$\begin{array}{clc}
&\Phi^{\circ^2}\Bigl(\dps
\int_0^th_2^2(s_2)h_1^{\otimes^2}1_C(.,s_2)1_{\widetilde
C}(s_2,.)d\mu(s_2)\Bigl) \\
=&\Phi^{\circ^2}\left(
\mathop{\lim}\limits_{{k\to\infty}\atop{l\to\infty}}\dps\sum_{k,l}
h_1^{\otimes^2}1_{C_1^k}1_{C_1^l}
\int_0^th_2^2(s_2)1_{C_2^k}(s_2)1_{C_2^l}(s_2)d\mu(s_2)\right) \\
= &\mathop{\lim}\limits_{{k\to\infty}\atop{l\to\infty}}
\Phi^{\circ^2}
\left(
\dps\sum_{k,l}h_1^{\otimes^2}1_{C_1^k}1_{C_1^l}
\int_0^th_2^2(s_2)1_{C_2^k}(s_2)1_{C_2^l}(s_2)d\mu(s_2)
\right) &
        \left\{\begin{array}{l}
        t\mbox{ fixé}\\
        \mbox{cont. de } \Phi^{\circ^2} \mbox{(Prop. 6.3.1.)}
        \end{array}\right.\\
= &\mathop{\lim}\limits_{ {k\to\infty}\atop{l\to\infty} }
\Phi^{\circ^2}
\left(
\dps\sum_{k,l}h_1^{\otimes^2}1_{C_1^k}1_{C_1^l}
\right)
\dps\int_0^th_2^2(s_2)1_{C_2^k}(s_2)1_{C_2^l}(s_2)d\mu(s_2)\\
&\mathop{\lim}\limits_{{k\to\infty}\atop{l\to\infty}}
\dps\int_0^t
\Phi^{\circ^2}
\left(
\dps\sum_{k,l}h_1^{\otimes^2}1_{C_1^k}1_{C_1^l}
\right)
h_2^2(s_2)1_{C_2^k}(s_2)1_{C_2^l}(s_2)d\mu(s_2)\\
\end{array}$$
On peut donc écrire :
$$E\Bigl[\sum_k\Phi^{\circ^2}(I_2^k)-
\int_0^th_2^2(s_2)\Phi^{\circ^2}
\bigl[\bigl(h_11_{]0,s_2]}\bigl)^{\circ^2}\bigl]d\mu(s_2)
\Bigl]^2$$
$$=E\Biggl[\mathop{\lim}\limits_{{k\to\infty}\atop{l\to\infty}}
  \Bigl[\int_0^th_2^2(s_2)
   \Bigl(\Phi^{\circ^2}
   \bigl(
   h_1\otimes h_1 \sum_{(k,l)}
   1_{C_1^k\times C_2^k}(.,s_2)1_{C_2^l\times C_1^l}(s_2,.)
   \bigl)-$$

   $$\Phi^{\circ^2}
   \bigl(
   h_1\otimes h_1
   1_{C}(.,s_2)1_{\widetilde C}(s_2,.)
   \bigl)
   \Bigl)d\mu(s_2)
  \Bigl]
\Biggl]^2$$
Dans la suite on va noter :
\begin{itemize}
\item
$Z_k(.,.,s_2) :=\Phi^{\circ^2}
   \bigl(
   h_1\otimes h_1\dps \sum_{(k,l)}
   1_{C_1^k\times C_2^k}(.,s_2)1_{C_2^l\times C_1^l}(s_2,.)
   \bigl)$
\item
$Z(.,.,s_2):=\Phi^{\circ^2}
   \bigl(
   h_1\otimes h_1
   1_{C}(.,s_2)1_{\widetilde C}(s_2,.)
   \bigl)$
\item
$Y_k := \vert Z_k-Z\vert$
\end{itemize}
$$\begin{array}{cl}
& E\Biggl[
  \Bigl[\dps\int_0^th_2^2(s_2)
   \Bigl(\Phi^{\circ^2}
   \bigl(
   h_1\otimes h_1 \sum_{(k,l)}
   1_{C_1^k\times C_2^k}(.,s_2)1_{C_2^l\times C_1^l}(s_2,.)
   \bigl)-\\
   &\hspace{1cm}\Phi^{\circ^2}
   \bigl(
   h_1\otimes h_1
   1_{C}(.,s_2)1_{\widetilde C}(s_2,.)
   \bigl)
   \Bigl)d\mu(s_2)
  \Bigl]
\Biggl]^2\\
\leq & E\left( \dps\int_0^t h_2^2(s_2)
\left\vert
Z_k(.,.,s_2)-Z(.,.,s_2)
\right\vert d\mu(s_2)\right)^2\\
=&E\left( \dps\int_0^t\int_0^t \bigl(h_2^2\bigl)^{\otimes^2}(s_2,t_2)
Y_k(s_2)Y_k(t_2)d\mu(s_2,t_2)\right)\\
=&\dps\int_0^t\int_0^t \bigl(h_2^2\bigl)^{\otimes^2}(s_2,t_2)
E\left( Y_k(s_2)Y_k(t_2)\right)d\mu(s_2,t_2)\\
\leq &\dps\int_0^t\int_0^t \bigl(h_2^2\bigl)^{\otimes^2}(s_2,t_2)
\left(E \left\vert Y_k(s_2)\right\vert^2\right)^{\frac{1}{2}}
\left(E \left\vert Y_k(t_2)\right\vert^2\right)^{\frac{1}{2}}
d\mu(s_2,t_2)\\
\end{array}$$
Pour finir il suffit de voir que :
$$\begin{array}{lcl}
E \left\vert Y_k(s_2)\right\vert^2 & = &
E\Bigl[
\Phi^{\circ^2}\bigl(h_1\otimes h_1\dps\sum_{(k,l)}
1_{C_1^k\times C_2^k}(.,s_2)1_{C_2^l\times C_1^l}(s_2,.)-\\
&&\hspace{1cm}h_1\otimes h_1\dps\sum_{(k,l)}
1_{C}(.,s_2)1_{\widetilde C}(s_2,.)
\bigl)\Bigl]^2\\
 & = &\left\Vert h_1\otimes h_1
 \Bigl[
 1_{C}(.,s_2)1_{\widetilde C}(s_2,.)-
\dps\sum_{(k,l)}
1_{C_1^k\times C_2^k}(.,s_2)1_{C_2^l\times C_1^l}(s_2,.)
\Bigl]
 \right\Vert^2
\end{array}$$
\end{enumerate} \normalsize$\triangle$
%
%
%
%
\newpage
\chapter{Existence de modifications continues}
\section{Méthode}
\subsection{Un résultat technique}
On va utiliser la variante, un peu plus générale, du théorème de
Kolmogorov-Centsov qui suit (confère {\bf [12]}) :
\begin{lm}
On suppose qu'un processus $X:=\bigl(X_t\bigl)_{0\leq t\leq T}$
sur un espace probabilisé $\bigl(\Omega,\F,P\bigl)$ satisfasse la
condition :
$$E\left\vert X_t-X_s\right\vert^{\alpha}\leq C
\left\vert F(t)-F(s)\right\vert^{1+\beta}$$
Où :
\begin{itemize}
\item
$\alpha$, $\beta$ et C sont des constantes strictement
positives.
\item
F est une fonction réelle, continue sur $[0,T]$ et à variation
bornée.
\end{itemize}
Alors il existe des modifications continues de X
\end{lm}
\subsection*{\underline{\textbf{Démonstration}}}
\begin{enumerate}
\item[]
\small
Pour l'essentiel on reprend la démonstration de {\bf [12]}. On a
deux points à vérifier :
\begin{enumerate}
\item
$P\bigl(\bigl[\vert X_t-X_s\vert>\varepsilon\bigl]\bigl)
\leq
\dps\frac{E\vert X_t-X_s\vert^{\alpha}}{\varepsilon^{\alpha}}
\leq
C\varepsilon^{-\alpha}\vert F(t)-F(s)\vert^{1+\beta}$\\
et la continuité de F entraîne :
$X_s\mathop{\longrightarrow}\limits_{s\to t}X_t$ en probabilité.
\item
On prend $t=\frac{k}{2^n}$, $s=\frac{k-1}{2^n}$ et
$\varepsilon=2^{-\gamma n}$ $(0<\gamma<\frac{\beta}{\alpha})$\\
$P\bigl(\bigl[\vert X_{\frac{k}{2^n}}-X_{\frac{k-1}{2^n}}\vert
\geq
2^{-\gamma n}\vert\bigl]\bigl)
\leq
C2^{-\gamma n}\vert F(\frac{k}{2^n})-F(\frac{k-1}{2^n})\vert$\\
$P\Biggl(\Bigl[\dps\max_{0\leq k\leq2^n}
\vert X_{\frac{k}{2^n}}-X_{\frac{k-1}{2^n}}\vert\Bigl]\bigl)
\geq
2^{-\gamma n}\Biggl)
\leq
P\left( \dps\bigcup_{k=0}^{2^n}
\Bigl[\vert X_{\frac{k}{2^n}}-X_{\frac{k-1}{2^n}}\vert\geq 2^{-\gamma n}
\Bigl]\right)$\\
$\leq C2^{-\gamma n}\dps\sum_{k=0}^{2^n}
\left\vert F\left(\frac{k}{2^{n}}\right)
-F\left(\frac{k-1}{2^n}\right)\right\vert
\leq C\frac{varF}{2^{\gamma n}}$\\
où l'on a noté varF pour la variation de F.
\end{enumerate}
Tous les autres arguments sont repris sans changement.
\end{enumerate} \normalsize$\triangle$
\subsection{Position du problème}
\begin{itemize}
\item
On choisit $h_1, h_2\in H$, $t\in\R$
\item
On pose $Z_t^n:=\left(\Phi^{\circ^2}+\Phi\circ a_1^2\right)_n
\left(h_1\otimes\Bigl(h_21_{]0,t]}\Bigl)1_C\right)$
\item
On sait que la suite $\bigl(Z_t^n\bigl)$ converge dans $L^2$ vers
 $\dps\int_0^t\Phi(h_1)_sd\Phi(h_2)_s=:Z_t$ et on veut montrer
 l'existence de modifications continues pour le processus
 $\bigl(Z_t\bigl)$.\\
 Pour ce faire, on va étudier $E\vert Z_t^n-Z_s^n\vert^4$ dans le
 but d'utiliser la variante précédente du théorème de
 Kolmogorov-Centsov.\\
 On peut remarquer que
$Z_t^n-Z_s^n=\left(\Phi^{\circ^2}+\Phi\circ a_1^2\right)_n
\left(h_1\otimes\Bigl(h_21_{]s,t]}\Bigl)1_C\right)$.\\
Dans le but de simplifier les écritures nous noterons
$f_{st}:=h_1\otimes\Bigl(h_21_{]s,t]}\Bigl)1_C$.
\end{itemize}
\subsection{Méthode}
Pour calculer $E\vert Z_t^n-Z_s^n\vert^4$ nous allons utiliser les
calculs de $\left(\Phi^{\circ^2}+\Phi\circ a_1^2\right)^2$ et
l'orthogonalité des termes d'ordre différent. Les
identifications établies précédemment nous permettrons de faire
cette étude sur chacun des ordres, ordre après ordre. A cet effet
on rappel :
$$\left\vert\begin{array}{lcl}\hline
O_4\bigl(\vert Z_t^n-Z_s^n\vert^2\bigl) & = &
\Phi^{\circ^4}_n\bigl(f_{st}^{\circ}\circ f_{st}^{\circ}\bigl)\\
O_3\bigl(\vert Z_t^n-Z_s^n\vert^2\bigl) & = &
\left(\Phi^{\circ^3}_n\circ a_1^2\right)_n
\bigl(f_{st}^{\circ}\circ f_{st}^{\circ}\bigl)\\
O_2\bigl(\vert Z_t^n-Z_s^n\vert^2\bigl) & = &
\Phi^{\circ^2}_n\left(
4f_{st}^{\circ}{\mathop{\sim}\limits_{1}}f_{st}^{\circ}+
a^4_1\bigl(f_{st}^{\circ}\circ f_{st}^{\circ}\bigl)
\right)\\
O_1\bigl(\vert Z_t^n-Z_s^n\vert^2\bigl) & = &
\Phi_n\left(
a^4_3\bigl(f_{st}^{\circ}\circ f_{st}^{\circ}\bigl)+
4a^4_1\bigl(f_{st}^{\circ}{\mathop{\sim}\limits_{1}}f_{st}^{\circ}\bigl)
-6a^4_1\circ\pi_1\bigl(f_{st}^{\circ}{\mathop{\sim}\limits_{1}}f_{st}^{\circ}\bigl)
\right)\\
O_0\bigl(\vert Z_t^n-Z_s^n\vert^2\bigl) & = &
2\n f_{st}\n^2_A+
a^4_4\bigl(f_{st}^{\circ}\circ f_{st}^{\circ}\bigl)\\ \hline
\end{array}\right\vert$$
\section{Termes d'ordre 4}
\begin{lm}
$$E\Bigl(O_4\bigl(\vert Z_t^n-Z_s^n\vert^2\bigl)\Bigl)^2
\leq
\n h_1\n_A^4\n h_21_{]s,t]}\n^4_A$$
\end{lm}
\subsection*{\underline{\textbf{Démonstration}}}
\begin{enumerate}
\item[]
\small
$E^{\nu}\Bigl(O_4\bigl(\vert Z_t^n-Z_s^n\vert^2\bigl)\Bigl)^2
\leq E^{\mu_A}\bigl(f_{st}^{\circ}\circ f_{st}^{\circ}\bigl)^2$\\
mais $$\begin{array}{lcl}
f_{st}^{\circ}\circ f_{st}^{\circ} & = &
f_{st}^{\circ}\otimes f_{st}^{\circ}\\
& = &\frac{1}{4}\left(
h_1\otimes\bigl(h_21_{]s,t]}\bigl)1_C+
\bigl(h_21_{]s,t]}\otimes h_1\bigl)1_{\widetilde C}
\right)^{\otimes^2}\\
& = &\frac{1}{4}\Bigl(
h_1\otimes h_2^{st}\otimes
h_1\otimes h_2^{st}    1_C\otimes1_C+
h_1\otimes h_2^{st}\otimes
h_2^{st}\otimes h_1    1_C\otimes1_{\widetilde C}\\
&&\hspace{1cm}+
h_2^{st}\otimes h_1\otimes
h_1\otimes h_2^{st}1_{\widetilde C}\otimes1_C+
h_2^{st}\otimes h_1\otimes
h_2^{st}\otimes h_11_{\widetilde C}\otimes1_{\widetilde C}\Bigl)\\
\end{array}$$
$$\begin{array}{lcl}
E^{\nu}\left(O_4\bigl(\vert Z_t^n-Z_s^n\vert^2\bigl)\right)^2
&\leq&
\frac{1}{4}\Bigl[
\n h_1\otimes h_2^{st}\otimes
h_1\otimes h_2^{st}    1_C\otimes1_C\n_A^2\\
&&\hspace{0,5cm}+\n h_1\otimes h_2^{st}\otimes
h_2^{st}\otimes h_1    1_C\otimes1_{\widetilde C}\n_A^2\\
&&\hspace{0,5cm}+\n h_2^{st}\otimes h_1\otimes
h_1\otimes h_2^{st}1_{\widetilde C}\otimes1_C\n_A^2\\
&&\hspace{0,5cm}+\n h_2^{st}\otimes h_1\otimes
h_2^{st}\otimes h_11_{\widetilde C}\otimes1_{\widetilde C}\Bigl)\n_A^2\Bigl]
\end{array}$$
$$\begin{array}{lcl}
\n h_1\otimes h_2^{st}\otimes
h_1\otimes h_2^{st}    1_C\otimes1_C\n_A^2 & \leq&
\dps \int\bigl(h_1\otimes h_1\bigl)^2(x_1,x_3)
\bigl(h_2^{st}\otimes
h_2^{st}\bigl)^2(x_2,x_4)d\mu_A^{\otimes^4}(x_1,x_2,x_3,x_4)\\
& = &\n h_1 \n^4_A\n h_2^{st}\n_A^4
\end{array}$$
Les trois autres termes se traitent de la même façon et donnent le
même résultat, ce qui établit le lemme.
\end{enumerate} \normalsize$\triangle$
\section{Termes d'ordre 3}
\begin{lm}
$$E\Bigl(O_3\bigl(\vert Z_t^n-Z_s^n\vert^2\bigl)\Bigl)^2
\leq C_{4,1}
\n h_1\n_A^4\n h_21_{]s,t]}\n^4_A$$
\end{lm}
\subsection*{\underline{\textbf{Démonstration}}}
\begin{enumerate}
\item[]
\small
$$\begin{array}{lcl}
E^{\nu}\Bigl(O_3\bigl(\vert Z_t^n-Z_s^n\vert^2\bigl)\Bigl)^2 &
\leq & C_{4,1}  \n a_1^4\bigl(f_{st}^{\circ}\circ
f_{st}^{\circ}\bigl)\n^2_{A}\\
& \leq & C_{4,1}\n f_{st}^{\circ}\circ f_{st}^{\circ}\n^2_A\\
& \leq & C_{4,1}\n h_1\n^4_A\n h_21_{]s,t]} \n^4_A\\
\end{array}$$
\end{enumerate} \normalsize$\triangle$
\section{Termes d'ordre 2}
\begin{lm}
$$E\Bigl(O_2\bigl(\vert Z_t^n-Z_s^n\vert^2\bigl)\Bigl)^2
\leq\left(\frac{1}{2}+C_{4,2}\right)
\n h_1\n_A^4\n h_21_{]s,t]}\n^4_A$$
\end{lm}
\subsection*{\underline{\textbf{Démonstration}}}
\begin{enumerate}
\item[]
\small
$$E\Bigl(O_2\bigl(\vert Z_t^n-Z_s^n\vert^2\bigl)\Bigl)^2
\leq
4\left\Vert
f_{st}^{\circ}\mathop{\sim}\limits_{1}f_{st}^{\circ}
\right\Vert+
\left\Vert
a_2^4\bigl(f_{st}^{\circ}\circ f_{st}^{\circ} \bigl)
\right\Vert
$$
$$\begin{array}{lcl}
\left\Vert
f_{st}^{\circ}\mathop{\sim}\limits_{1}f_{st}^{\circ}
\right\Vert & = &
\dps\int\left[\int
 f_{st}^{\circ}(s_1,s_2)f_{st}^{\circ}(s_2,s_3)d\mu_A(s_2)\right]^2
 d\mu_A^{\otimes^2}(s_2,s_3)\\
 & \mathop{\leq}\limits_{\mbox{Hölder}} &\dps\int
 \left(\int\Bigl(
 f_{st}^{\circ}(s_1,s_2)\Bigl)^2d\mu_A(s_2)\right)
\left(\int\Bigl( f_{st}^{\circ}(s_2,s_3)\Bigl)^2d\mu_A(s_2)
d\mu_A^{\otimes^2}(s_2,s_3)\right)\\
& \mathop{\leq}\limits_{\mbox{Fubini}} &\n f_{st}^{\circ}\n^4_A\\
\n f_{st}^{\circ}\n^4_A & = &\dps\int
\left[
\Bigl(
h_1\otimes h_21_{]s,t]}
\Bigl)^{\circ^2}
\right]^2d\mu_A^{\otimes^2}\\
&\leq & \frac{1}{2}\n h_1\n^4_A\n h_21{]s,t]}\n^4_A
\end{array}$$
$\begin{array}{lcl}
\n a_2^4 \bigl( f_{st}^{\circ}\circ f_{st}^{\circ}\bigl)\n_A^2 &
\leq&
C_{4,2}\n h_1\n_A^4\n h_21_{]s,t]}\n_A^4
\end{array}$
d'où le résultat.
\end{enumerate} \normalsize$\triangle$
\section{Termes d'ordre 1}
\begin{lm}
$$E\Bigl(O_1\bigl(\vert Z_t^n-Z_s^n\vert^2\bigl)\Bigl)^2
\leq\left(C_{4,3}+\frac{5}{2}C_{4,1}\right)
\n h_1\n_A^4\n h_21_{]s,t]}\n^4_A$$
\end{lm}
\subsection*{\underline{\textbf{Démonstration}}}
\begin{enumerate}
\item[]
\small
$$\begin{array}{lcl}
\left\Vert
a_3^4 \bigl( f_{st}^{\circ}\circ f_{st}^{\circ}\bigl)
\right\Vert_A^2 &
\leq&
C_{4,3}\n h_1\n_A^4\n h_21_{]s,t]}\n_A^4\\
\left\Vert
a_1^4 \bigl( f_{st}^{\circ}\circ f_{st}^{\circ}\bigl)
\right\Vert_A^2 & \leq &
C_{4,1}\n
f_{st}^{\circ}\mathop{\sim}\limits_{1} f_{st}^{\circ}\n_A^2
\leq
\frac{1}{2}C_{4,1}\n h_1\n_A^4\n h_21_{]s,t]}\n_A^4\\
\left\Vert
a_1^4 \circ\pi_1\bigl(
 f_{st}^{\circ} \mathop{\sim}\limits_{1}f_{st}^{\circ}\bigl)
\right\Vert_A^2 &
\leq&
C_{4,1}\left\Vert
\pi_1\bigl(f_{st}^{\circ}
\mathop{\sim}\limits_{1}f_{st}^{\circ}\bigl)
\right\Vert_A^4
\leq
C_{4,1}\left\Vert
f_{st}^{\circ} \mathop{\sim}\limits_{1}f_{st}^{\circ}
\right\Vert_A^4\\
&\leq & \n h_1\n_A^4\n h_21_{]s,t]}\n_A^4\\
\end{array}$$
\end{enumerate} \normalsize$\triangle$
\section{Termes d'ordre 0}
\begin{lm}
$$E\Bigl(O_0\bigl(\vert Z_t^n-Z_s^n\vert^2\bigl)\Bigl)^2
\leq\left(C_{4,4}+\frac{1}{2}\right)
\n h_1\n_A^4\n h_21_{]s,t]}\n^4_A$$
\end{lm}
\subsection*{\underline{\textbf{Démonstration}}}
\begin{enumerate}
\item[]
\small
$$\left\Vert
a_4^4 \bigl(
 f_{st}^{\circ} \circ f_{st}^{\circ}\bigl)
\right\Vert_A^2
\leq C_{4,4}
\left\Vert
f_{st}^{\circ} \circ f_{st}^{\circ}
\right\Vert_A^2
\leq C_{4,4}\n h_1\n_A^4\n h_21_{]s,t]}\n^4_A$$
$$\left\Vert
f_{st}^{\circ}
\right\Vert_A^2 = \dps\int
\Bigl[
\Bigl(
h_1\otimes \bigl(h_21_{]s,t]}\bigl)1_C
\Bigl)^{\circ^2}
\Bigl]^2d\mu^{\otimes^2}_A
\leq
\frac{1}{2}\n h_1\n_A^4\n h_21_{]s,t]}\n^4_A$$
\end{enumerate} \normalsize$\triangle$
\section{Modifications continues}
La réunion des résultats précédents donne le lemme :
\begin{lm}
$$E\left(\vert Z_t^n-Z_s^n\vert^2\right)^2
\leq(\frac{7}{2}C_{4,1}+C_{4,2}+C_{4,3}+C_{4,4}+2)
\n h_1\n_A^4\n h_21_{]s,t]}\n^4_A$$
\end{lm}
\begin{thm}
Le processus $\left(\dps\int_0^t\Phi(h_1)_sd\Phi(h_2)_s\right)_{t}$
admet des modifications continues.
\end{thm}
\subsection*{\underline{\textbf{Démonstration}}}
\begin{enumerate}
\item[]
\small
$\n h_21_{]s,t]}\n^2_A = \dps\int h_2^21_{]s,t]}d\mu_A$\\
On applique le lemme 8.1.1. avec $F(t)=\dps\int h_2^21_{]s,t]}d\mu_A$
\end{enumerate} \normalsize$\triangle$
%
%
%
%
\newpage
\chapter{Variation quadratique}
\section{Position du problème}
On veut mettre en évidence la variation quadratique du processus :
$$\left(\dps\int_0^t\Phi(h_1)_sd\Phi(h_2)_s\right)_{t\geq0}.$$
Nous conservons la notation $Z_t:=\int_0^t\Phi(h_1)_sd\Phi(h_2)_s$.\\
Plus précisément, étant donnée une suite de partitions
$0=t_0<t_1<\cdots<t_n=t$ de [0,t], nous cherchons la limite dans $L^2$
de la suite :
$$\left(\sum_{k=0}^n
\left\vert Z_{t_{k+1}}-Z_{t_k}
\right\vert^2
\right)_{n\in{\tiny\N}}$$
Nous avons déjà établi que (en notant $f_k$ pour $f_{t_kt_{k+1}}$)
$$\begin{array}{lcl}
\left\vert Z_{t_{k+1}}-Z_{t_k}
\right\vert^2 & = &
\Phi^{\circ^4}_n\bigl(f_{k}^{\circ}\circ f_{k}^{\circ}\bigl)\\
 & + &
\left(\Phi^{\circ^3}\circ a_1^2\right)_n
\bigl(f_{k}^{\circ}\circ f_{k}^{\circ}\bigl)\\
&+ &
\Phi^{\circ^2}_n\left(
4f_{k}^{\circ}{\mathop{\sim}\limits_{1}}f_{k}^{\circ}+
a^4_1\bigl(f_{k}^{\circ}\circ f_{k}^{\circ}\bigl)
\right)\\
&+&
\Phi_n\left(
a^4_3\bigl(f_{k}^{\circ}\circ f_{kt}^{\circ}\bigl)+
4a^4_1\bigl(f_{k}^{\circ}{\mathop{\sim}\limits_{1}}f_{k}^{\circ}\bigl)
-6a^4_1\circ\pi_1\bigl(f_{k}^{\circ}{\mathop{\sim}\limits_{1}}f_{k}^{\circ}\bigl)
\right)\\
&+ &
2\n f_{st}\n^2_A+
a^4_4\bigl(f_{k}^{\circ}\circ f_{k}^{\circ}\bigl)\\
\end{array}$$
La linéarité et la continuité des opérateurs  nous
incite, avant d'identifier la variation quadratique,
 à étudier les deux limites de $L^2$ suivantes :
$$\left\{\begin{array}{l}
\dps\lim_{k\to\infty}\sum_{k}
f_{k}^{\circ}\circ f_{k}^{\circ}\\
\dps\lim_{k\to\infty}

4a^4_1\left(
\sum_{k}\Bigl(f_{k}^{\circ}{\mathop{\sim}\limits_{1}}f_{k}^{\circ}
-6\pi_1\bigl(f_{k}^{\circ}{\mathop{\sim}\limits_{1}}f_{k}^{\circ}\bigl)
\Bigl)
\right)\\
\end{array}\right.$$
Pour ce qui concerne la deuxième on se reportera au lemme 7.3.1.
\section{Lemmes techniques}
\begin{lm}
$$\lim_{k\to\infty}^{L^2}\sum_{k}
f_{k}^{\circ}\circ f_{k}=0$$
\end{lm}
\subsection*{\underline{\textbf{Démonstration}}}
\begin{enumerate}
\item[]
\small
On a déjà vu que :\\
$f_k^{\circ}\circ f_k^{\circ}= \frac{1}{4}\Bigl(
h_1\otimes h^k_2\otimes h_1\otimes h^k_2
1_C\otimes1_C +
h_1\otimes h^k_2\otimes h^k_2\otimes h_1
1_C\otimes1_{\widetilde C}+$\\
$\hspace*{2,5cm}
h^k_2\otimes h_1\otimes h_1\otimes h^k_2
1_{\widetilde C}\otimes1_C +
h^k_2\otimes h_1\otimes h^k_2\otimes h_1
1_{\widetilde C}\otimes1_{\widetilde C}
\Bigl)$\\
$$\begin{array}{lcl}
\left\Vert
\dps\sum_k f_k^{\circ}\circ f_k^{\circ}
\right\Vert^2 & \leq &
\dps\int
\left(
 \sum_k\sum_{\sigma\in\Sigma_4}
 \Bigl(
  h_1\otimes h^k_2\otimes h_1\otimes h^k_2
  1_C\otimes1_C
 \Bigl)_{\sigma}
\right)^2d\mu^{\otimes^4}\\
&\leq &
(4!)^2\dps\sum_{\sigma\in\Sigma_4}
\Biggl[
 \sum_k\int
 \Bigl(
  h_1\otimes h^k_2\otimes h_1\otimes h^k_2
  1_C\otimes1_C
 \Bigl)_{\sigma}^2d\mu^{\otimes^4}\\
 &&+2\dps\sum_{k<j}\int
 \Bigl(
  h_1\otimes h^k_2\otimes h_1\otimes h^k_2
  1_C\otimes1_C
 \Bigl)_{\sigma}
 \Bigl(
  h_1\otimes h^j_2\otimes h_1\otimes h^j_2
  1_C\otimes1_C
 \Bigl)_{\sigma}d\mu^{\otimes^4}
\Biggl]\\

&&\hspace{3cm}\mbox{(inégalité de convexité)}\\
& \leq &
(4!)^2\dps\sum_{\sigma\in\Sigma_4}
\Biggl[
 \sum_k\int
 \Bigl[
  h_1^{\otimes^4}
 \Bigl]_{\sigma}d\mu^{\otimes^4}
 \int
 \Bigl[
  h_2^{\otimes^4}\Bigl(1_{]t_k,t_{k+1}]}\Bigl)^{\otimes^4}
 \Bigl]_{\sigma}d\mu^{\otimes^4}\\
&&+2\dps\sum_{k<j}\int
h_1^{\otimes^4}(x_{\sigma(1)},x_{\sigma(3)},y_{\sigma(1)},x_{\sigma(3)})
h_2^{\otimes^4}(x_{\sigma(2)},x_{\sigma(4)},y_{\sigma(2)},x_{\sigma(4)})\\
&&\hspace*{2cm}1_{C^{\times^4}}(x_{\sigma},y_{\sigma})
\Bigl(
 1_{]t_k,t_{k+1}]\times]t_j,t_{j+1}]}
 (x_{\sigma(2)},x_{\sigma(4)},y_{\sigma(2)},x_{\sigma(4)})
\Bigl)^{\otimes^2}\Biggl]\\
&&\hspace*{4cm}d\mu^{\otimes^8}(x_{\sigma},y_{\sigma})\\
\end{array}$$
On sait que
$$\begin{array}{lcl}
\dps\sum_k\Bigl(1_{]t_k,t_{k+1}]} \Bigl)^{\otimes^4}&
\mathop{\longrightarrow}\limits^{L^2\bigl(\mu^{\otimes^4}\bigl)}&
1_{[x_{\sigma(2)}=x_{\sigma(4)}=y_{\sigma(2)}=y_{\sigma(4)}]}\\
\dps\sum_{k<j}
\Bigl(
 1_{]t_k,t_{k+1}]\times]t_j,t_{j+1}]}
\Bigl)^{\otimes^2}&
\mathop{\longrightarrow}\limits^{L^2\bigl(\mu^{\otimes^4}\bigl)}&
1_{[x_{\sigma(2)}<x_{\sigma(4)},x_{\sigma(2)}=y_{\sigma(2)},
x_{\sigma(4)}=y_{\sigma(4)}]}\\
\end{array}$$
Pour établir le dernier point il suffit de remarquer :\\
\begin{itemize}
\item d'une part que
$\left\{\begin{array}{lcl}
1_{[x=y]} & = &
\dps \lim_{k\to\infty}^{L^2}
\sum_k1_{]t_k,t_{k+1}]\times]t_k,t_{k+1}]}(x,y)\\
1_{[x<y]} & = &
\dps \lim_{j\to\infty}^{L^2}
\sum_{k<j}1_{]t_k,t_{k+1}]\times]t_j,t_{j+1}]}(x,y)\\
\end{array}\right.$
\item et d'autre part que
$\left\{\begin{array}{l}
1_{[x=y]} 1_{[z=t]}1_{[x<z]}1_{[y<t]} \\
= \dps \lim_{k,j\to\infty}^{L^2}\Biggl[
\sum_k1_{]t_k,t_{k+1}]^2}(x,y)\cdot
\sum_k1_{]t_k,t_{k+1}]^2}(z,t)\cdot\\
\dps\sum_{k<j}1_{]t_k,t_{k+1}]\times]t_j,t_{j+1}]}(x,z)\cdot
\sum_{k<j}1_{]t_k,t_{k+1}]\times]t_j,t_{j+1}]}(y,t)\Biggl]
\\
\end{array}\right.$
\item car
$\left\{\begin{array}{l}
\Biggl[\dps\sum_k1{]t_k,t_{k+1}]^2}(x,y)\cdot
\sum_k1{]t_k,t_{k+1}]^2}(z,t)\cdot\\
\dps\sum_{k<j}1{]t_k,t_{k+1}]\times]t_j,t_{j+1}]}(x,z)\cdot
\sum_{k<j}1{]t_k,t_{k+1}]\times]t_j,t_{j+1}]}(y,t)\Biggl]\\
=\dps\sum_{k<j}1{]t_k,t_{k+1}]^2}\otimes1_{]t_j,t_{j+1}]^2}(x,y,z,t)\\
\end{array}\right.$
\end{itemize}
Finalement, en utilisant le théorème de convergence dominée de
Lebesgue et la mesure nulle de ces ensembles on a le résultat
annoncé.
\end{enumerate} \normalsize$\triangle$
\begin{lm}
$$\Biggl(\Phi\Bigl(h_11_{]0,s]}\Bigl)\Biggl)^2=
\Phi^{\circ^2}\Bigl(\Bigl[h_11_{]0,s]}\Bigl]^{\circ^2}\Bigl)+
\Phi\circ a_1^2\Bigl(\Bigl[h_11_{]0,s]}\Bigl]^{\circ^2}\Bigl)+
\left\Vert
h_11_{]0,s]}
\right\Vert^2$$
\end{lm}
\subsection*{\underline{\textbf{Démonstration}}}
\begin{enumerate}
\item[]
\small
On utilise la proposition 6.2.1. :
$$\Phi\Bigl(h_11_{]0,s]}\Bigl)^2 =
\varphi^{(2)}\Bigl( h_11_{]0,s]}\otimes h_11_{]0,s]}\Bigl)+
\varphi^{(1,1)}\Bigl(2 h_11_{]0,s]}\otimes h_11_{]0,s]}\Bigl)+
\left\Vert h_11_{]0,s]} \right\Vert^2$$
La proposition 6.3.2. nous permet alors de conclure
\end{enumerate} \normalsize$\triangle$
\section{Variation Quadratique}
Les lemmes et calculs précédents donnent le résultat final,
 à savoir :
\begin{thm}
$$\lim_{k\to\infty}^{L^2}\sum_k\left\vert Z_{t_{k+1}}-Z_{t_{k}}\right\vert^2
= \int_0^th_2^2(s)\Phi\bigl(h_11_{[0,s]}\bigl)^2d\mu(s)+
\int_0^th_2^2(s)
\Phi\circ a_1^2\circ \pi_1
\bigl(  h_11_{]0,s]}  \bigl)^{\circ^2}d\mu(s)$$
\end{thm}
\subsection*{Remarques}
\begin{itemize}
\item
On peut écrire $\dps\int_0^th_2^2(s)\Phi\bigl(h_11_{[0,s]}\bigl)^2d\mu(s)=
\dps\int_0^th_2^2(s)\Biggl(
\int_0^sh_1(u)d\mu(u)\Biggl)^2d\mu(s)$
\item
Dans le cas où les variables $\bigl(X_k\bigl)_k$ suivent une loi normale
on retrouve la variation quadratique gaussienne puisqu'alors
 $a_1^2\circ\pi_1\bigl(h_11_{]0,s]}\bigl)^{\circ^2}=0$
\end{itemize}
%
%
\newpage
\begin{center}
{\LARGE BIBLIOGRAPHIE}
\end{center}
\vspace{3cm}
\Large
{\bf [1]} SIMON : \textit{The ${\cal P}\bigl(\phi\bigl)_2$ Euclidean (Quantum) Field Theorie  }
\newline
{\bf [2]}   SZEGO :   \textit{Orthogonal polynomial}, 1939
\newline
{\bf [3]}   WILF,Herbert S. :   \textit{GeneratingFunctionology},
Academic Press,London, 1990.
\newline
{\bf [4]} BOCCARA   :   \textit{Analyse fonctionnelle }, Masson, Paris, 19xx
\newline
{\bf [5]} BREZIS, Haïm    :   \textit{Analyse fonctionnelle et
applications}, Masson, Paris, 1983.
\newline
{\bf [6]} DACUNHA-CASTEL\&DUFFLO   :   \textit{Probabilités et
statistiques}, Ellipse, Paris, 1983.
\newline
{\bf [7]}    SCHRAFTETTER, Eric :   \textit{Quelques aspects des séries aléatoires
définies par une mesure vectorielle, application aux EDPS}, Thèse
de doctorat, Angers, 1998.
\newline
{\bf [8]} SCHWARTZ, Laurent   :   \textit{Théorie des
distributions}, Herman, Paris, 1966.
\newline
{\bf [9]} STROMBERG, K. R.   :   \textit{Probability for
Analysts}, Chapman\&Hall, New-York, 1994
\newline
{\bf [10]} MEYER, Paul André    :   \textit{Quantum Probability for
Probabilists}, Springer, LNM 1533, 1995.
\newline
{\bf [11]}  I. KARATZAS and S.E. SHREVE  :   \textit{Brownian Motion and Stochastic Calculus}, Springer-Verlag , Volume 113 in the series "Graduate Texts in Mathematics", 1988.
\end{document}